\pdfoutput=1

\documentclass[12pt,a4paper]{article}

\usepackage{amsfonts,amssymb,amsthm,amsmath,latexsym}
\usepackage{eqsection,indent}
\usepackage{subeqnarray, eepicemu, url,cite,bm}
\usepackage{multirow}  
\usepackage{tensor}  
\usepackage{graphicx,graphbox}
\usepackage{float}  

\usepackage{tikz}

\usepackage{color}

\usetikzlibrary{arrows,shapes,decorations,automata,backgrounds,petri,bending}
\usetikzlibrary{matrix,cd,decorations.markings,graphs}
\usetikzlibrary{positioning}

\date{November 28, 2023}

\DeclareFontFamily{U}{mathx}{}
\DeclareFontShape{U}{mathx}{m}{n}{<-> mathx10}{}
\DeclareSymbolFont{mathx}{U}{mathx}{m}{n}
\DeclareMathAccent{\widecheck}{0}{mathx}{"71}

\tikzset{->-/.style={decoration={markings,mark=at position #1 with {\arrow{>}}},postaction={decorate}}}

\begin{document}

\title{\vspace*{-1cm}
       Lattice paths and branched continued fractions  \\[3mm]
       III.~Generalizations of the \\ Laguerre, rook and Lah polynomials
      }

\author{
      \hspace*{-1.8cm}
      {\large Bishal Deb${}^{1,2}$, Alexander Dyachenko${}^3$,
                 Mathias P\'etr\'eolle${}^1$ and Alan D.~Sokal${}^{1,4}$}
   \\[5mm]
     \hspace*{-1.3cm}
      \normalsize
           ${}^1$Department of Mathematics, University College London,
                    London WC1E 6BT, UK   \\
     \hspace*{-1.15cm}
      \normalsize
           ${}^2$Sorbonne Universit\'e and Universit\'e Paris Cit\'e, CNRS,
         Laboratoire de Probabilit\'es, \\[-1mm]
     \hspace*{-6cm}
      \normalsize
     Statistique et Mod\'elisation, 75005 Paris, France \\
     \hspace*{-2cm}
      \normalsize
           ${}^3$Keldysh Institute of Applied Mathematics,
                    Russian Academy of Sciences, Moscow 125047, Russia  \\
     \hspace*{-2.95cm}
      \normalsize
           ${}^4$Department of Physics, New York University,
                    New York, NY 10003, USA
       \\[4mm]
     \hspace*{-1cm}
     {\tt bishal@gonitsora.com},
     {\tt diachenko@sfedu.ru}, \\
     \hspace*{-1cm}
     {\tt mathias.petreolle@gmail.com}, {\tt sokal@nyu.edu}
}

\maketitle
\thispagestyle{empty}   

\begin{abstract}
We introduce a triangular array $\widehat{\sf L}^{(\alpha)}$
of 5-variable homogeneous polynomials that enumerate Laguerre digraphs
(digraphs in which each vertex has out-degree 0~or~1 and in-degree 0~or~1)
with separate weights for peaks, valleys, double ascents, double descents,
and loops.
These polynomials generalize the classical Laguerre polynomials
as well as the rook and Lah polynomials.
We show that this triangular array is totally positive
and that the sequence of its row-generating polynomials
is Hankel-totally positive,
under suitable restrictions on the values given to the indeterminates.
This implies, in particular, the coefficientwise Hankel-total positivity
of the monic unsigned univariate Laguerre polyomials.
Our proof uses the method of production matrices
as applied to exponential Riordan arrays.
Our main technical lemma concerns the total positivity
of a large class of quadridiagonal production matrices;
it generalizes the tridiagonal comparison theorem.
In some cases these polynomials are given by a branched continued fraction.
Our constructions are motivated in part by recurrences for
the multiple orthogonal polynomials associated to weights
based on modified Bessel functions of the first kind $I_\alpha$.
\end{abstract}


\medskip
\noindent
{\bf Mathematics Subject Classification (MSC 2020) codes:}
05A15 (Primary);
05A05, 05A19, 05A20, 05C30, 05C38, 15B48, 30B70, 33C10, 33C45, 42C05 (Secondary).



\newtheorem{theorem}{Theorem}[section]
\newtheorem{proposition}[theorem]{Proposition}
\newtheorem{lemma}[theorem]{Lemma}
\newtheorem{corollary}[theorem]{Corollary}
\newtheorem{definition}[theorem]{Definition}
\newtheorem{conjecture}[theorem]{Conjecture}
\newtheorem{question}[theorem]{Question}
\newtheorem{problem}[theorem]{Problem}
\newtheorem{openproblem}[theorem]{Open Problem}
\newtheorem{example}[theorem]{Example}
\newtheorem{remark}[theorem]{Remark}

\renewcommand{\theenumi}{\alph{enumi}}
\renewcommand{\labelenumi}{(\theenumi)}
\def\eop{\hbox{\kern1pt\vrule height6pt width4pt
depth1pt\kern1pt}\medskip}
\def\prf{\par\noindent{\bf Proof.\enspace}\rm}
\def\rmk{\par\medskip\noindent{\bf Remark\enspace}\rm}

\newcommand{\textbfit}[1]{\textbf{\textit{#1}}}

\newcommand{\bigdash}{%
\smallskip\begin{center} \rule{5cm}{0.1mm} \end{center}\smallskip}

\newcommand{\safepar}{ {\protect\hfill\protect\break\hspace*{5mm}} }

\newcommand{\be}{\begin{equation}}
\newcommand{\ee}{\end{equation}}
\newcommand{\<}{\langle}
\renewcommand{\>}{\rangle}
\newcommand{\widebar}{\overline}
\def\reff#1{(\protect\ref{#1})}
\def\spose#1{\hbox to 0pt{#1\hss}}
\def\ltapprox{\mathrel{\spose{\lower 3pt\hbox{$\mathchar"218$}}
    \raise 2.0pt\hbox{$\mathchar"13C$}}}
\def\gtapprox{\mathrel{\spose{\lower 3pt\hbox{$\mathchar"218$}}
    \raise 2.0pt\hbox{$\mathchar"13E$}}}
\def\textprime{${}^\prime$}
\def\proof{\par\medskip\noindent{\sc Proof.\ }}
\def\firstproof{\par\medskip\noindent{\sc First Proof.\ }}
\def\secondproof{\par\medskip\noindent{\sc Second Proof.\ }}
\def\alternateproof{\par\medskip\noindent{\sc Alternate Proof.\ }}
\def\algebraicproof{\par\medskip\noindent{\sc Algebraic Proof.\ }}
\def\graphicalproof{\par\medskip\noindent{\sc Graphical Proof.\ }}
\def\combinatorialproof{\par\medskip\noindent{\sc Combinatorial Proof.\ }}
\def\proofof#1{\bigskip\noindent{\sc Proof of #1.\ }}
\def\firstproofof#1{\bigskip\noindent{\sc First Proof of #1.\ }}
\def\secondproofof#1{\bigskip\noindent{\sc Second Proof of #1.\ }}
\def\thirdproofof#1{\bigskip\noindent{\sc Third Proof of #1.\ }}
\def\algebraicproofof#1{\bigskip\noindent{\sc Algebraic Proof of #1.\ }}
\def\combinatorialproofof#1{\bigskip\noindent{\sc Combinatorial Proof of #1.\ }}
\def\completionofproofof#1{\bigskip\noindent{\sc Completion of the Proof of #1.\ }}
\def\sketchofproof{\par\medskip\noindent{\sc Sketch of proof.\ }}
\renewcommand{\qed}{ $\square$ \bigskip}
\newcommand{\myendremark}{ $\blacksquare$ \bigskip}
\def\half{ {1 \over 2} }
\def\third{ {1 \over 3} }
\def\twothird{ {2 \over 3} }
\def\smfrac#1#2{{\textstyle{#1\over #2}}}
\def\smhalf{ {\smfrac{1}{2}} }
\newcommand{\real}{\mathop{\rm Re}\nolimits}
\renewcommand{\Re}{\mathop{\rm Re}\nolimits}
\newcommand{\imag}{\mathop{\rm Im}\nolimits}
\renewcommand{\Im}{\mathop{\rm Im}\nolimits}
\newcommand{\sgn}{\mathop{\rm sgn}\nolimits}
\newcommand{\tr}{\mathop{\rm tr}\nolimits}
\newcommand{\supp}{\mathop{\rm supp}\nolimits}
\newcommand{\disc}{\mathop{\rm disc}\nolimits}
\newcommand{\diag}{\mathop{\rm diag}\nolimits}
\newcommand{\tridiag}{\mathop{\rm tridiag}\nolimits}
\newcommand{\AZ}{\mathop{\rm AZ}\nolimits}
\newcommand{\EAZ}{\mathop{\rm EAZ}\nolimits}
\newcommand{\NC}{\mathop{\rm NC}\nolimits}
\newcommand{\PF}{{\rm PF}}
\newcommand{\rk}{\mathop{\rm rk}\nolimits}
\newcommand{\perm}{\mathop{\rm perm}\nolimits}
\def\hboxscript#1{ {\hbox{\scriptsize\em #1}} }
\renewcommand{\emptyset}{\varnothing}
\newcommand{\eqdef}{\stackrel{\rm def}{=}}

\newcommand{\restrict}{\upharpoonright}

\newcommand{\compinv}{{\langle -1 \rangle}}   

\newcommand{\scra}{{\mathcal{A}}}
\newcommand{\scrb}{{\mathcal{B}}}
\newcommand{\scrc}{{\mathcal{C}}}
\newcommand{\scrd}{{\mathcal{D}}}
\newcommand{\scrdtilde}{{\widetilde{\mathcal{D}}}}
\newcommand{\scre}{{\mathcal{E}}}
\newcommand{\scrf}{{\mathcal{F}}}
\newcommand{\scrg}{{\mathcal{G}}}
\newcommand{\scrh}{{\mathcal{H}}}
\newcommand{\scri}{{\mathcal{I}}}
\newcommand{\scrj}{{\mathcal{J}}}
\newcommand{\scrk}{{\mathcal{K}}}
\newcommand{\scrl}{{\mathcal{L}}}
\newcommand{\scrlbar}{{\overline{\mathcal{L}}}}
\newcommand{\scrm}{{\mathcal{M}}}
\newcommand{\scrn}{{\mathcal{N}}}
\newcommand{\scro}{{\mathcal{O}}}
\newcommand\scroo{
  \mathchoice
    {{\scriptstyle\mathcal{O}}}
    {{\scriptstyle\mathcal{O}}}
    {{\scriptscriptstyle\mathcal{O}}}
    {\scalebox{0.6}{$\scriptscriptstyle\mathcal{O}$}}
  }
\newcommand{\scrp}{{\mathcal{P}}}
\newcommand{\scrq}{{\mathcal{Q}}}
\newcommand{\scrr}{{\mathcal{R}}}
\newcommand{\scrs}{{\mathcal{S}}}
\newcommand{\scrt}{{\mathcal{T}}}
\newcommand{\scrv}{{\mathcal{V}}}
\newcommand{\scrw}{{\mathcal{W}}}
\newcommand{\scrz}{{\mathcal{Z}}}
\newcommand{\SP}{{\mathcal{SP}}}
\newcommand{\ST}{{\mathcal{ST}}}

\newcommand{\bfa}{{\mathbf{a}}}
\newcommand{\bfb}{{\mathbf{b}}}
\newcommand{\bfc}{{\mathbf{c}}}
\newcommand{\bfd}{{\mathbf{d}}}
\newcommand{\bfe}{{\mathbf{e}}}
\newcommand{\bfff}{{\mathbf{f}}}
\newcommand{\bfg}{{\mathbf{g}}}
\newcommand{\bfh}{{\mathbf{h}}}
\newcommand{\bfj}{{\mathbf{j}}}
\newcommand{\bfi}{{\mathbf{i}}}
\newcommand{\bfk}{{\mathbf{k}}}
\newcommand{\bfl}{{\mathbf{l}}}
\newcommand{\bfL}{{\mathbf{L}}}
\newcommand{\bfm}{{\mathbf{m}}}
\newcommand{\bfn}{{\mathbf{n}}}
\newcommand{\bfp}{{\mathbf{p}}}
\newcommand{\bfr}{{\mathbf{r}}}
\newcommand{\bft}{{\mathbf{t}}}
\newcommand{\bfu}{{\mathbf{u}}}
\newcommand{\bfv}{{\mathbf{v}}}
\newcommand{\bfw}{{\mathbf{w}}}
\newcommand{\bfx}{{\mathbf{x}}}
\newcommand{\bfX}{{\mathbf{X}}}
\newcommand{\bfy}{{\mathbf{y}}}
\newcommand{\bfz}{{\mathbf{z}}}
\renewcommand{\k}{{\mathbf{k}}}
\newcommand{\n}{{\mathbf{n}}}
\newcommand{\vv}{{\mathbf{v}}}
\newcommand{\w}{{\mathbf{w}}}
\newcommand{\x}{{\mathbf{x}}}
\newcommand{\y}{{\mathbf{y}}}
\newcommand{\cc}{{\mathbf{c}}}
\newcommand{\zero}{{\mathbf{0}}}
\newcommand{\one}{{\mathbf{1}}}
\newcommand{\bmm}{{\mathbf{m}}}

\newcommand{\ahat}{{\widehat{a}}}
\newcommand{\Ghat}{{\widehat{G}}}
\newcommand{\Zhat}{{\widehat{Z}}}

\newcommand{\Pcheck}{{\widecheck{P}}}
\newcommand{\pcheck}{{\widecheck{p}}}

\newcommand{\C}{{\mathbb C}}
\newcommand{\D}{{\mathbb D}}
\newcommand{\Z}{{\mathbb Z}}
\newcommand{\N}{{\mathbb N}}
\newcommand{\Q}{{\mathbb Q}}
\newcommand{\PP}{{\mathbb P}}
\newcommand{\R}{{\mathbb R}}
\newcommand{\RR}{{\mathbb R}}
\newcommand{\E}{{\mathbb E}}

\newcommand{\Sym}{{\mathfrak{S}}}
\newcommand{\SymB}{{\mathfrak{B}}}
\newcommand{\Alt}{{\mathrm{Alt}}}

\newcommand{\germanA}{{\mathfrak{A}}}
\newcommand{\germanB}{{\mathfrak{B}}}
\newcommand{\germanQ}{{\mathfrak{Q}}}
\newcommand{\germanh}{{\mathfrak{h}}}

\newcommand{\myle}{\preceq}
\newcommand{\myge}{\succeq}
\newcommand{\mygt}{\succ}

\newcommand{\B}{{\sf B}}
\newcommand{\OB}{B^{\rm ord}}
\newcommand{\OS}{{\sf OS}}
\newcommand{\OO}{{\sf O}}
\newcommand{\OSP}{{\sf OSP}}
\newcommand{\Eu}{{\sf Eu}}
\newcommand{\ERR}{{\sf ERR}}
\newcommand{\sfB}{{\sf B}}
\newcommand{\sfD}{{\sf D}}
\newcommand{\sfE}{{\sf E}}
\newcommand{\sfG}{{\sf G}}
\newcommand{\sfJ}{{\sf J}}
\newcommand{\sfL}{{\sf L}}
\newcommand{\sfLhat}{{\widehat{{\sf L}}}}
\newcommand{\sfLcheck}{{\widecheck{{\sf L}}}}
\newcommand{\sfLtilde}{{\widetilde{{\sf L}}}}
\newcommand{\sfM}{{\sf M}}
\newcommand{\sfP}{{\sf P}}
\newcommand{\sfQ}{{\sf Q}}
\newcommand{\sfS}{{\sf S}}
\newcommand{\sfT}{{\sf T}}
\newcommand{\sfW}{{\sf W}}
\newcommand{\sfMV}{{\sf MV}}
\newcommand{\AMV}{{\sf AMV}}
\newcommand{\BM}{{\sf BM}}
\newcommand{\emIB}{B^{\rm irr}}
\newcommand{\emIP}{P^{\rm irr}}
\newcommand{\emOB}{B^{\rm ord}}
\newcommand{\emCB}{B^{\rm cyc}}
\newcommand{\emSC}{P^{\rm cyc}}

\newcommand{\cyc}{{\rm cyc}}
\newcommand{\e}{{\rm e}}
\newcommand{\ecyc}{{\rm ecyc}}
\newcommand{\epa}{{\rm epa}}
\newcommand{\iv}{{\rm iv}}
\newcommand{\pa}{{\rm pa}}
\newcommand{\pk}{{\rm p}}
\newcommand{\val}{{\rm v}}
\newcommand{\da}{{\rm da}}
\newcommand{\dd}{{\rm dd}}
\newcommand{\fp}{{\rm fp}}
\newcommand{\pkcyc}{{\rm pcyc}}
\newcommand{\valcyc}{{\rm vcyc}}
\newcommand{\dacyc}{{\rm dacyc}}
\newcommand{\ddcyc}{{\rm ddcyc}}
\newcommand{\pkpa}{{\rm ppa}}
\newcommand{\valpa}{{\rm vpa}}
\newcommand{\dapa}{{\rm dapa}}
\newcommand{\ddpa}{{\rm ddpa}}

\newcommand{\yp}{{y_\pk}}
\newcommand{\yptilde}{{\widetilde{y}_\pk}}
\newcommand{\yphat}{{\widehat{y}_\pk}}
\newcommand{\yv}{{y_\val}}
\newcommand{\yiv}{{y_\iv}}
\newcommand{\yda}{{y_\da}}
\newcommand{\ydatilde}{{\widetilde{y}_\da}}
\newcommand{\ydd}{{y_\dd}}
\newcommand{\yddtilde}{{\widetilde{y}_\dd}}
\newcommand{\yfp}{{y_\fp}}
\newcommand{\zp}{{z_\pk}}
\newcommand{\zv}{{z_\val}}
\newcommand{\zda}{{z_\da}}
\newcommand{\zdd}{{z_\dd}}

\newcommand{\lev}{{\rm lev}}
\newcommand{\stat}{{\rm stat}}
\newcommand{\umiv}{{\rm umiv}}
\newcommand{\umpa}{{\rm umpa}}
\newcommand{\mysteryone}{{\rm mys1}}
\newcommand{\mysterytwo}{{\rm mys2}}
\newcommand{\Asc}{{\rm Asc}}
\newcommand{\asc}{{\rm asc}}
\newcommand{\Des}{{\rm Des}}
\newcommand{\des}{{\rm des}}
\newcommand{\Exc}{{\rm Exc}}
\newcommand{\exc}{{\rm exc}}
\newcommand{\Wex}{{\rm Wex}}
\newcommand{\wex}{{\rm wex}}
\newcommand{\Fix}{{\rm Fix}}
\newcommand{\fix}{{\rm fix}}
\newcommand{\lrmax}{{\rm lrmax}}
\newcommand{\rlmax}{{\rm rlmax}}
\newcommand{\Rec}{{\rm Rec}}
\newcommand{\rec}{{\rm rec}}
\newcommand{\Arec}{{\rm Arec}}
\newcommand{\arec}{{\rm arec}}
\newcommand{\ERec}{{\rm ERec}}
\newcommand{\erec}{{\rm erec}}
\newcommand{\EArec}{{\rm EArec}}
\newcommand{\earec}{{\rm earec}}
\newcommand{\recarec}{{\rm recarec}}
\newcommand{\nonrec}{{\rm nonrec}}
\newcommand{\Cpeak}{{\rm Cpeak}}
\newcommand{\cpeak}{{\rm cpeak}}
\newcommand{\Cval}{{\rm Cval}}
\newcommand{\cval}{{\rm cval}}
\newcommand{\Cdasc}{{\rm Cdasc}}
\newcommand{\cdasc}{{\rm cdasc}}
\newcommand{\Cddes}{{\rm Cddes}}
\newcommand{\cddes}{{\rm cddes}}
\newcommand{\cdrise}{{\rm cdrise}}
\newcommand{\cdfall}{{\rm cdfall}}
\newcommand{\Peak}{{\rm Peak}}
\newcommand{\peak}{{\rm peak}}
\newcommand{\Val}{{\rm Val}}
\newcommand{\Dasc}{{\rm Dasc}}
\newcommand{\dasc}{{\rm dasc}}
\newcommand{\Ddes}{{\rm Ddes}}
\newcommand{\ddes}{{\rm ddes}}
\newcommand{\inv}{{\rm inv}}
\newcommand{\maj}{{\rm maj}}
\newcommand{\rs}{{\rm rs}}
\newcommand{\cross}{{\rm cr}}
\newcommand{\crosshat}{{\widehat{\rm cr}}}
\newcommand{\nest}{{\rm ne}}
\newcommand{\rodd}{{\rm rodd}}
\newcommand{\reven}{{\rm reven}}
\newcommand{\lodd}{{\rm lodd}}
\newcommand{\leven}{{\rm leven}}
\newcommand{\sg}{{\rm sg}}
\newcommand{\bl}{{\rm bl}}
\newcommand{\tran}{{\rm tr}}
\newcommand{\area}{{\rm area}}
\newcommand{\ret}{{\rm ret}}
\newcommand{\peaks}{{\rm peaks}}
\newcommand{\hl}{{\rm hl}}
\newcommand{\sll}{{\rm sl}}
\newcommand{\negg}{{\rm neg}}
\newcommand{\imp}{{\rm imp}}
\newcommand{\osg}{{\rm osg}}
\newcommand{\ons}{{\rm ons}}
\newcommand{\isg}{{\rm isg}}
\newcommand{\ins}{{\rm ins}}
\newcommand{\LL}{{\rm LL}}
\newcommand{\height}{{\rm ht}}
\newcommand{\as}{{\rm as}}

\newcommand{\ba}{{\bm{a}}}
\newcommand{\bahat}{{\widehat{\bm{a}}}}
\newcommand{\sfa}{{{\sf a}}}
\newcommand{\bb}{{\bm{b}}}
\newcommand{\bc}{{\bm{c}}}
\newcommand{\bchat}{{\widehat{\bm{c}}}}
\newcommand{\bd}{{\bm{d}}}
\newcommand{\bee}{{\bm{e}}}
\newcommand{\beh}{{\bm{eh}}}
\newcommand{\bff}{{\bm{f}}}
\newcommand{\bg}{{\bm{g}}}
\newcommand{\bh}{{\bm{h}}}
\newcommand{\bll}{{\bm{\ell}}}
\newcommand{\bp}{{\bm{p}}}
\newcommand{\bq}{{\bm{q}}}
\newcommand{\br}{{\bm{r}}}
\newcommand{\bs}{{\bm{s}}}
\newcommand{\bt}{{\bm{t}}}
\newcommand{\bu}{{\bm{u}}}
\newcommand{\bv}{{\bm{v}}}
\newcommand{\bw}{{\bm{w}}}
\newcommand{\bx}{{\bm{x}}}
\newcommand{\by}{{\bm{y}}}
\newcommand{\bz}{{\bm{z}}}
\newcommand{\bA}{{\bm{A}}}
\newcommand{\bB}{{\bm{B}}}
\newcommand{\bC}{{\bm{C}}}
\newcommand{\bE}{{\bm{E}}}
\newcommand{\bF}{{\bm{F}}}
\newcommand{\bG}{{\bm{G}}}
\newcommand{\bH}{{\bm{H}}}
\newcommand{\bI}{{\bm{I}}}
\newcommand{\bJ}{{\bm{J}}}
\newcommand{\bL}{{\bm{L}}}
\newcommand{\bLhat}{{\widehat{\bm{L}}}}
\newcommand{\bM}{{\bm{M}}}
\newcommand{\bN}{{\bm{N}}}
\newcommand{\bP}{{\bm{P}}}
\newcommand{\bQ}{{\bm{Q}}}
\newcommand{\bR}{{\bm{R}}}
\newcommand{\bS}{{\bm{S}}}
\newcommand{\bT}{{\bm{T}}}
\newcommand{\bW}{{\bm{W}}}
\newcommand{\bX}{{\bm{X}}}
\newcommand{\bY}{{\bm{Y}}}
\newcommand{\bIB}{{\bm{B}^{\rm irr}}}
\newcommand{\bOB}{{\bm{B}^{\rm ord}}}
\newcommand{\bOS}{{\bm{OS}}}
\newcommand{\bERR}{{\bm{ERR}}}
\newcommand{\bSP}{{\bm{SP}}}
\newcommand{\bMV}{{\bm{MV}}}
\newcommand{\bBM}{{\bm{BM}}}
\newcommand{\balpha}{{\bm{\alpha}}}
\newcommand{\balphapre}{{\bm{\alpha}^{\rm pre}}}
\newcommand{\bbeta}{{\bm{\beta}}}
\newcommand{\bgamma}{{\bm{\gamma}}}
\newcommand{\bdelta}{{\bm{\delta}}}
\newcommand{\bkappa}{{\bm{\kappa}}}
\newcommand{\bmu}{{\bm{\mu}}}
\newcommand{\bomega}{{\bm{\omega}}}
\newcommand{\bsigma}{{\bm{\sigma}}}
\newcommand{\btau}{{\bm{\tau}}}
\newcommand{\bphi}{{\bm{\phi}}}
\newcommand{\bphihat}{{\skew{3}\widehat{\vphantom{t}\protect\smash{\bm{\phi}}}}}
\newcommand{\bpsi}{{\bm{\psi}}}
\newcommand{\bxi}{{\bm{\xi}}}
\newcommand{\bzeta}{{\bm{\zeta}}}
\newcommand{\bone}{{\bm{1}}}
\newcommand{\bzero}{{\bm{0}}}

\newcommand{\Cbar}{{\overline{C}}}
\newcommand{\Dbar}{{\overline{D}}}
\newcommand{\dbar}{{\overline{d}}}
\newcommand{\Pbar}{{\bar{P}}}
\newcommand{\pbar}{{\bar{p}}}
\newcommand{\Lbar}{{\bar{L}}}
\newcommand{\Ubar}{{\bar{U}}}
\def\Btilde{{\widetilde{B}}}
\def\Ctilde{{\widetilde{C}}}
\def\Ftilde{{\widetilde{F}}}
\def\Gtilde{{\widetilde{G}}}
\def\Htilde{{\widetilde{H}}}
\def\Lhat{{\widehat{L}}}
\def\Ltilde{{\widetilde{L}}}
\def\Ptilde{{\widetilde{P}}}
\def\ptilde{{\widetilde{p}}}
\def\Chat{{\widehat{C}}}
\def\ctilde{{\widetilde{c}}}
\def\zbar{{\overline{Z}}}
\def\pitilde{{\widetilde{\pi}}}
\def\omegahat{{\widehat{\omega}}}

\newcommand{\sech}{{\rm sech}}

%
%
\newcommand{\sn}{{\rm sn}}
\newcommand{\cn}{{\rm cn}}
\newcommand{\dn}{{\rm dn}}
\newcommand{\sm}{{\rm sm}}
\newcommand{\cm}{{\rm cm}}

%
%
\newcommand{\zfz}{ {{}_0 \! F_0} }
\newcommand{\zfo}{ {{}_0  F_1} }
\newcommand{\ofz}{ {{}_1 \! F_0} }
\newcommand{\ofo}{ {{}_1 \! F_1} }
\newcommand{\oft}{ {{}_1 \! F_2} }

%
%
\newcommand{\FHyper}[2]{ {\tensor[_{#1 \!}]{F}{_{#2}}\!} }
\newcommand{\FHYPER}[5]{ {\FHyper{#1}{#2} \!\biggl(
   \!\!\begin{array}{c} #3 \\[1mm] #4 \end{array}\! \bigg|\, #5 \! \biggr)} }
\newcommand{\tfo}{ {\FHyper{2}{1}} }
\newcommand{\tfz}{ {\FHyper{2}{0}} }
\newcommand{\threefz}{ {\FHyper{3}{0}} }
\newcommand{\FHYPERbottomzero}[3]{ {\FHyper{#1}{0} \hspace*{-0mm}\biggl(
   \!\!\begin{array}{c} #2 \\[1mm] \hbox{---} \end{array}\! \bigg|\, #3 \! \biggr)} }
\newcommand{\FHYPERtopzero}[3]{ {\FHyper{0}{#1} \hspace*{-0mm}\biggl(
   \!\!\begin{array}{c} \hbox{---} \\[1mm] #2 \end{array}\! \bigg|\, #3 \! \biggr)} }

\newcommand{\phiHyper}[2]{ {\tensor[_{#1}]{\phi}{_{#2}}} }
\newcommand{\psiHyper}[2]{ {\tensor[_{#1}]{\psi}{_{#2}}} }
\newcommand{\PhiHyper}[2]{ {\tensor[_{#1}]{\Phi}{_{#2}}} }
\newcommand{\PsiHyper}[2]{ {\tensor[_{#1}]{\Psi}{_{#2}}} }
\newcommand{\phiHYPER}[6]{ {\phiHyper{#1}{#2} \!\left(
   \!\!\begin{array}{c} #3 \\ #4 \end{array}\! ;\, #5, \, #6 \! \right)\!} }
\newcommand{\psiHYPER}[6]{ {\psiHyper{#1}{#2} \!\left(
   \!\!\begin{array}{c} #3 \\ #4 \end{array}\! ;\, #5, \, #6 \! \right)} }
\newcommand{\PhiHYPER}[5]{ {\PhiHyper{#1}{#2} \!\left(
   \!\!\begin{array}{c} #3 \\ #4 \end{array}\! ;\, #5 \! \right)\!} }
\newcommand{\PsiHYPER}[5]{ {\PsiHyper{#1}{#2} \!\left(
   \!\!\begin{array}{c} #3 \\ #4 \end{array}\! ;\, #5 \! \right)\!} }
\newcommand{\zerophizero}{ {\phiHyper{0}{0}} }
\newcommand{\ophizero}{ {\phiHyper{1}{0}} }
\newcommand{\zphio}{ {\phiHyper{0}{1}} }
\newcommand{\ophio}{ {\phiHyper{1}{1}} }
\newcommand{\tphio}{ {\phiHyper{2}{1}} }
\newcommand{\tphiz}{ {\phiHyper{2}{0}} }
\newcommand{\tPhio}{ {\PhiHyper{2}{1}} }
\newcommand{\opsio}{ {\psiHyper{1}{1}} }

%
%
\newcommand{\stirlingsubset}[2]{\genfrac{\{}{\}}{0pt}{}{#1}{#2}}
\newcommand{\stirlingcycle}[2]{\genfrac{[}{]}{0pt}{}{#1}{#2}}
\newcommand{\assocstirlingsubset}[3]{{\genfrac{\{}{\}}{0pt}{}{#1}{#2}}_{\! \ge #3}}
\newcommand{\genstirlingsubset}[4]{{\genfrac{\{}{\}}{0pt}{}{#1}{#2}}_{\! #3,#4}}
\newcommand{\irredstirlingsubset}[2]{{\genfrac{\{}{\}}{0pt}{}{#1}{#2}}^{\!\rm irr}}
\newcommand{\euler}[2]{\genfrac{\langle}{\rangle}{0pt}{}{#1}{#2}}
\newcommand{\eulergen}[3]{{\genfrac{\langle}{\rangle}{0pt}{}{#1}{#2}}_{\! #3}}
\newcommand{\eulersecond}[2]{\left\langle\!\! \euler{#1}{#2} \!\!\right\rangle}
\newcommand{\eulersecondgen}[3]{{\left\langle\!\! \euler{#1}{#2} \!\!\right\rangle}_{\! #3}}
\newcommand{\binomvert}[2]{\genfrac{\vert}{\vert}{0pt}{}{#1}{#2}}
\newcommand{\binomsquare}[2]{\genfrac{[}{]}{0pt}{}{#1}{#2}}
\newcommand{\doublebinom}[2]{\left(\!\! \binom{#1}{#2} \!\!\right)}
\newcommand{\lahnum}[2]{\genfrac{\lfloor}{\rfloor}{0pt}{}{#1}{#2}}
\newcommand{\rlahnum}[3]{{\genfrac{\lfloor}{\rfloor}{0pt}{}{#1}{#2}}_{\! #3}}

%
%
\newcommand{\Lna}{{L_n^{(\alpha)}}}
\newcommand{\scrlna}{{\scrl_n^{(\alpha)}}}
\newcommand{\scrlbarna}{{\scrlbar_n^{(\alpha)}}}
\newcommand{\scrlhatna}{{\widehat{\scrl}_n^{(\alpha)}}}
\newcommand{\scrlhatbarna}{{\overline{\widehat{\scrl}}_n^{(\alpha)}}}
\newcommand{\scrlbarnka}{{\scrlbar_{n,k}^{(\alpha)}}}
\newcommand{\scrlbarnzeroa}{{\scrlbar_{n,0}^{(\alpha)}}}
\newcommand{\Lnb}{{L_n^{[\beta]}}}
\newcommand{\scrlnb}{{\scrl_n^{[\beta]}}}
\newcommand{\scrlbarnb}{{\scrlbar_n^{[\beta]}}}
\newcommand{\Lah}{{\textrm{Lah}}}
\newcommand{\LD}{{\mathbf{LD}}}
\newcommand{\MLD}{{\mathbf{MLD}}}


\newenvironment{sarray}{
             \textfont0=\scriptfont0
             \scriptfont0=\scriptscriptfont0
             \textfont1=\scriptfont1
             \scriptfont1=\scriptscriptfont1
             \textfont2=\scriptfont2
             \scriptfont2=\scriptscriptfont2
             \textfont3=\scriptfont3
             \scriptfont3=\scriptscriptfont3
           \renewcommand{\arraystretch}{0.7}
           \begin{array}{l}}{\end{array}}

\newenvironment{scarray}{
             \textfont0=\scriptfont0
             \scriptfont0=\scriptscriptfont0
             \textfont1=\scriptfont1
             \scriptfont1=\scriptscriptfont1
             \textfont2=\scriptfont2
             \scriptfont2=\scriptscriptfont2
             \textfont3=\scriptfont3
             \scriptfont3=\scriptscriptfont3
           \renewcommand{\arraystretch}{0.7}
           \begin{array}{c}}{\end{array}}


\newcommand*\circled[1]{\tikz[baseline=(char.base)]{
  \node[shape=circle,draw,inner sep=1pt] (char) {#1};}}
\newcommand{\ostar}{{\circledast}}
\newcommand{\ostarN}{{\,\circledast_{\vphantom{\dot{N}}N}\,}}
\newcommand{\ostarPsi}{{\,\circledast_{\vphantom{\dot{\Psi}}\Psi}\,}}
\newcommand{\starN}{{\,\ast_{\vphantom{\dot{N}}N}\,}}
\newcommand{\starpsi}{{\,\ast_{\vphantom{\dot{\bpsi}}\!\bpsi}\,}}
\newcommand{\starone}{{\,\ast_{\vphantom{\dot{1}}1}\,}}
\newcommand{\startwo}{{\,\ast_{\vphantom{\dot{2}}2}\,}}
\newcommand{\starinfty}{{\,\ast_{\vphantom{\dot{\infty}}\infty}\,}}
\newcommand{\starT}{{\,\ast_{\vphantom{\dot{T}}T}\,}}

\newcommand*{\Scale}[2][4]{\scalebox{#1}{$#2$}}

\newcommand*{\Scaletext}[2][4]{\scalebox{#1}{#2}} 

\clearpage 

\enlargethispage{6\baselineskip}
\tableofcontents
\thispagestyle{empty}   

\clearpage

\section{Introduction and statement of main results}

In a seminal 1980 paper, Flajolet \cite{Flajolet_80}
showed that the coefficients in the Taylor expansion
of the generic Stieltjes-type (resp.\ Jacobi-type) continued fraction
--- which he called the {\em Stieltjes--Rogers}\/
 (resp.\ {\em Jacobi--Rogers}\/) {\em polynomials}\/ ---
can be interpreted as the generating polynomials
for Dyck (resp.\ Motzkin) paths with specified height-dependent weights.
More recently it was independently discovered by several authors
\cite{Fusy_15,Oste_15,Josuat-Verges_18,Sokal_totalpos}
that Thron-type continued fractions also have an interpretation of this kind:
namely, their Taylor coefficients
--- which we call, by analogy, the {\em Thron--Rogers polynomials}\/ ---
can be interpreted as the generating polynomials 
for Schr\"oder paths with specified height-dependent weights.

In a recent paper \cite{latpath_SRTR}
we presented an infinite sequence
of generalizations of the Stieltjes--Rogers and Thron--Rogers polynomials,
which are parametrized by an integer $m \ge 1$
and reduce to the classical Stieltjes--Rogers and Thron--Rogers polynomials
when $m=1$;
they are the generating polynomials of $m$-Dyck and $m$-Schr\"oder paths,
respectively, with height-dependent weights,
and are also the Taylor coefficients of certain branched continued fractions.
We proved that these generalizations all possess the fundamental property of
coefficientwise Hankel-total positivity \cite{Sokal_flajolet,Sokal_totalpos},
jointly in all the (infinitely many) indeterminates.
These facts were known when $m = 1$ \cite{Sokal_flajolet,Sokal_totalpos}
but were new when $m > 1$.
By specializing the indeterminates we were able to give many examples
of Hankel-totally positive sequences
whose generating functions do not possess nice classical continued fractions.
(The concept of Hankel-total positivity \cite{Sokal_flajolet,Sokal_totalpos}
 will be explained in more detail later in this Introduction.)
Similar considerations apply to the $m$-Jacobi--Rogers polynomials,
which are the generating polynomials of $m$-\L{}ukasiewicz paths;
but here the Hankel-total positivity is a more delicate matter,
for which a sufficient (but not necessary) condition
is the total positivity of a lower-Hessenberg production matrix.

In the present paper we will apply the $m$-Jacobi--Rogers theory with $m=2$
to prove the coefficientwise Hankel-total positivity
for some multivariate generalizations of the Laguerre, rook and Lah polynomials.
Our constructions are motivated in part by the work of
Coussement and Van Assche \cite{Coussement_03}
on the multiple orthogonal polynomials associated to weights
based on modified Bessel functions of the first kind $I_\alpha$;
we will explain this unexpected connection in Section~\ref{sec.MOP}.
Our method
--- to use production matrices for a binomial row-generating matrix,
exploiting the theory of exponential Riordan arrays ---
was employed previously by ourselves
\cite{latpath_lah,forests_totalpos,Chen-Sokal_trees_totalpos}
and Zhu \cite{Zhu_21,Zhu_22};
but the present case is more delicate than these previous applications
because the total positivity of the quadridiagonal production matrix
(Section~\ref{subsec.prodmat.TP.multivariate.quadridiagonal})
is decidedly nontrivial.

The starting point of this work was the fact,
due to Gantmakher and Krein \cite{Gantmakher_37},
that a Hankel matrix of real numbers is totally positive
if and only if the underlying sequence is a Stieltjes moment sequence
(see Section~\ref{subsec.TP} below),
combined with the classical fact that,
for each $\alpha \ge -1$ and $x \ge 0$,
the sequence $(\scrlna(x))_{n \ge 0}$ of monic unsigned Laguerre polynomials
[defined in \reff{def.scrlna} below]
is a Stieltjes moment sequence [cf.~\reff{eq.scrlna.stieltjes} below].
This led Sylvie Corteel and one of the authors (A.D.S.)\ to conjecture,
a few years ago \cite{Corteel-Sokal_conj_Laguerre},
that the sequence $(\scrlna(x))_{n \ge 0}$ is
{\em coefficientwise}\/ Hankel-totally positive
(see Section~\ref{subsec.TP} for the definition).
The purpose of this paper is to prove this conjecture,
along with a vast multivariate generalization.
\enlargethispage*{3\baselineskip}
The univariate conjecture was also recently proven by Zhu \cite{Zhu_21,Zhu_22}.


\subsection{Univariate Laguerre polynomials}

The Laguerre polynomials $\Lna$ are conventionally defined as
\cite{Rainville_60,Szego_75,Andrews_99,Ismail_05}
\be
  \Lna(x)  \;=\;
  \sum_{k=0}^n \binom{n+\alpha}{n-k} \, {(-x)^k \over k!}
 \label{def.Lna}
\ee
and have the generating function
\be
  \sum_{n=0}^\infty \Lna(x) \: t^n
  \;=\;
  (1-t)^{-(1+\alpha)} \, e^{-xt/(1-t)}
  \;.
 \label{eq.genfn.1}
\ee
However, for our purposes it is more convenient to work with the
\textbfit{monic unsigned Laguerre polynomials}
\begin{eqnarray}
   \scrlna(x)
   \;\eqdef\;
   n! \, \Lna(-x)
   & = &
   \sum_{k=0}^n {n! \over k!} \, \binom{n+\alpha}{n-k} \, x^k
   \;=\;
   \sum_{k=0}^n \binom{n}{k} \, (n+\alpha)^{\underline{n-k}} \, x^k
        \nonumber \\[2mm]
   & = &
   (1+\alpha)^{\overline{n}} \; \FHYPER{1}{1}{-n}{1+\alpha}{-x}
 \label{def.scrlna}
\end{eqnarray}
[where $\rho^{\underline{n}} \eqdef \rho(\rho-1) \cdots (\rho-n+1)$
 and $\rho^{\overline{n}} \eqdef \rho(\rho+1) \cdots (\rho+n-1)$],
or the \textbfit{reversed monic unsigned Laguerre polynomials}
\begin{eqnarray}
   \scrlbarna(x)
   \;\eqdef\; 
   n! \, x^n \, \Lna(-1/x)
   & = &
   \sum_{k=0}^n \binom{n}{k} \, \binom{n+\alpha}{k} \, k! \, x^k
   \;=\;
   \sum_{k=0}^n \binom{n}{k} \, (n+\alpha)^{\underline{k}} \, x^k
        \nonumber \\[2mm]
   & = &
   \FHYPERbottomzero{2}{-n,\,-n-\alpha}{x}
   \;.
 \label{def.scrlbarna}
\end{eqnarray}
Note that $\scrlna(x)$ and $\scrlbarna(x)$ are polynomials
jointly in $x$ and $\alpha$, with nonnegative integer coefficients.
The first few $\scrlna(x)$ are thus
\begin{subeqnarray}
   \scrl_0^{(\alpha)}(x)  & = &  1              \\
   \scrl_1^{(\alpha)}(x)  & = &  (1+\alpha) \,+\, x            \\
   \scrl_2^{(\alpha)}(x)  & = &
        (1+\alpha)(2+\alpha) \,+\, 2(2+\alpha)x \,+\, x^2   \\
   \scrl_3^{(\alpha)}(x)  & = &
        (1+\alpha)(2+\alpha)(3+\alpha) \,+\, 3(2+\alpha)(3+\alpha)x
        \,+\, 3(3+\alpha) x^2 \,+\, x^3
     \qquad
\end{subeqnarray}
It follows from \reff{eq.genfn.1} that these polynomials
have the exponential generating functions
\begin{eqnarray}
  \sum_{n=0}^\infty \scrlna(x) \: {t^n \over n!}
  & = &
  (1-t)^{-(1+\alpha)} \, e^{xt/(1-t)}
     \label{eq.egf.scrlna}  \\[2mm]
  \sum_{n=0}^\infty \scrlbarna(x) \: {t^n \over n!}
  & = &
  (1-xt)^{-(1+\alpha)} \, e^{t/(1-xt)}
     \label{eq.egf.scrlbarna}
\end{eqnarray}

A major role will be played in what follows by the
\textbfit{coefficient matrix} of the monic unsigned Laguerre polynomials:
it is the unit-lower-triangular matrix $\sfL^{(\alpha)}$ with entries
\be
   (\sfL^{(\alpha)})_{n,k}
   \;=\;
   {n! \over k!} \, \binom{n+\alpha}{n-k}
   \;=\;
   \binom{n}{k} \, (n+\alpha)^{\underline{n-k}}
   \;=\;
   \binom{n}{k} \, (1+\alpha+k)^{\overline{n-k}}
   \;.
 \label{def.coeffmat}
\ee
This matrix is an exponential Riordan array
(see~Section~\ref{sec.exp_riordan})
$\scrr[F,G]$ with $F(t) = (1-t)^{-(1+\alpha)}$ and $G(t) = t/(1-t)$:
compare \reff{eq.egf.scrlna} with \reff{eq.exp_riordan.bivariate_egf} below.
It is also the matrix of generalized Stieltjes--Rogers polynomials
of the first kind
(see Section~\ref{subsec.SR})
for the Stieltjes-type continued fraction associated to the sequence
${\bigl( (1+\alpha)^{\overline{n}} \bigr)_{n \ge 0}}$:
see Proposition~\ref{prop.L.SR}.

An important role will also be played by the
\textbfit{binomial row-generating matrix}
$\sfL^{(\alpha)} B_x$,
where $B_x$ is the weighted binomial matrix
\be
   (B_x)_{ij}  \;=\;  \binom{i}{j} \, x^{i-j}
 \label{def.Bx}
\ee
(note that it too is unit-lower-triangular).
The matrix $\sfL^{(\alpha)} B_x$ is an exponential Riordan array
$\scrr[e^{xG} F,G]$ with $F$ and $G$ as above
(see Corollary~\ref{cor.exp_riordan.Bx}).
The zeroth column of $\sfL^{(\alpha)} B_x$
consists of the monic unsigned Laguerre polynomials $\scrl_{n}^{(\alpha)}(x)$.
More generally, it can be shown
(see Proposition~\ref{prop.lagcoeff.Bx.NEW}) that
\be
   (\sfL^{(\alpha)} B_x)_{n,k}
   \;=\;
   {1 \over k!} \, {d^k \over dx^k} \, \scrlna(x)
   \;=\;
   \binom{n}{k} \, \scrl_{n-k}^{(\alpha+k)}(x)
   \;.
 \label{eq.L.Bx}
\ee

\bigskip

Several specializations of the Laguerre polynomials
to integer values of $\alpha$
correspond to well-known combinatorial objects:

1) The \textbfit{rook polynomial} $R_B(x)$ of a chessboard~$B$
is the generating polynomial for placements of
zero or more nonattacking rooks on~$B$,
with a weight~$x$ for each rook.\footnote{
   See \cite[chapters~7 and 8]{Riordan_58}
       \cite[sections~2.3 and 2.4]{Stanley_12}
   for further discussion of rook polynomials.
}
In particular, for an $m \times n$ rectangular chessboard,
the number of ways of placing $k$ non-attacking rooks is
$\binom{m}{k} \binom{n}{k} k!$:
we choose $k$ rows, $k$ columns, and a permutation connecting them.
It follows that the rook polynomial of an $n \times n$ square chessboard
is the reversed monic unsigned Laguerre polynomial
$\scrlbar_n^{(0)}(x) = \sum_{k=0}^n \binom{n}{k}^2 k! \, x^k$;
and more generally, for any integer $\alpha \ge -n$,
the rook polynomial of an $n \times (n+\alpha)$ rectangular chessboard is
$\scrlbar_n^{(\alpha)}(x) =
 \sum_{k=0}^n \binom{n}{k} \binom{n+\alpha}{k} k! \, x^k$.\footnote{
    See \cite[A144084]{OEIS} for the triangular array
    corresponding to the rook polynomial of the $n \times n$ chessboard,
    and see \cite[A002720]{OEIS} for the row sums.

    \ \ 
    Note also that a rook configuration of an $m \times n$ chessboard
    can equivalently be viewed as a matching of the complete bipartite graph
    $K_{m,n}$;  so the rook polynomial of the $m \times n$ chessboard
    is also the matching polynomial \cite[p.~334]{Lovasz_86} of $K_{m,n}$.
}
The first few rook polynomials for square chessboards are
\begin{subeqnarray}
   \scrlbar_0^{(0)}(x)  & = &  1              \\
   \scrlbar_1^{(0)}(x)  & = &  1 + x \\
   \scrlbar_2^{(0)}(x)  & = &  1 + 4 x + 2 x^2 \\
   \scrlbar_3^{(0)}(x)  & = &  1 + 9 x + 18 x^2 + 6 x^3 \\
   \scrlbar_4^{(0)}(x)  & = &  1 + 16 x + 72 x^2 + 96 x^3 + 24 x^4
     \qquad
\end{subeqnarray}

2) The {\em Lah number}\/ $\lahnum{n}{k}$ is
the number of partitions of an $n$-element set
into $k$ nonempty linearly ordered blocks (also called {\em lists}\/);
we set $\lahnum{0}{k} = \delta_{k0}$.
The \textbfit{Lah polynomial} is then defined as
$\Lah_n(x) = \sum_{k=0}^n \lahnum{n}{k} \: x^k$.
It is easy to see that the Lah numbers have the explicit expression
\be
   \lahnum{n}{k}  \;=\;  {n! \over k!} \binom{n-1}{n-k}
           \;=\;  \begin{cases}
                      \delta_{k0}       & \textrm{if $n=0$} \\[2mm]
                      {n! \over k!}  \binom{n-1}{k-1}
                                        & \textrm{if $n \ge 1$}
                  \end{cases}
\ee
It follows that $\Lah_n(x)$
is the monic unsigned Laguerre polynomial $\scrl^{(-1)}_n(x)$.\footnote{
   See \cite[A105278/A008297/A066667]{OEIS} for further information
   on the Lah numbers and Lah polynomials.
}
The first few Lah polynomials are
\begin{subeqnarray}
   \scrl_0^{(-1)}(x)  & = &  1              \\
   \scrl_1^{(-1)}(x)  & = &  \hphantom{0 +} x \\
   \scrl_2^{(-1)}(x)  & = &  \hphantom{0 +} 2 x + x^2  \\
   \scrl_3^{(-1)}(x)  & = &  \hphantom{0 +} 6 x + 6 x^2 + x^3  \\
   \scrl_4^{(-1)}(x)  & = &  \hphantom{0 +} 24 x + 36 x^2 + 12 x^3 + x^4
     \qquad
\end{subeqnarray}

3) More generally, for any integer $r \ge 0$,
the {\em $r$-Lah number}\/ $\rlahnum{n}{k}{r}$
is the number of partitions of an $(n+r)$-element set
into $k+r$ nonempty linearly ordered blocks,
with the restriction that $r$ distinguished elements
must belong to distinct blocks \cite{Nyul_15}.
The \textbfit{$\bm{r}$-Lah polynomial} is then defined as
$\Lah^{[r]}_n(x) = \sum_{k=0}^n \rlahnum{n}{k}{r} \: x^k$.
Clearly, $\Lah^{[0]}_n(x) = \Lah_n(x)$
and $\Lah^{[1]}_n(x) = \Lah_{n+1}(x)/x$.
It is not difficult to show \cite[Theorem~3.7]{Nyul_15}
that the $r$-Lah numbers have the explicit expression
\be
   \rlahnum{n}{k}{r}
   \;=\;
   {n! \over k!} \binom{n+2r-1}{n-k}
   \;=\;
   \begin{cases}
       \delta_{k0}                             & \textrm{if $n=0$} \\[2mm]
       {n! \over k!}  \binom{n+2r-1}{k+2r-1}   & \textrm{if $n \ge 1$}
   \end{cases}
\ee
It follows that $\Lah^{[r]}_n(x)$
is the monic unsigned Laguerre polynomial $\scrl^{(2r-1)}_n(x)$.\footnote{
   See \cite[A143497/A143498/A143499]{OEIS} for further information
   on the $r$-Lah numbers and $r$-Lah polynomials.
}

In summary, the rook polynomials correspond to
$\alpha = 0$ for square chessboards
and $\alpha = 1,2,3,\ldots$ for rectangular chessboards,
the Lah polynomials correspond to $\alpha = -1$,
and the $r$-Lah polynomials correspond to $\alpha = 2r-1$.

\subsection{Combinatorial interpretation in terms of Laguerre digraphs}

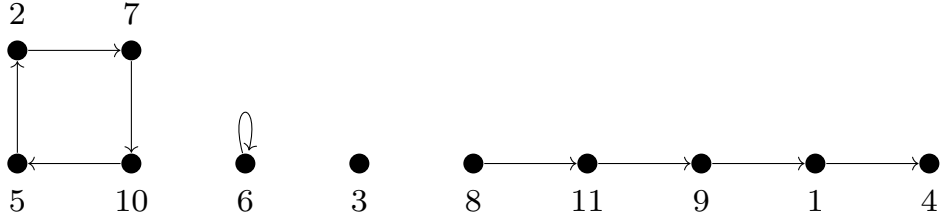
\begin{figure}
\centering
\begin{tikzpicture}[scale=1.5, every node/.style={scale=1.5}]

	\node[circle,fill=black,inner sep=1pt,minimum size=5pt] (e) at (0,0) {};
	\node[circle,fill=black,inner sep=1pt,minimum size=5pt] (j) at (1,0) {};
	\node[circle,fill=black,inner sep=1pt,minimum size=5pt] (g) at (1,1) {};
	\node[circle,fill=black,inner sep=1pt,minimum size=5pt] (b) at (0,1) {};

	\node[circle,fill=black,inner sep=1pt,minimum size=5pt] (f) at (2,0) {} edge [in=45,out=135, loop above] node {} ();

	\node[circle,fill=black,inner sep=1pt,minimum size=5pt] (c) at (3,0) {};

	\node[circle,fill=black,inner sep=1pt,minimum size=5pt] (h) at (4,0) {};
	\node[circle,fill=black,inner sep=1pt,minimum size=5pt] (k) at (5,0) {};
	\node[circle,fill=black,inner sep=1pt,minimum size=5pt] (i) at (6,0) {};
	\node[circle,fill=black,inner sep=1pt,minimum size=5pt] (a) at (7,0) {};
	\node[circle,fill=black,inner sep=1pt,minimum size=5pt] (d) at (8,0) {};

	\node[below = -0.2pt of a] {{\scriptsize 1}};
	\node[above = -0.2pt of b] {{\scriptsize 2}};
	\node[below = -0.2pt of c] {{\scriptsize 3}};
	\node[below = -0.2pt of d] {{\scriptsize 4}};
	\node[below = -0.2pt of e] {{\scriptsize 5}};
	\node[below = -0.2pt of f] {{\scriptsize 6}};
	\node[above = -0.2pt of g] {{\scriptsize 7}};
	\node[below = -0.2pt of h] {{\scriptsize 8}};
	\node[below = -0.2pt of i] {{\scriptsize 9}};
	\node[below = -0.2pt of j] {{\scriptsize 10}};
	\node[below = -0.2pt of k] {{\scriptsize 11}};

	\graph [multi] { 
	(b) -> (g) -> (j) -> (e);
	(e) -> (b);
	(f) -> [bend above] (f);
	(h) -> (k) -> (i) -> (a) -> (d);
	};

\end{tikzpicture}
\caption{A Laguerre digraph on $11$ vertices with $2$ cycles
   (one of which is a loop) and $2$ paths (one of which is an isolated vertex).}
\label{fig.laguerre.eg}
\end{figure}

Four decades ago, Foata and Strehl \cite{Foata_84}
introduced a beautiful combinatorial interpretation
of the Laguerre polynomials,
which will form the starting point for our work.\footnote{
   See also the recent paper of Strehl \cite{Strehl_16}
   for further applications of this formalism.
}
Let us define a \textbfit{Laguerre digraph}
to be a digraph in which each vertex has out-degree 0 or 1
and in-degree 0 or 1.
It follows that each weakly connected component of a Laguerre digraph
is either a directed path of some length $\ell \ge 0$
(where a path of length 0 is an isolated vertex)
or else a directed cycle of some length $\ell \ge 1$
(where a cycle of length 1 is a loop):
 see Figure~\ref{fig.laguerre.eg}.
For each integer $n \ge 0$,
let us write $\LD_n$ for the set of Laguerre digraphs
on the vertex set $[n] \eqdef \{1,\ldots,n\}$;
and for $n \ge k \ge 0$,
let us write $\LD_{n,k}$ for the set of Laguerre digraphs
on the vertex set $[n]$ with $k$ paths.
For $G \in \LD_n$,
we write $\e(G)$ [resp.\ $\cyc(G)$, $\pa(G)$]
for the number of edges (resp.\ cycles, paths) in $G$,
and observe that $\e(G) = n - \pa(G)$.
Foata and Strehl \cite{Foata_84} then showed that
\be
   \scrlna(x)
   \;=\;
   \sum_{G \in \LD_n}  x^{\pa(G)} \, (1+\alpha)^{\cyc(G)}
 \label{eq.foata.0}
\ee
or equivalently
\be
   \scrlbarna(x)
   \;=\;
   \sum_{G \in \LD_n}  x^{\e(G)} \, (1+\alpha)^{\cyc(G)}
   \;.
 \label{eq.foata}
\ee
Indeed, the proof of \reff{eq.foata.0}/\reff{eq.foata} is an easy argument
using the exponential formula \cite[chapter~5]{Stanley_99},
or equivalently, the theory of species \cite{Bergeron_98}.
Let us do \reff{eq.foata.0}:
The number of directed paths on $n \ge 1$ vertices is $n!$,
so with a weight $x$ per path they have exponential generating function
$xt/(1-t)$.
The number of directed cycles on $n \ge 1$ vertices is $(n-1)!$,
so with a weight $1+\alpha$ per cycle
they have exponential generating function $- (1+\alpha) \log(1-t)$.
A Laguerre digraph is a disjoint union of paths and cycles,
so by the exponential formula it has exponential generating function
\be
   \exp\Bigl[ {xt \over 1-t} \,-\, (1+\alpha) \log(1-t) \Bigr]
   \;=\;
   (1-t)^{-(1+\alpha)} \, e^{xt/(1-t)}
   \;,
\ee
which coincides with \reff{eq.egf.scrlna}.\footnote{
   Foata and Strehl \cite{Foata_84} also gave a direct combinatorial proof
   that \reff{eq.foata.0}/\reff{eq.foata} is equivalent to
   \reff{def.scrlna}/\reff{def.scrlbarna};
   this requires a bit more work \cite[Lemma~2.1]{Foata_84}.
}
%
%
The principal virtue of the combinatorial model
\reff{eq.foata.0}/\reff{eq.foata}
is that it treats $\alpha$ as an indeterminate;
it need not be an integer.\footnote{
   The formulae \reff{eq.foata.0}/\reff{eq.foata} also make clear
   (as Rota \cite[p.~206]{Mullin_70} realized many years ago)
   that the natural parameter for the Laguerre polynomials
   is $\lambda \eqdef 1+\alpha$, not $\alpha$.
   But it would be too confusing to try to change the notation
   at this late date.
}

The identity \reff{eq.foata.0} can equivalently be understood
as a combinatorial representation for the coefficient matrix $\sfL^{(\alpha)}$:
\be
   (\sfL^{(\alpha)})_{n,k}
   \;=\;
   \sum_{G \in \LD_{n,k}}  (1+\alpha)^{\cyc(G)}
   \;.
 \label{eq.foata.coeff}
\ee
This formula will be the starting point for our multivariate generalizations,
in which we will enumerate Laguerre digraphs with additional statistics.

\bigskip

{\bf Remarks.}
1.  There is a one-to-one correspondence
between Laguerre digraphs and rook configurations on an $n \times n$
chessboard: for a rook at position $(i,j)$ we draw an edge $i \to j$.
This explains why $\alpha = 0$ corresponds to the rook polynomial,
for which there is no special weighting of cycles.
However, we prefer the digraph point of view,
where the cycles (and therefore the role of $\alpha$) can be seen more clearly.

2.  Note also from \reff{eq.foata.0} that the case $\alpha = -1$
leads to Laguerre digraphs with no cycles
--- that is, only directed paths ---
with a weight $x$ per path.
These digraphs are in obvious one-to-one correspondence
with partitions of $[n]$ into linearly ordered blocks;
this explains why $\alpha = -1$ corresponds to the Lah polynomials.
\myendremark

\subsection{A first multivariate generalization}
   \label{subsec.multivariate.first}


It is now natural to generalize the Laguerre polynomials
by introducing further statistics into the digraph model.
Here is one way:
Let us write $\e_-(G), \e_0(G), \e_+(G)$
for the number of edges $i \to j$ in $G$
where $j < i$, $j = i$ or $j > i$, respectively.
We then generalize the coefficient matrix \reff{eq.foata.coeff}
by introducing separate variables $v_-, v_0, v_+$ for the three types of edges:
\be
   \sfL^{(\alpha)}(v_-, v_0, v_+)_{n,k}
   \;\eqdef\;
   \sum_{G \in \LD_{n,k}}  v_-^{\e_-(G)} \, v_0^{\e_0(G)} \, v_+^{\e_+(G)} 
                              \, (1+\alpha)^{\cyc(G)}
   \;.
 \label{def.coeffmat.gen}
\ee
This polynomial is homogeneous of degree $n-k$ in $v_-, v_0, v_+$.
We call the matrix $\sfL^{(\alpha)}(v_-, v_0, v_+)$
the \textbfit{(first) multivariate Laguerre coefficient matrix}.
In Section~\ref{subsec.egf.egf}
we will compute the bivariate exponential generating function
for this matrix: see \reff{eq.first_multivariate.egf}.
By specializing $v_- = v_0 = v_+ = v$,
we recover a rescaled version of the univariate Laguerre coefficient matrix
\be
   \sfL^{(\alpha)}(v,v,v)_{n,k}
   \;=\;
   (\sfL^{(\alpha)})_{n,k} \, v^{n-k}
   \;.
\ee

Let us then introduce the row-generating polynomials
of the multivariate Laguerre coefficient matrix:
\begin{eqnarray}
   \scrlna(x; v_-, v_0, v_+)
   & \eqdef &
   \sum_{k=0}^n \sfL^{(\alpha)}(v_-, v_0, v_+)_{n,k} \: x^k
       \label{def.multivariate.rowgen.scrl}
          \\[2mm]
   \scrlbarna(x; v_-, v_0, v_+)
   & \eqdef &
   \sum_{k=0}^n \sfL^{(\alpha)}(v_-, v_0, v_+)_{n,k} \: x^{n-k}
       \label{def.multivariate.rowgen.scrlbar}
\end{eqnarray}
Of course, by homogeneity \reff{def.multivariate.rowgen.scrlbar}
can also be understood as
\be
   \scrlbarna(x; v_-, v_0, v_+)
   \;=\;
   \sum_{k=0}^n \sfL^{(\alpha)}(xv_-, xv_0, xv_+)_{n,k}
   \;.
       \label{def.multivariate.rowgen.scrlbar.bis}
\ee


\bigskip

{\bf Remark.}
In the rook model, the three $v$ variables correspond
to giving different weights for rooks
below the diagonal, on the diagonal, or above the diagonal.
To our knowledge this type of weighting has not previously been considered,
even in the rook case $\alpha  = 0$.
However, two special cases generalize previously known results:

\medskip

1)  If we set $v_+ = 0$, then the only possible connected components of $G$
are decreasing paths (of length $\ell \ge 0$) and loops.
It follows that
\be
   \sfL^{(\alpha)}(v_-, v_0, 0)_{n,k}
   \;=\;
   \sum_{i=0}^n \binom{n}{i} \, [(1+\alpha)v_0]^i  \:
      \stirlingsubset{n-i}{k} \, v_-^{n-i-k}
   \;,
 \label{eq.rooks_decreasing}
\ee
where $\stirlingsubset{n}{k}$ denotes
the number of partitions of an $n$-element set into $k$ nonempty blocks:
we choose $i$ loops and then partition the remaining vertices
into $k$ nonempty blocks.
%
%
%
Further specializing to $\alpha = 0$ and $v_0 = 0$
yields the well-known \cite[Corollary~2.4.2]{Stanley_12}
formula for the counting of rook configurations
on an $n \times n$ triangular board:
\be
    \sfL^{(0)}(1,0,0)_{n,k} \;=\; \stirlingsubset{n}{k}
   \;,
  \label{eq.rook.triangular.1}
\ee
where a configuration with $k$ paths has $n-k$ edges and hence $n-k$ rooks.
Similarly, for $\alpha = 0$ and $v_0 = 1$ we get
\be
    \sfL^{(0)}(1,1,0)_{n,k} \;=\; \stirlingsubset{n+1}{k+1}
    \;,
  \label{eq.rook.triangular.2}
\ee
which follows from \reff{eq.rooks_decreasing}
by a standard identity \cite[eq.~(6.15)]{Graham_94}.

\medskip

2) At the other extreme, when $k=0$, paths are forbidden,
and we obtain the generating polynomial for permutations of $[n]$
with a weight $v_+$ for each cycle ascent (= excedance),
$v_-$ for each cycle descent (= anti-excedance),
$v_0$ for each fixed point, and $1+\alpha$ for each cycle.
For the corresponding exponential generating function,
see \reff{eq.first_multivariate.egf} specialized to $u=0$.
In particular, when $v_- = v_0 = v_+ = v$ we obtain
\be
   \sfL^{(\alpha)}(v,v,v)_{n,0}
   \;=\;
   (1+\alpha)^{\overline{n}} \, v^n
   \;.
\ee
And when $\alpha=0$ and $v_0 = v_-$, we obtain
\be
   \sfL^{(0)}(v_-,v_-,v_+)_{n,0}
   \;=\;
   \sum_{j=0}^{n-1} \euler{n}{j} \, v_+^j \, v_-^{n-j}
   \;,
\ee
where the Eulerian number $\euler{n}{j}$
is the number of permutations of $[n]$ with $j$ excedances (or $j$ descents).
More generally, these two specializations can be combined:
we obtain the generating polynomial for permutations
by number of excedances and number of cycles
\cite[eq.~(2.12)]{Sokal-Zeng_masterpoly}
\cite[Proposition~3.3]{Salas-Sokal_GKP}
\be
   \sfL^{(\alpha)}(v_-,v_-,v_+)_{n,0}
   \;=\;
   \sum_{j=0}^{n} \stirlingsubset{n}{j} \, (v_+ - v_-)^{n-j} \, v_-^j \, (1+\alpha)^{\overline{j}}
   \;.
\ee
And most generally, we can introduce additional fixed points, to obtain
\be
   \sfL^{(\alpha)}(v_-,v_0,v_+)_{n,0}
   \;=\;
   \sum_{i=0}^n  \binom{n}{i} \, [(1+\alpha)(v_0 - v_-)]^i \:
   \sfL^{(\alpha)}(v_-,v_-,v_+)_{n-i,0}
   \;.
\ee
\myendremark

\subsection{A second multivariate generalization}
   \label{subsec.multivariate.second}

But we can go farther:  instead of weighting {\em edges}\/ according
to whether they are increasing, decreasing or fixed points,
we can weight {\em pairs of successive edges}\/
according to whether they are peaks ($+-$), valleys ($-+$),
double ascents ($++$), double descents ($--$) or fixed points ($00$).
To define these concepts for a Laguerre digraph,
we first need to make a convention about {\em boundary conditions}\/
at the two ends of a path.
We will here use 0--0 boundary conditions:
that is, we extend the Laguerre digraph $G$ on the vertex set $[n]$
to a digraph $\widehat{G}$ on the vertex set $[n] \cup \{0\}$
by decreeing that any vertex $i \in [n]$
that has in-degree (resp.~out-degree) 0 in $G$
will receive an incoming (resp.~outgoing) edge
from (resp.~to) the vertex 0.
In this way each vertex $i \in [n]$
will have a unique predecessor $p(i) \in [n] \cup \{0\}$
and a unique successor $s(i) \in [n] \cup \{0\}$.
We then say that a vertex $i \in [n]$ is a
\begin{itemize}
   \item {\em peak}\/ (p) if $p(i) < i > s(i)$;
   \item {\em valley}\/ (v) if $p(i) > i < s(i)$;
   \item {\em double ascent}\/ (da) if $p(i) < i < s(i)$;
   \item {\em double descent}\/ (dd) if $p(i) > i > s(i)$;
   \item {\em fixed point}\/ (fp) if $p(i) = i = s(i)$.
\end{itemize}
(Note that ``fixed point'' is a synonym of ``loop''.)
When these concepts are applied to the cycles of a Laguerre digraph,
we obtain the usual \textbfit{cycle classification}
of indices in a permutation
as cycle peaks, cycle valleys, cycle double rises, cycle double falls
and fixed points \cite{Zeng_93,Sokal-Zeng_masterpoly}.
When applied to the paths of a Laguerre digraph,
we obtain the usual \textbfit{linear classification}
of indices in a permutation (written in word form)
as peaks, valleys, double ascents or double descents \cite[p.~45]{Stanley_12}.
Note that, because of the 0--0 boundary conditions,
an isolated vertex is always a peak,
the initial vertex of a path is always a peak or double ascent,
and the final vertex of a path is always a peak or double descent;
moreover, each path contains at least one peak.

We write $\pk(G), \val(G), \da(G), \dd(G), \fp(G)$
for the number of vertices $i \in [n]$ that are, respectively,
peaks, valleys, double ascents, double descents or fixed points.
We then introduce the
\textbfit{second multivariate Laguerre coefficient matrix}
\be
   \sfLhat^{(\alpha)}(\yp,\yv,\yda,\ydd,\yfp)_{n,k}
   \;\eqdef\;
   \sum_{G \in \LD_{n,k}}
       y_{\rm p}^{\pk(G)} y_{\rm v}^{\val(G)}
       y_{\rm da}^{\da(G)} y_{\rm dd}^{\dd(G)} y_{\rm fp}^{\fp(G)}
            \, (1+\alpha)^{\cyc(G)}
   \:.
 \label{def.coeffmat.gen.second}
\ee
This polynomial is homogeneous of degree $n$ in $\yp,\yv,\yda,\ydd,\yfp$.
In Section~\ref{subsec.egf.egf}
we will compute the bivariate exponential generating function
for this matrix (and indeed for a matrix that further generalizes it).

Because each path contains at least one peak,
we can, if we wish, remove a factor $y_{\rm p}^{k}$
and define a unit-lower-triangular matrix by
\be
   \sfLhat^{(\alpha)\flat}(\yp,\yv,\yda,\ydd,\yfp)_{n,k}
   \;\eqdef\;
   \sfLhat^{(\alpha)}(\yp,\yv,\yda,\ydd,\yfp)_{n,k} / y_{\rm p}^{k}
   \;.
 \label{def.coeffmat.gen.second.flat}
\ee
This polynomial is homogeneous of degree $n-k$ in $\yp,\yv,\yda,\ydd,\yfp$.

We can recover the first multivariate Laguerre coefficient matrix
by looking only at the second step of each pair:
specializing $\yp = \ydd = v_-$, $\yv = \yda = v_+$ and $\yfp = v_0$,
we obtain the previous counting with an extra weight $v_-$
associated to the final vertex of each path:
\be
   \sfLhat^{(\alpha)}(v_-,v_+,v_+,v_-,v_0)_{n,k}
   \;=\;
   \sfL^{(\alpha)}(v_-, v_0, v_+)_{n,k} \: v_-^k
 \label{eq.second-to-first.1}
\ee
or equivalently
\be
   \sfLhat^{(\alpha)\flat}(v_-,v_+,v_+,v_-,v_0)_{n,k}
   \;=\;
   \sfL^{(\alpha)}(v_-, v_0, v_+)_{n,k}
   \;.
 \label{eq.second-to-first.1a}
\ee
Alternatively, we can look only at the first step of each pair:
specializing $\yv = \ydd = v_-$, $\yp = \yda = v_+$ and $\yfp = v_0$,
we obtain the previous counting with an extra weight $v_+$
associated to the initial vertex of each path:
\be
   \sfLhat^{(\alpha)}(v_+,v_-,v_+,v_-,v_0)_{n,k}
   \;=\;
   \sfL^{(\alpha)}(v_-, v_0, v_+)_{n,k} \: v_+^k
 \label{eq.second-to-first.2}
\ee
or equivalently
\be
   \sfLhat^{(\alpha)\flat}(v_+,v_-,v_+,v_-,v_0)_{n,k}
   \;=\;
   \sfL^{(\alpha)}(v_-, v_0, v_+)_{n,k}
   \;.
 \label{eq.second-to-first.2a}
\ee

Once again we can introduce the row-generating polynomials:
\begin{eqnarray}
   \scrlhatna(x; \yp,\yv,\yda,\ydd,\yfp)
   & \eqdef &
   \sum_{k=0}^n \sfL^{(\alpha)}(\yp,\yv,\yda,\ydd,\yfp)_{n,k} \: x^k
       \label{def.multivariate.rowgen.scrl.second}
          \\[2mm]
   \scrlhatbarna(x; \yp,\yv,\yda,\ydd,\yfp)
   & \eqdef &
   \sum_{k=0}^n \sfL^{(\alpha)}(\yp,\yv,\yda,\ydd,\yfp)_{n,k} \: x^{n-k}
       \label{def.multivariate.rowgen.scrlbar.second}
\end{eqnarray}

\subsection{Total positivity}   \label{subsec.TP}

The main results of this paper concern the total positivity
of various matrices associated to the univariate and multivariate
Laguerre polynomials.
Recall first that a finite or infinite matrix of real numbers is called
{\em totally positive}\/ (TP) if all its minors are nonnegative,
and {\em totally positive of order~$r$} (TP${}_r$)
if all its minors of size $\le r$ are nonnegative.
Background information on totally positive matrices can be found
in \cite{Karlin_68,Gantmacher_02,Pinkus_10,Fallat_11};
they have application to many fields of pure and applied mathematics.\footnote{
   See \cite[footnote~4]{forests_totalpos} for many references.
}
In particular, it is known
\cite[Th\'eor\`eme~9]{Gantmakher_37} \cite[section~4.6]{Pinkus_10}
that an infinite Hankel matrix $(a_{i+j})_{i,j \ge 0}$
of real numbers is totally positive if and only if the underlying sequence
$(a_n)_{n \ge 0}$ is a Stieltjes moment sequence,
i.e.\ the moments of a positive measure on $[0,\infty)$.

But this is only the beginning of the story,
because we are here principally concerned,
not with sequences and matrices of real numbers,
but with sequences and matrices of polynomials
(with integer or real coefficients) in one or more indeterminates~$\bfx$:
they will typically be generating polynomials that enumerate
some combinatorial objects with respect to one or more statistics.
We equip the polynomial ring $\R[\bfx]$ with the coefficientwise
partial order:  that is, we say that $P$ is nonnegative
(and write $P \myge 0$)
in case $P$ is a polynomial with nonnegative coefficients.
We then say that a matrix with entries in $\R[\bfx]$ is
\textbfit{coefficientwise totally positive}
if all its minors are polynomials with nonnegative coefficients;
and analogously for coefficientwise total positivity of order~$r$.
We say that a sequence $\ba = (a_n)_{n \ge 0}$ with entries in $\R[\bfx]$
is \textbfit{coefficientwise Hankel-totally positive}
if its associated infinite Hankel
matrix
is coefficientwise totally positive;
and likewise for the version of order $r$.
Similar definitions apply to the formal-power-series ring $\R[[\bfx]]$.
Most generally, we can consider sequences and matrices
with values in an arbitrary partially ordered commutative ring
(a precise definition will be given in Section~\ref{subsec.totalpos.prelim});
total positivity and Hankel-total positivity
are then defined in the obvious way.
Coefficientwise Hankel-total positivity of a sequence of polynomials
$(P_n(\bfx))_{n \ge 0}$ {\em implies}\/ the pointwise Hankel-total positivity
(i.e.\ the Stieltjes moment property) for all $\bfx \ge 0$,
but it is vastly stronger.

For instance, it is known (see Section~\ref{sec.MOP})
that, for each $\alpha \ge -1$ and $x \ge 0$,
the sequence $(\scrlna(x))_{n \ge 0}$ of monic unsigned Laguerre polynomials
is a Stieltjes moment sequence.
So every minor of the Hankel matrix
$(\scrl_{i+j}^{(\alpha)}(x))_{i,j \ge 0}$
is a polynomial in $x$ and $\alpha$ that is nonnegative
whenever $\alpha \ge -1$ and $x \ge 0$.
But much more turns out to be true:
every minor of the Hankel matrix 
$(\scrl_{i+j}^{(\alpha)}(x))_{i,j \ge 0}$
is in~fact a polynomial in $x$ and $\alpha$ that is
{\em coefficientwise}\/ nonnegative in $x$ and $\lambda \eqdef 1+\alpha$.
That is, the sequence $(\scrl_n^{(-1+\lambda)}(x))_{n \ge 0}$
is {\em coefficientwise}\/ Hankel-totally positive in $x$ and $\lambda$:
this will be our first main result (Theorem~\ref{thm.univariate}(c)).
We will then extend this result to the multivariate polynomials
introduced in Sections~\ref{subsec.multivariate.first}
and \ref{subsec.multivariate.second}.

Our proofs of total positivity and Hankel-total positivity
will be based on the method of
{\em production matrices}\/ \cite{Deutsch_05,Deutsch_09}.
We will review this theory
in Sections~\ref{subsec.production} and \ref{subsec.totalpos.prodmat},
so now we state only the bare-bones definitions.
Let $P = (p_{ij})_{i,j \ge 0}$ be an infinite matrix
with entries in a commutative ring~$R$;
we assume that $P$ is either row-finite
(i.e.\ has only finitely many nonzero entries in each row)
or column-finite.
Now define an infinite matrix $A = (a_{nk})_{n,k \ge 0}$ by
\be
   a_{nk}  \;=\;  (P^n)_{0k}
   \;.
\ee
We call $P$ the \textbfit{production matrix}
and $A$ the \textbfit{output matrix}, and we write $A = \scro(P)$.
The two key facts here are the following \cite{Sokal_totalpos}:
if $R$ is a partially ordered commutative ring,
$P$ is totally positive of order $r$,
and $a$ is a nonnegative element of~$R$,
then $\scro(P+aI)$ is totally positive of order $r$
and the zeroth column of $\scro(P+aI)$ is Hankel-totally positive of order $r$.
See Section~\ref{subsec.totalpos.prodmat} for precise statements and proofs.

\subsection{Statement of main results: Univariate case}

Our first main result is the following:

\begin{theorem}[Total positivity of the univariate Laguerre polynomials]
   \label{thm.univariate}
\hfill\break\noindent
\vspace*{-6mm}
\begin{itemize}
   \item[(a)] The Laguerre coefficient matrix $\sfL^{(-1+\lambda)}$
is totally positive in the ring $\Z[\lambda]$
equipped with the coefficientwise order.
   \item[(b)] The binomial row-generating matrix $\sfL^{(-1+\lambda)} \, B_x$
is totally positive in the ring $\Z[x,\lambda]$
equipped with the coefficientwise order.
   \item[(c)] The sequence of monic unsigned Laguerre polynomials
$(\scrl_n^{(-1+\lambda)}(x))_{n \ge 0}$
is Hankel-totally positive in the ring $\Z[x,\lambda]$
equipped with the coefficientwise order.
   \item[(d)] The sequence of reversed monic unsigned Laguerre polynomials
$(\scrlbar_n^{(-1+\lambda)}(x))_{n \ge 0}$
is Hankel-totally positive in the ring $\Z[x,\lambda]$
equipped with the coefficientwise order.
\end{itemize}
\end{theorem}

Theorem~\ref{thm.univariate}(c,d) proves a conjecture made a few years ago
by Sylvie Corteel and one of us \cite{Corteel-Sokal_conj_Laguerre}.
Theorem~\ref{thm.univariate} was also proven recently
by Zhu \cite[Proposition~4.14]{Zhu_21} \cite[Proposition~5.11]{Zhu_22};
his methods are very similar to ours.\footnote{
   Zhu \cite[Proposition~4.14]{Zhu_21} \cite[Proposition~5.11]{Zhu_22}
   states a slightly weaker form of Theorem~\ref{thm.univariate}
   in which $\lambda$ is a nonnegative real number
   rather than an indeterminate.
   But his methods actually prove the stronger result claimed here.
}

The result in part~(a) is fairly easy, and we will give two proofs
using different methods:
one direct proof (Section~\ref{subsec.univariate_coeff_TP}),
and one proof using production matrices.
Part~(b) is then an immediate consequence of part~(a)
together with the coefficientwise total positivity of the binomial matrix $B_x$
(Lemma~\ref{lemma.binomialmatrix.TP}).
Also, part~(d) is trivially equivalent to part~(c) by virtue of the relation
$\scrlbar_n^{(-1+\lambda)}(x) = x^n \, \scrl_n^{(-1+\lambda)}(1/x)$.
So our main effort will be directed to proving part~(c):
we will do this by constructing the production matrix
for the binomial row-generating matrix $\sfL^{(-1+\lambda)} \, B_x$
and then proving its total positivity.
More precisely, we will prove the following:

\begin{proposition}[Production matrices for the univariate Laguerre polynomials]
   \label{prop.prodmat.univariate.NEW}
\hfill\break\noindent
\vspace*{-6mm}
\begin{itemize}
   \item[(a)] The production matrix of the Laguerre coefficient array
$\sfL^{(\alpha)}$ is the tridiagonal unit-lower-Hessenberg matrix
$P^\circ = (p^\circ_{ij})_{i,j \ge 0}$ defined by
\begin{subeqnarray}
   p^\circ_{n,n+1}  & = &   1   \\[1mm]
   p^\circ_{n,n}    & = &   2n+1+\alpha   \\[1mm]
   p^\circ_{n,n-1}  & = &   n(n+\alpha)   \\[1mm]
   p^\circ_{n,k}    & = &   0 \qquad\textrm{if $k < n-1$ or $k > n+1$}
 \label{eq.prop.prodmat.univariate.NEW.1}
\end{subeqnarray} 
   \item[(b)] The production matrix of the binomial row-generating matrix
$\sfL^{(\alpha)} B_x$ is the quadridiagonal unit-lower-Hessenberg matrix
$P = (p_{ij})_{i,j \ge 0}$ defined by
\begin{subeqnarray}
   p_{n,n+1}  & = &   1   \\[1mm]
   p_{n,n}    & = &   (2n+1+\alpha) \,+\, x   \\[1mm]
   p_{n,n-1}  & = &   n(n+\alpha) \,+\, 2nx \\[1mm]
   p_{n,n-2}  & = &   n(n-1) x  \\[1mm]
   p_{n,k}    & = &   0 \qquad\textrm{if $k < n-2$ or $k > n+1$}
 \label{eq.prop.prodmat.univariate.NEW.2}
\end{subeqnarray} 
\end{itemize}
\end{proposition}

{\bf Remarks.}
1.  When $\alpha \in \{-1,0,1\}$, the production matrix
\reff{eq.prop.prodmat.univariate.NEW.2}
arises from a 2-branched S-fraction \cite{latpath_SRTR}
--- more precisely, as the production matrix for the
generalized $m$-Stieltjes--Rogers polynomials of type $j$
with $m=2$ and $j \in \{0,1,2\}$.
See Appendix~\ref{subsec.mSR.gen} for the theory of these
generalized $m$-Stieltjes--Rogers polynomials,
and Appendix~\ref{subsec.mSR.univariate} for the application to 
\reff{eq.prop.prodmat.univariate.NEW.2}.

2.  It will follow easily from our general theory
(Lemma~\ref{lemma.production.AB})
that the production matrices of
Proposition~\ref{prop.prodmat.univariate.NEW}(a,b)
are related by $P = B_x^{-1} P^\circ B_x$.
What is far from obvious is why $P$ is quadridiagonal,
since it is {\em not}\/ in general true that
$B_x^{-1} T B_x$ is quadridiagonal whenever $T$ is tridiagonal.
In Appendix~\ref{app.banded} we explain why $P$ is quadridiagonal
in the present case, by answering the more general question:
Which lower-Hessenberg matrices $P$ have the property that
$B_\xi^{-1} P B_\xi$ is $(r,1)$-banded?
\myendremark

\begin{proposition}[Total positivity of the univariate production matrices]
   \label{prop.prodmat.TP.univariate.NEW}
\hfill\break\noindent
\vspace*{-6mm}
\begin{itemize}
   \item[(a)] The matrix $P^\circ = (p^\circ_{ij})_{i,j \ge 0}$
defined by \reff{eq.prop.prodmat.univariate.NEW.1},
with the change of variable $\alpha = -1 + \lambda$,
is totally positive in the ring $\Z[\lambda]$
equipped with the coefficientwise order.
   \item[(b)] The matrix $P = (p_{ij})_{i,j \ge 0}$
defined by \reff{eq.prop.prodmat.univariate.NEW.2},
with the change of variable $\alpha = -1 + \lambda$,
is totally positive in the ring $\Z[x,\lambda]$
equipped with the coefficientwise order.
\end{itemize}
\end{proposition}

Combining Propositions~\ref{prop.prodmat.univariate.NEW}
and \ref{prop.prodmat.TP.univariate.NEW}
with the general theory of production matrices
proves Theorem~\ref{thm.univariate}.

We will prove Proposition~\ref{prop.prodmat.univariate.NEW}
in Section~\ref{sec.prodmat.1var},
and Proposition~\ref{prop.prodmat.TP.univariate.NEW}
in Section~\ref{sec.prodmat.TP.univariate}.


\subsection{Statement of main results: Multivariate case}

We now state the corresponding results
for the first and second multivariate Laguerre polynomials
that were introduced in Sections~\ref{subsec.multivariate.first}
and \ref{subsec.multivariate.second}, respectively.
We begin by considering the second multivariate Laguerre coefficient matrix,
for which it is convenient to use the form
$\sfLhat^{(\alpha)\flat}(\yp,\yv,\yda,\ydd,\yfp)$
that was defined in \reff{def.coeffmat.gen.second.flat}.
Our main result, generalizing Theorem~\ref{thm.univariate}, is the following:

\enlargethispage*{3\baselineskip}

\begin{theorem}[Total positivity of the second multivariate Laguerre polynomials]
   \label{thm.multivariate}
\hfill\break\noindent
Let $\lambda,\yp,\yv,\yda,\ydd,\yfp$ be elements
of a partially ordered commutative ring $R$
that satisfy $\lambda \ge 0$,
$\lambda\yfp \ge \lambda\yp$, $\yp \ge 0$, $\yv \ge 0$
and $\yda+\ydd \ge \yp + \yv$;
and let $x$ be an indeterminate.
Then:
\begin{itemize}
   \item[(a)] The second multivariate Laguerre coefficient array
$\sfLhat^{(-1+\lambda)\flat}(\yp,\yv,\yda,\ydd,\yfp)$
is totally positive in the ring $R$.
   \item[(b)] The binomial row-generating matrix
$\sfLhat^{(-1+\lambda)\flat}(\yp,\yv,\yda,\ydd,\yfp) \, B_x$
is totally positive in the ring $R[x]$ equipped with the coefficientwise order.
   \item[(c)] The sequence of row-generating polynomials
$\widehat{\scrl}_n^{(-1+\lambda)}(x; \yp,\yv,\yda,\ydd,\yfp)$
defined by \reff{def.multivariate.rowgen.scrl.second}
is Hankel-totally positive in the ring $R[x]$
equipped with the coefficientwise order.
   \item[(d)] The sequence of reversed row-generating polynomials
$\overline{\widehat{\scrl}}_n^{(-1+\lambda)}(x; \yp,\yv,\yda,\ydd,\yfp)$
defined by \reff{def.multivariate.rowgen.scrlbar.second}
is Hankel-totally positive in the ring $R[x]$
equipped with the coefficientwise order.
\end{itemize}
\nopagebreak
These statements are also true when the roles of
$\yp$ and $\yv$ are interchanged.
\end{theorem}

As in the univariate case,
part~(b) is an immediate consequence of part~(a),
and part~(d) is trivially equivalent to part~(c).
We therefore concentrate on proving (a) and (c).

As before, we will prove Theorem~\ref{thm.multivariate}
by exhibiting the production matrices and then proving their total positivity.

\begin{proposition}[Production matrices for the multivariate Laguerre polynomials]
   \label{prop.prodmat.multivariate.NEW}
\hfill\break\noindent
\vspace*{-6mm}
\begin{itemize}
   \item[(a)] The production matrix of the multivariate Laguerre
coefficient array
\hfill\break
$\sfLhat^{(\alpha)\flat}(\yp,\yv,\yda,\ydd,\yfp)$
is the tridiagonal unit-lower-Hessenberg matrix
$P^{\circ\flat} = (p^{\circ\flat}_{ij})_{i,j \ge 0}$ defined by
\begin{subeqnarray}
   p^{\circ\flat}_{n,n+1}  & = &   1   \\[1mm]
   p^{\circ\flat}_{n,n}    & = &   (1+\alpha)\yfp \,+\, n(\yda+\ydd)   \\[1mm]
   p^{\circ\flat}_{n,n-1}  & = &   n(n+\alpha)\yp\yv   \\[1mm]
   p^{\circ\flat}_{n,k}    & = &   0 \qquad\textrm{if $k < n-1$ or $k > n+1$}
 \label{eq.prop.multivariate_second.prodmat.flat.BIS0}
\end{subeqnarray}
   \item[(b)] The production matrix of the binomial row-generating matrix
$\sfLhat^{(\alpha)\flat}(\yp,\yv,\yda,\ydd,\yfp) B_x$
is the quadridiagonal unit-lower-Hessenberg matrix
$P^\flat = (p^\flat_{ij})_{i,j \ge 0}$ defined by
\begin{subeqnarray}
   p^\flat_{n,n+1}  & = &   1   \\[1mm]
   p^\flat_{n,n}    & = &   (1+\alpha)\yfp \,+\, n(\yda+\ydd) \,+\, x \\[1mm]
   p^\flat_{n,n-1}  & = &   n(n+\alpha)\yp\yv \,+\, n(\yda+\ydd) x   \\[1mm]
   p^\flat_{n,n-2}  & = &   n(n-1) \yp\yv x   \\[1mm]
   p^\flat_{n,k}    & = &   0 \qquad\textrm{if $k < n-2$ or $k > n+1$}
 \label{eq.prop.multivariate_second.rowgen.prodmat.flat.BIS0}
\end{subeqnarray}
\end{itemize}
\end{proposition}

\begin{proposition}[Total positivity of the multivariate production matrices]
   \label{prop.prodmat.TP.multivariate.NEW}
\hfill\break\noindent
\vspace*{-6mm}
\begin{itemize}
   \item[(a)] The matrix $P^{\circ\flat} = (p^{\circ\flat}_{ij})_{i,j \ge 0}$
defined by \reff{eq.prop.multivariate_second.prodmat.flat.BIS0},
with the variables substituted to elements
of a partially ordered commutative ring $R$
that satisfy $\alpha \ge -1$,
${(1+\alpha)\yfp \ge (1+\alpha)\yp}$, $\yp \ge 0$, $\yv \ge 0$
and $\yda+\ydd \ge \yp + \yv$,
is totally positive in the ring $R$.
\item[(b)] The matrix $P^\flat = (p^\flat_{ij})_{i,j \ge 0}$
defined by \reff{eq.prop.multivariate_second.rowgen.prodmat.flat.BIS0},
with the variables substituted to elements
of a partially ordered commutative ring $R$
that satisfy $\alpha \ge -1$,
${(1+\alpha)\yfp = (1+\alpha)\yp}$, $\yp \ge 0$, $\yv \ge 0$,
$\yda+\ydd \ge \yp + \yv$ and $x \ge 0$,
is totally positive in the ring $R$.
\end{itemize}
These statements are also true when the roles of
$\yp$ and $\yv$ are interchanged.
\end{proposition}

We will prove Proposition~\ref{prop.prodmat.multivariate.NEW}
in Section~\ref{subsec.egf.prodmat},
and Proposition~\ref{prop.prodmat.TP.multivariate.NEW}
in Section~\ref{sec.prodmat.TP.multivariate}.
Here the most difficult part is
Proposition~\ref{prop.prodmat.TP.multivariate.NEW}(b):
namely, proving the total positivity of the quadridiagonal matrix
\reff{eq.prop.multivariate_second.rowgen.prodmat.flat.BIS0}
under the specified conditions.
We will do this by proving the total positivity of a
much more general quadridiagonal matrix
(Theorem~\ref{thm.prodmat.TP.bis.gen.new2});
and see also Appendix~\ref{sec.prodmat.TP.multivariate.quadridiagonal.OLD}
(Theorem~\ref{thm.quadmat2})
for a variant of this result.

\bigskip

{\bf Remarks.}
1.  The zeroth column of the multivariate Laguerre coefficient array
$\sfLhat^{(\alpha)\flat}(\yp,\yv,\yda,\ydd,\yfp)$
is the sequence of generating polynomials for permutations
weighted according to the cycle classification,
and the tridiagonal production matrix
\reff{eq.prop.multivariate_second.prodmat.flat.BIS0}
is the production matrix for the corresponding J-fraction,
namely, \cite[Theorem~2.4]{Sokal-Zeng_masterpoly}
specialized to $x_1 = u_1 = \yp$, $x_2 = u_2 = \ydd$,
$y_1 = v_1 = \yv$, $y_2 = v_2 = \yda$, $w_n = \yfp$, $\lambda = 1 + \alpha$.

2.  Note that in part~(a) we require only the {\em inequality}\/
\break
${(1+\alpha)\yfp \ge (1+\alpha)\yp}$,
while in part~(b) we require the {\em equality}\/
${(1+\alpha)\yfp = (1+\alpha)\yp}$.
However, this will not preclude us from proving
Theorem~\ref{thm.multivariate} assuming only the inequality,
as our general theory
(Corollaries~\ref{cor.iteration.homo.aI+bP} and
 \ref{thm.iteration2bis.aI+bP})
entitles us to use the production matrix $P^\flat + aI$
with $a = \lambda (\yfp-\yp) \ge 0$.
\myendremark

\medskip

We can now recover the first multivariate Laguerre coefficient matrix
$\sfL^{(\alpha)}(v_-, v_0, v_+)$ defined in \reff{def.coeffmat.gen},
and its corresponding row-generating polynomials
$\scrlna(x; v_-, v_0, v_+)$ defined in \reff{def.multivariate.rowgen.scrl},
by specializing either
\be
  \hbox{$\yp = \ydd = v_-$, $\;\yv = \yda = v_+$, $\;\yfp = v_0\;$
as in \reff{eq.second-to-first.1}/\reff{eq.second-to-first.1a}}
  \hphantom{\;.}
 \label{eq.first.specialization.1}
\ee
or
\be
   \hbox{$\yv = \ydd = v_-$, $\;\yp = \yda = v_+$, $\;\yfp = v_0\;$
as in \reff{eq.second-to-first.2}/\reff{eq.second-to-first.2a}}
  \;.
 \label{eq.first.specialization.2}
\ee
By applying these specializations to Theorem~\ref{thm.multivariate}, we obtain:

\begin{corollary}[Total positivity of the first multivariate Laguerre polynomials]
   \label{cor.multivariate}
\hfill\break\noindent
\vspace*{-6mm}
\begin{itemize}
   \item[(a)] The first multivariate Laguerre coefficient matrix
      $\sfL^{(-1+\lambda)}(v_-, v_0, v_+)$,
      defined in \reff{def.coeffmat.gen},
      is totally positive when the variables are substituted to elements
      of a partially ordered commutative ring $R$
      that satisfy $\lambda \ge 0$, $v_- \ge 0$, $v_+ \ge 0$,
      and either $v_0 \ge v_-$ or $v_0 \ge v_+$.

      \quad
      Equivalently, the matrices
      $\sfL^{(-1+\lambda)}(v_-, v_- + w, v_+)$
      and $\sfL^{(-1+\lambda)}(v_-, v_+ + w, v_+)$
      are totally positive in the ring $\Z[v_-,v_+,w,\lambda]$
      equipped with the coefficientwise order.
   \item[(b)]  The sequences
      $\bigl( \scrlna(x; v_-, v_-, v_+) \bigr)_{n \ge 0}$
      and $\bigl( \scrlna(x; v_-, v_+, v_+) \bigr)_{n \ge 0}$
      of row-generating polynomials,
      defined in \reff{def.multivariate.rowgen.scrl},
      are Hankel-totally positive,
      in the ring $R[x]$ equipped with the coefficientwise order,
      when the variables are substituted to elements
      of a partially ordered commutative ring $R$
      that satisfy $\lambda \ge 0$, $v_- \ge 0$, $v_+ \ge 0$,
      and either $v_0 \ge v_-$ or $v_0 \ge v_+$.

      \quad
      Equivalently, the sequences
      $\bigl( \scrlna(x; v_-, {v_- +w}, v_+) \bigr)_{n \ge 0}$
      and
      \break
      $\bigl( \scrlna(x; v_-, v_+ +w, v_+) \bigr)_{n \ge 0}$
      are Hankel-totally positive in the ring $\Z[x, v_-,v_+,w,\lambda]$
      equipped with the coefficientwise order.
\end{itemize}
\end{corollary}

The multivariate case
(Theorem~\ref{thm.multivariate} and
 Propositions~\ref{prop.prodmat.multivariate.NEW}
 and \ref{prop.prodmat.TP.multivariate.NEW})
obviously subsumes the univariate case
(Theorem~\ref{thm.univariate} and
 Propositions~\ref{prop.prodmat.univariate.NEW}
 and \ref{prop.prodmat.TP.univariate.NEW}),
so it is in principle redundant to consider the latter.
But since the proofs in the univariate case are quite a bit simpler,
and since this case exhibits some special features that do not carry over
to the fully multivariate case,
we think it useful to present the univariate case first.


%

%
%

\subsection{Plan of this paper}

Although the present paper is a follow-up to our papers
\cite{latpath_SRTR,latpath_lah},
we have endeavored, for the convenience of the reader,
to make it as self-contained as possible.
We have therefore begun, in Section~\ref{sec.prelim},
with a brief review of the key definitions and results from
\cite{latpath_SRTR,latpath_lah}
(plus a few other things) that will be needed in the sequel.
We then proceed as follows:
In Section~\ref{sec.prodmat.1var}
we determine the production matrices for the univariate
Laguerre coefficient matrix and binomial row-generating matrix.
In Section~\ref{sec.egf} we compute the exponential generating function
for the multivariate Laguerre coefficient matrix
\reff{def.coeffmat.gen.second}
and employ this to determine the corresponding production matrices,
using the theory of exponential Riordan arrays.
In Sections~\ref{sec.prodmat.TP.univariate}
and \ref{sec.prodmat.TP.multivariate}
we prove
the total positivity of the univariate and multivariate production matrices,
respectively.
In particular, in Section~\ref{subsec.prodmat.TP.multivariate.quadridiagonal}
(Theorem~\ref{thm.prodmat.TP.bis.gen.new2})
we prove the coefficientwise total positivity
of a rather general class of quadridiagonal matrices;
we think that this result, and its method of proof,
may be of some independent interest.
In Section~\ref{sec.MOP} we explain the unexpected connection with
multiple orthogonal polynomials;
in particular, we explain how we were led to guess
the production matrices \reff{eq.prop.prodmat.univariate.NEW.2}
and \reff{eq.prop.multivariate_second.rowgen.prodmat.flat.BIS0}.

In Appendix~\ref{app.modified}
we review the theory of classical and branched S-fractions,
and then introduce the generalized $m$-Stieltjes--Rogers polynomials
and modified $m$-Stieltjes--Rogers polynomials of type $j \ge 0$;
this expands and supersedes some of the discussion in
\cite[sections~5, 7 and 9]{latpath_SRTR}.
We also apply this theory to the
univariate Laguerre production matrix \reff{eq.prop.prodmat.univariate.NEW.2}.
In Appendix~\ref{app.banded} we answer the question:
Which lower-Hessenberg matrices $P$ have the property that
$B_\xi^{-1} P B_\xi$ is $(r,1)$-banded?
In Appendix~\ref{sec.prodmat.TP.multivariate.quadridiagonal.OLD}
we prove the total positivity of a class of quadridiagonal matrices,
giving a variant of what was done in
Section~\ref{subsec.prodmat.TP.multivariate.quadridiagonal}.

The first proof of the production matrices in this paper
was found in 2019 by one of us (M.P.)\ using
a bijection from Laguerre digraphs to labeled 2-\L{}ukasiewicz paths.
We hope to revisit this powerful method in a subsequent paper.

\section{Preliminaries}   \label{sec.prelim}

In this section we review some definitions and results from
\cite{latpath_SRTR,latpath_lah,forests_totalpos,Chen-Sokal_trees_totalpos,Sokal_totalpos}
that will be needed in what follows.
After a brief introduction to total positivity
in a partially ordered commutative ring (Section~\ref{subsec.totalpos.prelim}),
we provide a pr\'ecis of the theory of production matrices
(Section~\ref{subsec.production})
and its application to total positivity
(Section~\ref{subsec.totalpos.prodmat});
these latter results form the theoretical foundation for our work.
Then
we introduce the concept of binomial row-generating matrices
(Section~\ref{subsec.rowgen}).
Finally, we review the theory of exponential Riordan arrays
(Section~\ref{sec.exp_riordan});
this theory will be our main technical tool.

For the benefit of readers who are familiar with the papers
\cite{latpath_SRTR,latpath_lah,forests_totalpos,Chen-Sokal_trees_totalpos},
let us mention two things that are new here (and crucial for the present paper):
the tridiagonal comparison theorem (Proposition~\ref{prop.comparison}),
and adding a multiple of the identity to the production matrix
(Corollaries~\ref{cor.iteration.homo.aI+bP} and \ref{thm.iteration2bis.aI+bP}).

\subsection{Partially ordered commutative rings and total positivity}
   \label{subsec.totalpos.prelim}

In this paper all rings will be assumed to have an identity element 1
and to be nontrivial ($1 \ne 0$).

A \textbfit{partially ordered commutative ring} is a pair $(R,\scrp)$ where
$R$ is a commutative ring and $\scrp$ is a subset of $R$ satisfying
\begin{itemize}
   \item[(a)]  $0,1 \in \scrp$.
   \item[(b)]  If $a,b \in \scrp$, then $a+b \in \scrp$ and $ab \in \scrp$.
   \item[(c)]  $\scrp \cap (-\scrp) = \{0\}$.
\end{itemize}
We call $\scrp$ the {\em nonnegative elements}\/ of $R$,
and we define a partial order on $R$ (compatible with the ring structure)
by writing $a \le b$ as a synonym for $b-a \in \scrp$.
Please note that, unlike the practice in real algebraic geometry
\cite{Brumfiel_79,Lam_84,Prestel_01,Marshall_08},
we do {\em not}\/ assume here that squares are nonnegative;
indeed, this property fails completely for our prototypical example,
the ring of polynomials with the coefficientwise order,
since $(1-x)^2 = 1-2x+x^2$ is not coefficientwise nonnegative.

Now let $(R,\scrp)$ be a partially ordered commutative ring
and let $\bfx = \{x_i\}_{i \in I}$ be a collection of indeterminates.
In the polynomial ring $R[\bfx]$ and the formal-power-series ring $R[[\bfx]]$,
let $\scrp[\bfx]$ and $\scrp[[\bfx]]$ be the subsets
consisting of polynomials (resp.\ series) with nonnegative coefficients.
Then $(R[\bfx],\scrp[\bfx])$ and $(R[[\bfx]],\scrp[[\bfx]])$
are partially ordered commutative rings;
we refer to this as the \textbfit{coefficientwise order}
on $R[\bfx]$ and $R[[\bfx]]$.

A finite or infinite matrix with entries in a
partially ordered commutative ring
is called \textbfit{totally positive} (TP) if all its minors are nonnegative;
it is called \textbfit{totally positive of order~$\bm{r}$} (TP${}_r$)
if all its minors of size $\le r$ are nonnegative.
It follows immediately from the Cauchy--Binet formula that
the product of two TP (resp.\ TP${}_r$) matrices is TP
(resp.\ TP${}_r$).\footnote{
   For infinite matrices, we need some condition to ensure that
   the product is well-defined.
   For instance, the product $AB$ is well-defined whenever
   $A$ is row-finite (i.e.\ has only finitely many nonzero entries in each row)
   or $B$ is column-finite.
}
This fact is so fundamental to the theory of total positivity
that we shall henceforth use it without comment.

We say that a sequence $\ba = (a_n)_{n \ge 0}$
with entries in a partially ordered commutative ring
is \textbfit{Hankel-totally positive} 
(resp.\ \textbfit{Hankel-totally positive of order~$\bm{r}$})
if its associated infinite Hankel matrix
$H_\infty(\ba) = (a_{i+j})_{i,j \ge 0}$
is TP (resp.\ TP${}_r$).
We say that $\ba$
is \textbfit{Toeplitz-totally positive} 
(resp.\ \textbfit{Toeplitz-totally positive of order~$\bm{r}$})
if its associated infinite Toeplitz matrix
$T_\infty(\ba) = (a_{i-j})_{i,j \ge 0}$
(where $a_n \eqdef 0$ for $n < 0$)
is TP (resp.\ TP${}_r$).\footnote{
   When $R = \R$, Toeplitz-totally positive sequences are traditionally called
   {\em P\'olya frequency sequences}\/ (PF),
   and Toeplitz-totally positive sequences of order $r$
   are called {\em P\'olya frequency sequences of order $r$}\/ (PF${}_r$).
   See \cite[chapter~8]{Karlin_68} for a detailed treatment.
}

We will need a few easy facts about the total positivity of special matrices:

\begin{lemma}[Bidiagonal matrices]
  \label{lemma.bidiagonal}
Let $A$ be a matrix with entries in a partially ordered commutative ring,
with the property that all its nonzero entries belong to two consecutive
diagonals.
Then $A$ is totally positive if and only if all its entries are nonnegative.
\end{lemma}

\proof
The nonnegativity of the entries (i.e.\ TP${}_1$)
is obviously a necessary condition for TP.
Conversely, for a matrix of this type it is easy to see that
every nonzero minor is simply a product of some entries.
\qed

\begin{lemma}[Binomial matrix]
   \label{lemma.binomialmatrix.TP}
In the ring $\Z$, the binomial matrix
$B = \bigl( {\textstyle \binom{n}{k}} \bigr) _{n,k \ge 0}$
is totally positive.
More generally, the weighted binomial matrix
\be
   B_{x,y}
   \;=\;
   \bigl( x^{n-k} y^k {\textstyle \binom{n}{k}} \bigr) _{\! n,k \ge 0}
 \label{def.Bxy}
\ee
is totally positive in the ring $\Z[x,y]$ equipped with the
coefficientwise order.
\end{lemma}

We also write $B_x$ as a shorthand for $B_{x,1}$.

\proofof{Lemma~\ref{lemma.binomialmatrix.TP}}
It is well known that the binomial matrix $B$ is totally positive,
and this can be proven by a variety of methods:
e.g.\ using production matrices
\cite[pp.~136--137, Example~6.1]{Karlin_68}
\cite[pp.~108--109]{Pinkus_10},
by diagonal similarity to a totally positive Toeplitz matrix
\cite[p.~109]{Pinkus_10},
by exponentiation of a nonnegative lower-subdiagonal matrix
\cite[p.~63]{Fallat_11},
or by an application of the Lindstr\"om--Gessel--Viennot lemma
\cite[p.~24]{Fomin_00}.

Then $B_{x,y} = D B D'$ where $D = \diag\bigl( (x^n)_{n \ge 0} \bigr)$
and $D' = \diag\bigl( (x^{-k} y^k)_{k \ge 0} \bigr)$.
By Cauchy--Binet, $B_{x,y}$ is totally positive in the ring $\Z[x,x^{-1},y]$
equipped with the coefficientwise order.
But because $B$ is lower-triangular,
the elements of $B_{x,y}$ actually lie in the subring $\Z[x,y]$.
\qed

\noindent
See also Example~\ref{exam.binomial.matrix.TP} below
for an {\em ab initio}\/ proof of Lemma~\ref{lemma.binomialmatrix.TP}
using production matrices.

\begin{lemma}[Introducing indeterminates]
   \label{lemma.diagmult.TP}
Let $A = (a_{ij})_{i,j \ge 0}$ be a lower-triangular matrix
with entries in a partially ordered commutative ring $R$,
and let $\bx = (x_i)_{i \ge 1}$.
Define the lower-triangular matrix $B = (b_{ij})_{i,j \ge 0}$ by
\be
   b_{ij}  \;=\;  x_{j+1} x_{j+2} \cdots x_i \, a_{ij}
   \;.
\ee
Then:
\begin{itemize}
   \item[(a)] If $A$ is TP${}_r$ and $\bx$ are indeterminates,
      then $B$ is TP${}_r$ in the ring $R[\bx]$ equipped with
      the coefficientwise order.
   \item[(b)] If $A$ is TP${}_r$ and $\bx$ are nonnegative elements of $R$,
      then $B$ is TP${}_r$ in the ring $R$.
\end{itemize}
\end{lemma}

\proof
(a) Let $\bx = (x_i)_{i \ge 1}$ be commuting indeterminates,
and let us work in the ring $R[\bx,\bx^{-1}]$
equipped with the coefficientwise order.
Let $D = \diag(1,\, x_1,\, x_1 x_2,\, \ldots)$.
Then $D$ is invertible, and
$D^{-1} =
 \diag(1,\, x_1^{-1},\, x_1^{-1} x_2^{-1},\, \ldots)$
has nonnegative elements.
It follows that $B = D A D^{-1}$ is TP${}_r$ in the ring $R[\bx,\bx^{-1}]$
equipped with the coefficientwise order.
But the matrix elements $b_{ij}$
actually belong to the subring $R[\bx] \subseteq R[\bx,\bx^{-1}]$.
So $B$ is TP${}_r$ in the ring $R[\bx]$
equipped with the coefficientwise order.

(b) follows from (a) by specializing indeterminates.
\qed

\bigskip

Finally, we need some special facts about the total positivity
of tridiagonal matrices.
We recall that a {\em contiguous principal minor}\/
of a matrix $A$ is a minor $\det A_{II}$ for $I = \{r,{r+1},{r+2},\ldots,s\}$
for some $r \le s$.

\begin{lemma}
   \label{lemma.tridiagonal}
Every nonzero $k \times k$ minor of a tridiagonal matrix
is a product of off-diagonal elements and contiguous principal minors,
with sizes adding up to $k$.
\end{lemma}

\noindent
The standard proof for real matrices \cite[p.~98]{Pinkus_10}
is valid in any commutative ring.

\begin{corollary}[Total positivity of tridiagonal matrices]
   \label{cor.tridiagonal}
Let $A$ be a tridiagonal matrix with entries in a
partially ordered commutative ring.
Then $A$ is totally positive of order $r$
if and only~if all its off-diagonal elements
and all its contiguous principal minors of size $\le r$ are nonnegative.
\end{corollary}

\noindent
In other words, the well-known \cite[Theorem~4.3]{Pinkus_10}
criterion for total positivity of tridiagonal matrices with real entries
extends without change to matrices with entries in
an arbitrary partially ordered commutative ring.

This has the following important consequence:

\begin{proposition}[Tridiagonal comparison theorem, weak form]
   \label{prop.comparison}
Let $A$ and $D$ be matrices with entries in a
partially ordered commutative ring,
with $A$ being tridiagonal and $D$ being diagonal.
If $A$ is totally positive of order $r$ and $D$ is nonnegative,
then $A+D$ is totally positive of order $r$.
\end{proposition}

\proof
We apply Corollary~\ref{cor.tridiagonal}.
The off-diagonal elements of $A+D$ are the same as those of $A$.
The principal minors of $A+D$ are
\be
   \det\, (A+D)_{II}
   \;=\;
   \sum_{J \subseteq I} (\det A_{JJ}) \prod_{i \in I \setminus J} d_i
   \;,
\ee
where $d_i$ are the diagonal elements of $D$.
\qed

{\bf Remarks.}
1.  There is a stronger form of the tridiagonal comparison theorem
\cite{Sokal_totalpos} \cite[Proposition~3.1]{Zhu_21b}
in which we can increase the diagonal elements
{\em and}\/ decrease the off-diagonal elements
(while keeping them nonnegative).
But we will not need this stronger result here.

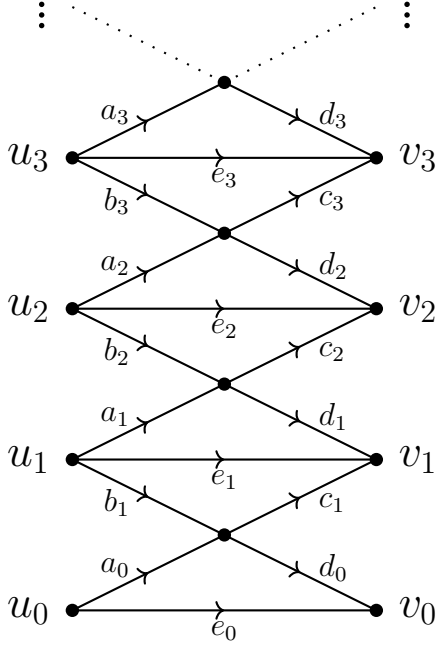
\begin{figure}[t]
\vspace{3mm}
\centering
\begin{tikzpicture}

\foreach \i in {0,...,3}
{
\pgfmathtruncatemacro{\rs}{2*\i+1};
\pgfmathtruncatemacro{\rt}{2*\i};
\path[draw=black,thick,->-=0.5] (0,\rt) to (2,\rs);
\path[draw=black,thick,->-=0.5] (0,\rt) to (4,\rt);
\path[draw=black,thick,->-=0.5] (2,\rs) to (4,\rt);
\node[circle,fill=black,inner sep=1pt,minimum size=5pt] at (0,\rt) {};
\node[circle,fill=black,inner sep=1pt,minimum size=5pt] at (2,\rs) {};
\node[circle,fill=black,inner sep=1pt,minimum size=5pt] at (4,\rt) {};
}

\foreach \i in {1,...,3}
{
\pgfmathtruncatemacro{\rt}{2*\i};
\pgfmathtruncatemacro{\rtt}{2*\i-1};
\path[draw=black,thick,->-=0.5] (0,\rt) to (2,\rtt);
\path[draw=black,thick,->-=0.5] (2,\rtt) to (4,\rt);
}

\foreach \i in {4}
{
\pgfmathtruncatemacro{\rt}{2*\i};
\pgfmathtruncatemacro{\rtt}{2*\i-1};
\path[draw=black,thick,loosely dotted] (0,\rt) to (2,\rtt);
\path[draw=black,thick,loosely dotted] (2,\rtt) to (4,\rt);
}

\foreach \i in {0,...,3}
{
\pgfmathtruncatemacro{\rt}{2*\i};

\node[left] at (-0.15,\rt) {\Large {$u_\i$}};

\node[right] at (4.15,\rt) {\Large {$v_\i$}};

\pgfmathsetmacro{\rtpp}{{2*\i+0.57}};

\node[left] at (0.9,\rtpp) {{$a_\i$}};

\node[right] at (3.1,\rtpp) {{$d_\i$}};

\node[below] at (2,\rt) {{$e_\i$}};

}

\foreach \i in {4}
{
\pgfmathtruncatemacro{\rt}{2*\i};

\node[left] at (-0.15,\rt) {\Huge {$\vdots$}};

\node[right] at (4.15,\rt) {\Huge {$\vdots$}};

}

\foreach \i in {1,...,3}
{

\pgfmathsetmacro{\rtp}{{2*\i-0.57}};

\node[left] at (0.9,\rtp) {{$b_\i$}};

\node[right] at (3.1,\rtp) {{$c_\i$}};

}

\end{tikzpicture}

\caption{Digraph representation of the matrix $LU+D$ where 
$L$ is the lower-bidiagonal matrix with 
	$a_0, a_1, \ldots$ on the diagonal and
	$b_1, b_2, \ldots$ on the subdiagonal,
$U$ is the upper-bidiagonal matrix with
        $c_1, c_2, \ldots$ on the superdiagonal and
        $d_0, d_1, \ldots$ on the diagonal,
and $D$ is a diagonal matrix with entries $e_0, e_1,\ldots\;$.
}
\vspace{5mm}
\label{fig.tridiag.LGV}
\end{figure}

2.  In the special case where $A = LU$
with $L$ nonnegative lower-bidiagonal
and $U$ nonnegative upper-bidiagonal,
Proposition~\ref{prop.comparison} can be proven by a
Lindstr\"om--Gessel--Viennot argument
\cite[Chapter~32]{Aigner_18}
using the digraph shown in Figure~\ref{fig.tridiag.LGV}.
\myendremark

\subsection{Production matrices}   \label{subsec.production}

The method of production matrices \cite{Deutsch_05,Deutsch_09}
has become in recent years an important tool in enumerative combinatorics.
In the special case of a tridiagonal production matrix,
this construction goes back to Stieltjes' \cite{Stieltjes_1889,Stieltjes_1894}
work on continued fractions:
the production matrix of a classical S-fraction or J-fraction is tridiagonal.
In~the present paper, by contrast,
we shall need production matrices that are lower-Hessenberg
(i.e.\ vanish above the first superdiagonal)
but are not in general tridiagonal.
We therefore begin by reviewing briefly
the basic theory of production matrices.
The important connection of production matrices with total positivity
will be treated in the next subsection.

Let $P = (p_{ij})_{i,j \ge 0}$ be an infinite matrix
with entries in a commutative ring $R$.
In~order that powers of $P$ be well-defined,
we shall assume that $P$ is either row-finite
(i.e.\ has only finitely many nonzero entries in each row)
or column-finite.

Let us now define an infinite matrix $A = (a_{nk})_{n,k \ge 0}$ by
\be
   a_{nk}  \;=\;  (P^n)_{0k}
 \label{def.iteration}
\ee
(in particular, $a_{0k} = \delta_{0k}$).
Writing out the matrix multiplications explicitly, we have
\be
   a_{nk}
   \;=\;
   \sum_{i_1,\ldots,i_{n-1}}
      p_{0 i_1} \, p_{i_1 i_2} \, p_{i_2 i_3} \,\cdots\,
        p_{i_{n-2} i_{n-1}} \, p_{i_{n-1} k}
   \;,
 \label{def.iteration.walk}
\ee
so that $a_{nk}$ is the total weight for all $n$-step walks in $\N$
from $i_0 = 0$ to $i_n = k$, in~which the weight of a walk is the
product of the weights of its steps, and a step from $i$ to $j$
gets a weight $p_{ij}$.
Yet another equivalent formulation is to define the entries $a_{nk}$
by the recurrence
\be
   a_{nk}  \;=\;  \sum_{i=0}^\infty a_{n-1,i} \, p_{ik}
   \qquad\hbox{for $n \ge 1$}
 \label{def.iteration.bis}
\ee
with the initial condition $a_{0k} = \delta_{0k}$.

We call $P$ the \textbfit{production matrix}
and $A$ the \textbfit{output matrix},
and we write $A = \scro(P)$.
Note that if $P$ is row-finite, then so is $\scro(P)$;
if $P$ is lower-Hessenberg, then $\scro(P)$ is lower-triangular;
if $P$ is lower-Hessenberg with invertible superdiagonal entries,
then $\scro(P)$ is lower-triangular with invertible diagonal entries;
and if $P$ is unit-lower-Hessenberg
(i.e.\ lower-Hessenberg with entries 1 on the superdiagonal),
then $\scro(P)$ is unit-lower-triangular.
In all the applications in this paper, $P$ will be lower-Hessenberg,
most often unit-lower-Hessenberg.

The matrix $P$ can also be interpreted as the adjacency matrix
for a weighted directed graph on the vertex set $\N$
(where the edge $i \to j$ is omitted whenever $p_{ij}  = 0$).
Then $P$ is row-finite (resp.\ column-finite)
if and only if every vertex has finite out-degree (resp.\ finite in-degree);
and $P$ is lower-Hessenberg if and only if
all edges $i \to j$ satisfy $j \le i+1$.

This iteration process can be given a compact matrix formulation.
Let us define the \textbfit{augmented production matrix}
\be
   \widetilde{P}
   \;\eqdef\;
   \left[
    \begin{array}{c}
         1 \;\; 0 \;\; 0 \;\; 0 \;\; \cdots \; \vphantom{\Sigma} \\
         \hline
         P
    \end{array}
    \right]
    \;,
 \label{def.prodmat_augmented}
\ee
or in other words
\be
   (\widetilde{P})_{n,k}
   \;=\;
   \begin{cases}
       \delta_{0k}  & \textrm{if $n=0$} \\
       p_{n-1,k}    & \textrm{if $n \ge 1$}
   \end{cases}
\ee
Then the recurrence \reff{def.iteration.bis}
together with the initial condition $a_{0k} = \delta_{0k}$ can be written as
\be
   A
   \;=\;
   \left[
    \begin{array}{c}
         1 \;\; 0 \;\; 0 \;\; 0 \;\; \cdots \; \vphantom{\Sigma} \\
         \hline
         AP
    \end{array}
    \right]
   \;=\;
   \left[
    \begin{array}{c|c}
         1 & \bzero  \\
         \hline
         \bzero & A
    \end{array}
    \right]
   \left[
    \begin{array}{c}
         1 \;\; 0 \;\; 0 \;\; 0 \;\; \cdots \; \vphantom{\Sigma} \\
         \hline
         P
    \end{array}
    \right]
   \;=\;
   \left[
    \begin{array}{c|c}
         1 & \bzero  \\
         \hline
         \bzero & A
    \end{array}
    \right]
   \widetilde{P}
   \;.
 \label{eq.prodmat.u}
\ee
This identity can be iterated to give the factorization
\be
   A
   \;=\;
   \cdots\,
   \left[
    \begin{array}{c|c}
         I_3 & \bzero  \\
         \hline
                &   \\[-4mm]
         \bzero & \widetilde{P}
    \end{array}
    \right]
   \left[
    \begin{array}{c|c}
         I_2 & \bzero  \\
         \hline
                &   \\[-4mm]
         \bzero & \widetilde{P}
    \end{array}
    \right]
   \left[
    \begin{array}{c|c}
         I_1 & \bzero  \\
         \hline
                &   \\[-4mm]
         \bzero & \widetilde{P}
    \end{array}
    \right]
    \widetilde{P}
 \label{eq.prodmat.u.iterated}
\ee
where $I_k$ is the $k \times k$ identity matrix;
and conversely, \reff{eq.prodmat.u.iterated} implies \reff{eq.prodmat.u}.

Now let $\Delta = (\delta_{i+1,j})_{i,j \ge 0}$
be the matrix with 1 on the superdiagonal and 0 elsewhere.
Then for any matrix $M$ with rows indexed by $\N$,
the product $\Delta M$ is simply $M$ with its zeroth row removed
and all other rows shifted upwards.
(Some authors use the notation $\overline{M} \eqdef \Delta M$.)
The recurrence \reff{def.iteration.bis} can then be written as
\be
   \Delta \, \scro(P)  \;=\;  \scro(P) \, P
   \;.
 \label{def.iteration.bis.matrixform}
\ee
It follows that if $A$ is a row-finite matrix
that has a row-finite inverse $A^{-1}$
and has first row $a_{0k} = \delta_{0k}$,
then $P = A^{-1} \Delta A$ is the unique matrix such that $A = \scro(P)$.
This holds, in particular, if $A$ is lower-triangular with
invertible diagonal entries and $a_{00} = 1$;
then $A^{-1}$ is lower-triangular
and $P = A^{-1} \Delta A$ is lower-Hessenberg.
And if $A$ is unit-lower-triangular,
then $P = A^{-1} \Delta A$ is unit-lower-Hessenberg.

Later we shall need the following easy but fundamental fact,
which shows how the production matrix transforms
when the output matrix $A$ is {\em right}\/-multiplied by another matrix $B$:

\begin{lemma}[Right-multiplication lemma]
   \label{lemma.production.AB}
Let $P = (p_{ij})_{i,j \ge 0}$ be a row-finite matrix
(with entries in a commutative ring $R$),
with output matrix $A = \scro(P)$;
and let $B = (b_{ij})_{i,j \ge 0}$
be a lower-triangular matrix with invertible (in $R$) diagonal entries.
Then
\be
   AB \;=\;  b_{00} \, \scro(B^{-1} P B)
   \;.
 \label{eq.lemma.production.AB}
\ee
That is, up to a factor $b_{00}$,
the matrix $AB$ has production matrix $B^{-1} P B$.
\end{lemma}

\proof
Since $P$ is row-finite, so is $A = \scro(P)$;
then the matrix products $AB$ and $B^{-1} P B$
arising in the lemma are well-defined.  Now
\be
   a_{nk}
   \;=\;
   \sum_{i_1,\ldots,i_{n-1}}
      p_{0 i_1} \, p_{i_1 i_2} \, p_{i_2 i_3} \,\cdots\,
        p_{i_{n-2} i_{n-1}} \, p_{i_{n-1} k}
   \;,
\ee
while
\be
   \scro(B^{-1} P B)_{nk}
   \;=\;
   \sum_{j,i_1,\ldots,i_{n-1},i_n}
      (B^{-1})_{0j} \,
      p_{j i_1} \, p_{i_1 i_2} \, p_{i_2 i_3} \,\cdots\,
        p_{i_{n-2} i_{n-1}} \, p_{i_{n-1} i_n} \, b_{i_n k}
   \;.
\ee
But $B$ is lower-triangular with invertible diagonal entries,
so $B$ is invertible and $B^{-1}$ is lower-triangular,
with $(B^{-1})_{0j} = b_{00}^{-1} \delta_{j0}$.
It follows that $AB = b_{00} \, \scro(B^{-1} P B)$.
\qed

{\bf Remark.}
If $b_{00} = 1$, then $(AB)_{00} = 1$
and hence $AB$ can be the output matrix of some production matrix;
and in this case \reff{eq.lemma.production.AB}
is an immediate consequence of \reff{def.iteration.bis.matrixform}.
Indeed, right-multiplying \reff{def.iteration.bis.matrixform} by $B$ yields
\be
   \Delta \, (\scro(P) \, B)
   \;=\;
   \scro(P) \, PB
   \;=\;
   (\scro(P) \, B) \, B^{-1} P B
   \;.
\ee
\myendremark

We will frequently apply this lemma
with $B$ taken to be the binomial matrix $B_x$.
Unfortunately we know very little about what happens to the
production matrix when the output matrix is {\em left-}\/multiplied
by another matrix $B$.
But we do know one special case: see Lemma~\ref{lemma.aI+bP} below.

\subsection{Production matrices and total positivity}
   \label{subsec.totalpos.prodmat}

Let $P = (p_{ij})_{i,j \ge 0}$ be a matrix with entries in a
partially ordered commutative ring $R$.
We will use $P$ as a production matrix;
let $A = \scro(P)$ be the corresponding output matrix.
As before, we assume that $P$ is either row-finite or column-finite.

When $P$ is totally positive, it turns out \cite{Sokal_totalpos}
that the output matrix $\scro(P)$ has {\em two}\/ total-positivity properties:
firstly, it is totally positive;
and secondly, its zeroth column is Hankel-totally positive.
More generally, the same properties hold whenever $P = P_0 + aI$,
where $P_0$ is totally positive and $a$ is nonnegative.
Since \cite{Sokal_totalpos} is not yet publicly available,
we shall present briefly here (with proof) the main results
that will be needed in the sequel.

The fundamental fact that drives the whole theory is the following:

\begin{proposition}[Minors of the output matrix]
   \label{prop.iteration.homo}
Every $k \times k$ minor of the output matrix $A = \scro(P)$
can be written as a sum of products of minors of size $\le k$
of the production matrix $P$.
\end{proposition}

In this proposition the matrix elements $\bfp = \{p_{ij}\}_{i,j \ge 0}$
should be interpreted in the first instance as indeterminates:
for instance, we can fix a row-finite or column-finite set
$S \subseteq \N \times \N$
and define the matrix $P^S = (p^S_{ij})_{i,j \in \N}$ with entries
\be
   p^S_{ij}
   \;=\;
   \begin{cases}
       p_{ij}  & \textrm{if $(i,j) \in S$} \\[1mm]
       0       & \textrm{if $(i,j) \notin S$}
   \end{cases}
\ee
Then the entries (and hence also the minors) of both $P$ and $A$
belong to the polynomial ring $\Z[\bfp]$,
and the assertion of Proposition~\ref{prop.iteration.homo} makes sense.
Of course, we can subsequently specialize the indeterminates $\bfp$
to values in any commutative ring $R$.

%
\proofof{Proposition~\ref{prop.iteration.homo}}
%
For any infinite matrix $X = (x_{ij})_{i,j \ge 0}$,
let us write $X_N = (x_{ij})_{0 \le i \le N-1 ,\, j \ge 0}$
for the submatrix consisting of the first $N$ rows
(and {\em all}\/ the columns) of $X$.
Every $k \times k$ minor of $A$ is of course
a $k \times k$ minor of $A_N$ for some $N$,
so it suffices to prove that the claim about minors holds for all the $A_N$.
But this is easy: the fundamental identity \reff{eq.prodmat.u} implies
\be
   A_N
   \;=\;
   \left[
    \begin{array}{c|c}
         1 & \bzero  \\
         \hline
         \bzero & A_{N-1}
    \end{array}
    \right]
   \,
   \left[
    \begin{array}{c}
         1 \;\; 0 \;\; 0 \;\; 0 \;\; \cdots \; \vphantom{\Sigma} \\
         \hline
         P
    \end{array}
    \right]
   \;.
 \label{eq.proof.prop.iteration.homo}
\ee
So the result follows by induction on $N$, using the Cauchy--Binet formula.
\qed

If we now specialize the indeterminates $\bfp$
to values in some partially ordered commutative ring $R$,
we can immediately conclude:

\begin{theorem}[Total positivity of the output matrix]
   \label{thm.iteration.homo}
Let $P$ be an infinite matrix that is either row-finite or column-finite,
with entries in a partially ordered commutative ring $R$.
If $P$ is totally positive of order~$r$, then so is $A = \scro(P)$.
\end{theorem}

\medskip

{\bf Remarks.}
1.  In the case $R = \R$, Theorem~\ref{thm.iteration.homo}
is due to Karlin \cite[pp.~132--134]{Karlin_68};
see also \cite[Theorem~1.11]{Pinkus_10}.
Karlin's proof is different from ours.

2.  Our quick inductive proof of Proposition~\ref{prop.iteration.homo}
follows an idea of Zhu \cite[proof of Theorem~2.1]{Zhu_13},
which was in turn inspired in part by Aigner \cite[pp.~45--46]{Aigner_99}.
The same idea recurs in recent work of several authors
\cite[Theorem~2.1]{Zhu_14}
\cite[Theorem~2.1(i)]{Chen_15a}
\cite[Theorem~2.3(i)]{Chen_15b}
\cite[Theorem~2.1]{Liang_16}
\cite[Theorems~2.1 and 2.3]{Chen_19}
\cite{Gao_non-triangular_transforms}.
However, all of these results concerned only special cases:
\cite{Aigner_99,Zhu_13,Chen_15b,Liang_16}
treated the case in which the production matrix $P$ is tridiagonal;
\cite{Zhu_14} treated a (special) case in which $P$ is upper-bidiagonal;
\cite{Chen_15a} treated the case in which
$P$ is the production matrix of a Riordan array;
\cite{Chen_19,Gao_non-triangular_transforms}
treated (implicitly) the case in which $P$ is upper-triangular and Toeplitz.
But the argument is in fact completely general, as we have just seen;
there is no need to assume any special form for the matrix $P$.

3. A slightly different version of this proof
was presented in \cite{latpath_SRTR,latpath_lah}.
The simplified reformulation
given here,
using the augmented production matrix,
is due to Mu and Wang \cite{Mu_20}.
\myendremark

\medskip

\begin{example}[Binomial matrix]
   \label{exam.binomial.matrix.TP}
\rm
Let $P$ be the upper-bidiagonal Toeplitz matrix
$xI + y\Delta$, where $x$ and $y$ are indeterminates.
By Lemma~\ref{lemma.bidiagonal}, $P$ is TP
in the ring $\Z[x,y]$ equipped with the coefficientwise order.
An easy computation shows that $\scro(xI + y\Delta) = B_{x,y}$,
the weighted binomial matrix
with entries $(B_{x,y})_{nk} = x^{n-k} y^k \binom{n}{k}$.
So Theorem~\ref{thm.iteration.homo} implies that $B_{x,y}$ is TP
in the ring $\Z[x,y]$ equipped with the coefficientwise order.
This gives an {\em ab initio}\/ proof of Lemma~\ref{lemma.binomialmatrix.TP}.
\myendremark
\end{example}

More generally, we have:

\begin{lemma}
   \label{lemma.aI+bP}
\quad
$\scro(aI+bP) \:=\: B_{a,b} \, \scro(P)$.
\end{lemma}

\noindent
Note that Example~\ref{exam.binomial.matrix.TP}
is the special case $P = \Delta$.

\proofof{Lemma~\ref{lemma.aI+bP}}
Since
\be
   (aI+bP)^n  \;=\;  \sum_{j=0}^n a^{n-j} \, b^j \, \binom{n}{j} \, P^j
   \;,
\ee
we have
\begin{subeqnarray}
   \scro(aI+bP)_{nk}
   \;\eqdef\;
   [(aI+bP)^n]_{0k}
   & = &
   \sum_{j=0}^n a^{n-j} \, b^j \, \binom{n}{j} \, (P^j)_{0k}
       \qquad \\[2mm]
   & = &
   \sum_{j=0}^n (B_{a,b})_{nj} \, \scro(P)_{jk}
      \slabel{eq.proof.lemma.aI+bP.b}  \\[2mm]
   & = &
   [B_{a,b} \, \scro(P)]_{nk}  \;.
       \\[-9mm] \nonumber
\end{subeqnarray}
\qed


Since the binomial matrix $B_{a,b}$ is totally positive
by Lemma~\ref{lemma.binomialmatrix.TP}, we conclude:

\begin{corollary}[Total positivity of output matrix, improved]
   \label{cor.iteration.homo.aI+bP}
Let $P$ be an infinite matrix that is either row-finite or column-finite,
with entries in a partially ordered commutative ring $R$.
Then $P$ is TP${}_r$ $\implies$ $\scro(P)$ is TP${}_r$
$\iff$ $\scro(aI+bP)$ is TP${}_r$ for all $a,b \ge 0$.
%
\end{corollary}

\noindent
Here $a$ and $b$ can be nonnegative elements of $R$;
or $a$ and $b$ can be indeterminates
and we work in the ring $R[a,b]$ equipped with the coefficientwise order.
Corollary~\ref{cor.iteration.homo.aI+bP} allows us to prove the
total positivity of the output matrix ${\scro(aI+bP)}$ whenever $P$ is TP,
{\em even if $aI+bP$ is not TP}\/.
This will play a crucial role in our proof of Theorem~\ref{thm.multivariate},
by allowing us to take $\lambda\yfp \ge \lambda\yp$
rather than just $\lambda\yfp = \lambda\yp$.

\bigskip

{\bf Remark.}
Lemma~\ref{lemma.aI+bP} and Corollary~\ref{cor.iteration.homo.aI+bP}
are a special case of an idea of Zhu \cite{Zhu_18};
see \cite{Sokal_totalpos} for a general version.
\myendremark

\medskip

\bigskip

Now define 
$\scroo_0(P)$ to be the zeroth-column sequence of $\scro(P)$, i.e.
\be
   \scroo_0(P)_n  \;\eqdef\;  \scro(P)_{n0}  \;\eqdef\;  (P^n)_{00}
   \;.
 \label{def.scroo0}
\ee
Then the Hankel matrix of $\scroo_0(P)$ has matrix elements
\begin{eqnarray}
   & &
   \!\!\!\!\!\!\!
   H_\infty(\scroo_0(P))_{nn'}
   \;=\;
   \scroo_0(P)_{n+n'}
   \;=\;
   (P^{n+n'})_{00}
   \;=\;
   \sum_{k=0}^\infty (P^n)_{0k} \, (P^{n'})_{k0}
   \;=\;
          \nonumber \\
   & &
   \sum_{k=0}^\infty (P^n)_{0k} \, ((P^{\rm T})^{n'})_{0k}
   \;=\;
   \sum_{k=0}^\infty \scro(P)_{nk} \, \scro(P^{\rm T})_{n'k}
   \;=\;
   \big[ \scro(P) \, {\scro(P^{\rm T})}^{\rm T} \big]_{nn'}
   \;.
   \qquad
\end{eqnarray}
(Note that the sum over $k$ has only finitely many nonzero terms:
 if $P$ is row-finite, then there are finitely many nonzero $(P^n)_{0k}$,
 while if $P$ is column-finite,
 there are finitely many nonzero $(P^{n'})_{k0}$.\footnote{
    Or to put it another way:
    If $P$ is row-finite, then $\scro(P)$ is row-finite;
    $\scro(P^{\rm T})$ need not be row- or column-finite,
    but the product $\scro(P) \, {\scro(P^{\rm T})}^{\rm T}$
    is anyway well-defined.
    Likewise, if $P$ is column-finite,
    then ${\scro(P^{\rm T})}^{\rm T}$ is column-finite;
    $\scro(P)$ need not be row- or column-finite,
    but the product $\scro(P) \, {\scro(P^{\rm T})}^{\rm T}$
    is anyway well-defined.
})
We have therefore proven:

\begin{lemma}[Identity for Hankel matrix of the zeroth column]
   \label{lemma.hankel.karlin}
Let $P$ be a row-finite or column-finite matrix
with entries in a commutative ring $R$.
Then
\be
   H_\infty(\scroo_0(P))
   \;=\;
   \scro(P) \, {\scro(P^{\rm T})}^{\rm T}
   \;.
\ee
\end{lemma}

Combining Proposition~\ref{prop.iteration.homo}
with Lemma~\ref{lemma.hankel.karlin} and the Cauchy--Binet formula,
we obtain:

\begin{corollary}[Hankel minors of the zeroth column]
   \label{cor.iteration2}
Every $k \times k$ minor of the infinite Hankel matrix
$H_\infty(\scroo_0(P)) = ((P^{n+n'})_{00})_{n,n' \ge 0}$
can be written as a sum of products
of the minors of size $\le k$ of the production matrix $P$.
\end{corollary}

And specializing the indeterminates $\bfp$
to nonnegative elements in a partially ordered commutative ring,
in such a way that $P$ is row-finite or column-finite,
we deduce:

\begin{theorem}[Hankel-total positivity of the zeroth column]
   \label{thm.iteration2bis}
Let $P = (p_{ij})_{i,j \ge 0}$ be an infinite row-finite or column-finite
matrix with entries in a partially ordered commutative ring $R$,
and define the infinite Hankel matrix
$H_\infty(\scroo_0(P)) = ((P^{n+n'})_{00})_{n,n' \ge 0}$.
If $P$ is totally positive of order~$r$, then so is $H_\infty(\scroo_0(P))$.
\end{theorem}

Once again we can we can improve these results
to replace $P$ by $aI+bP$.
It suffices to note that
\begin{subeqnarray}
   H_\infty(\scroo_0(aI+bP))
   & = &
   \scro(aI+bP) \, {\scro(aI + bP^{\rm T})}^{\rm T}
       \\[2mm]
   & = &
   B_{a,b} \, \scro(P) \, {\scro(P^{\rm T})}^{\rm T} \, (B_{a,b})^{\rm T}
       \\[2mm]
   & = &
   B_{a,b} \, H_\infty(\scroo_0(P)) \, (B_{a,b})^{\rm T}
\end{subeqnarray}
by Lemmas~\ref{lemma.aI+bP} and \ref{lemma.hankel.karlin}.
Since the binomial matrix $B_{a,b}$ is totally positive
by Lemma~\ref{lemma.binomialmatrix.TP}, we conclude:

\begin{corollary}[Hankel-total positivity of zeroth column, improved]
   \label{thm.iteration2bis.aI+bP}
Let $P$ be an infinite matrix that is either row-finite or column-finite,
with entries in a partially ordered commutative ring $R$.
Then $P$ is TP${}_r$ $\implies$ $\scro(P)$ and $\scro(P^{\rm T})$ are TP${}_r$
$\implies$ $\scroo_0(P)$ is Hankel-TP${}_r$
$\iff$ $\scroo_0(aI+bP)$ is Hankel-TP${}_r$ for all $a,b \ge 0$.
\end{corollary}

\noindent
Once again, $a$ and $b$ can be nonnegative elements of $R$,
or $a$ and $b$ can be indeterminates
and we work in the ring $R[a,b]$ equipped with the coefficientwise order.
Corollary~\ref{thm.iteration2bis.aI+bP} allows us to prove the
Hankel-total positivity of the output sequence ${\scroo_0(aI+bP)}$
whenever $P$ is TP, even if $aI+bP$ is not TP.
It will play an important role in our proof of Theorem~\ref{thm.multivariate}.

\medskip

{\bf Remark.}
We see from \reff{eq.proof.lemma.aI+bP.b} specialized to $k=0$
that $\scroo_0(aI+bP)$ is a binomial transform of $\scroo_0(P)$:
\be
   \scroo_0(aI+bP)_n
   \;=\;
   \sum_{j=0}^n a^{n-j} \, b^j \, \binom{n}{j} \, \scroo_0(P)_j
   \;.
\ee
Now, it is known \cite{Zhu_19,Sokal_totalpos}
that the binomial transform preserves Hankel-TP${}_r$ of arbitrary sequences
in a partially ordered commutative ring.
So Corollary~\ref{thm.iteration2bis.aI+bP} is just a special case
of this general result combined with Theorem~\ref{thm.iteration2bis}.
\myendremark

%

\subsection{Binomial row-generating matrices}  \label{subsec.rowgen}

Let $A = (a_{nk})_{n,k \ge 0}$ be a row-finite matrix
with entries in a commutative ring $R$.
(In most applications, including all those in the present paper,
 the matrix $A$ will be lower-triangular.)
We define its \textbfit{row-generating polynomials} in the usual way:
\be
   A_n(x)  \;\eqdef\;  \sum_{k=0}^\infty a_{nk} \, x^k
   \;,
 \label{def.An.0}
\ee
where the sum is actually finite because $A$ is row-finite.
More generally, let us define its
\textbfit{binomial partial row-generating polynomials}
\begin{subeqnarray}
   A_{n,k}(x)
   & \eqdef &
   \sum_{\ell=k}^\infty a_{n\ell} \, \binom{\ell}{k} \, x^{\ell-k}
         \\[2mm]
   & = &
   {1 \over k!} \, {d^k \over dx^k} \, A_n(x)
   \;.
 \label{def.Ank}
\end{subeqnarray}
(Note that the operator $(1/k!) \, d^k\!/\!dx^k$ has a well-defined action
 on the polynomial ring $R[x]$ even if $R$ does not contain the rationals,
 since $(1/k!) (d^k\!/\!dx^k) x^n = \binom{n}{k} x^{n-k}$.)
The polynomials $A_{n,k}(x)$ are the matrix elements of the
\textbfit{binomial row-generating matrix} $A B_x$:
\be
   (A B_x)_{nk}  \;=\;  A_{n,k}(x)  \;,
\ee
where $B_x = B_{x,1}$ is the weighted binomial matrix defined in \reff{def.Bx}.
The zeroth column of the matrix $A B_x$ consists
of the row-generating polynomials $A_n(x) = A_{n,0}(x)$.

In this paper the matrix $A$
will be either the Laguerre coefficient matrix $\sfL^{(\alpha)}$
or one of its multivariate generalizations.

We can now explain the method that we will use to prove
Theorems~\ref{thm.univariate} and \ref{thm.multivariate}:

\begin{proposition}
   \label{prop.method}
Let $P$ be a row-finite matrix
with entries in a partially ordered commutative ring $R$,
and let $A = \scro(P)$.
\begin{itemize}
   \item[(a)]  If $P$ is totally positive of order~$r$, then so is $A$.
   \item[(b)]  If the matrix $B_x^{-1} P B_x$ is totally positive of order~$r$
in the ring $R[x]$ equipped with the coefficientwise order,
then the sequence $(A_n(x))_{n \ge 0}$ of row-generating polynomials
is Hankel-totally positive of order~$r$
in the ring $R[x]$ equipped with the coefficientwise order.
\end{itemize}
\end{proposition}

\noindent
Indeed, (a) is just a restatement of Theorem~\ref{thm.iteration.homo};
and (b) is an immediate consequence of Lemma~\ref{lemma.production.AB}
and Theorem~\ref{thm.iteration2bis}
together with the fact that the zeroth column of the matrix $A B_x$ consists
of the row-generating polynomials $A_n(x)$.

\bigskip

{\bf Remark.}
The binomial row-generating matrix $A B_x$
can also be considered as a modified Wronskian matrix.
To see this, let $N \in \N \cup \{\infty\}$,
and fix formal power series $f_1(x), \ldots, f_N(x) \in R[[x]]$.
We define the $N \times \infty$ \textbfit{Wronskian matrix}
\be
   W(f_1,\ldots,f_N)(x)
   \;=\;
   \big( f_n^{(k)}(x) \bigr)_{1 \le n \le N,\, k \ge 0}
\ee
(where ${}^{(k)}$ denotes the $k$th derivative)
and the $N \times \infty$ \textbfit{modified Wronskian matrix}
\be
   \widetilde{W}(f_1,\ldots,f_N)(x)
   \;=\;
   \big( f_n^{(k)}(x)/k! \bigr)_{1 \le n \le N,\, k \ge 0}
   \;;
\ee
these are matrices with entries in $R[[x]]$.
Of course, we have
\be
   W(f_1,\ldots,f_N)(x)  \;=\;  \widetilde{W}(f_1,\ldots,f_N)(x) \: D
\ee
where $D = \diag\big( (k!)_{k \ge 0} \big)$.
The key fact is that
\be
   \widetilde{W}(f_1,\ldots,f_N)(x)
   \;=\;
   \widetilde{W}(f_1,\ldots,f_N)(0) \: B_x
   \;.
\ee
In particular, if the ring $R$ carries a partial order,
the total positivity of order~$r$
of the matrix $\widetilde{W}(f_1,\ldots,f_N)(0)$
in the ring $R$
implies the total positivity of order~$r$
of the matrix $\widetilde{W}(f_1,\ldots,f_N)(x)$
in the ring $R[[x]]$ equipped with the coefficientwise order.
\myendremark

\subsection{Exponential Riordan arrays}
   \label{sec.exp_riordan}

Let $R$ be a commutative ring containing the rationals,
and let $F(t) = \sum_{n=0}^\infty f_n t^n/n!$
and $G(t) = \sum_{n=1}^\infty g_n t^n/n!$ be formal power series
with coefficients in $R$; we set $g_0 = 0$.
Then the \textbfit{exponential Riordan array}
\cite{Deutsch_04,Deutsch_09,Barry_16}
associated to the pair $(F,G)$
--- or equivalently to the pair of sequences
$\bff = (f_n)_{n \ge 0}$ and $\bg = (g_n)_{n \ge 1}$ ---
is the infinite lower-triangular matrix
$\scrr[F,G] = (\scrr[F,G]_{nk})_{n,k \ge 0}$ defined by
\be
   \scrr[F,G]_{nk}
   \;=\;
   {n! \over k!} \:
   [t^n] \, F(t) G(t)^k
   \;.
 \label{def.exp_riordan}
\ee
That is, the $k$th column of $\scrr[F,G]$
has exponential generating function $F(t) G(t)^k/k!$.
It follows that the bivariate exponential generating function
of $\scrr[F,G]$ is
\be
   \sum_{n,k=0}^\infty \scrr[F,G]_{nk} \: {t^n \over n!} \, u^k
   \;=\;
   F(t) \, e^{u G(t)}
   \;.
 \label{eq.exp_riordan.bivariate_egf}
\ee
Please note that the diagonal elements of $\scrr[F,G]$
are $\scrr[F,G]_{nn} = f_0 g_1^n$,
so the matrix $\scrr[F,G]$ is invertible
in the ring $R^{\N \times \N}_{\rm lt}$ of lower-triangular matrices
if and only if $f_0$ and $g_1$ are invertible in $R$.

We shall use an easy but important result that is sometimes called
the \emph{fundamental theorem of exponential Riordan arrays} (FTERA):

\begin{lemma}[Fundamental theorem of exponential Riordan arrays]
   \label{lemma.FETRA}
Let $\bb = (b_n)_{n \ge 0}$ be a sequence with
exponential generating function $B(t) = \sum_{n=0}^\infty b_n t^n/n!$.
Considering $\bb$ as a column vector and letting $\scrr[F,G]$
act on it by matrix multiplication, we obtain a sequence $\scrr[F,G] \bb$
whose exponential generating function is $F(t) \, B(G(t))$.
\end{lemma}

\proof
We compute
\begin{subeqnarray}
   \sum_{k=0}^n \scrr[F,G]_{nk} \, b_k
   & = &
   \sum_{k=0}^\infty {n! \over k!} \, [t^n] \, F(t) G(t)^k \, b_k
            \\[2mm]
   & = &
   n! \: [t^n] \: F(t) \sum_{k=0}^\infty b_k \, {G(t)^k \over k!}
            \\[2mm]
   & = &
   n! \: [t^n] \: F(t) \, B(G(t))
   \;.
\end{subeqnarray}
\qed

Let us now consider the product of two exponential Riordan arrays
$\scrr[F_1,G_1]$ and $\scrr[F_2,G_2]$.
Applying the FTERA to the $k$th column of $\scrr[F_2,G_2]$,
whose exponential generating function is $F_2(t) G_2(t)^k/k!$,
we readily obtain:

\begin{lemma}[Product of two exponential Riordan arrays]
   \label{lemma.exp_riordan.product}
We have
\be
   \scrr[F_1,G_1] \, \scrr[F_2,G_2]
   \;=\;
   \scrr[ (F_2 \circ G_1) F_1 ,\, G_2 \circ G_1]
 \label{eq.exp_riordan.composition}
\ee
where $\circ$ denotes composition of formal power series.
\end{lemma}

In particular, if we let $\scrr[F_2,G_2]$
be the weighted binomial matrix $B_\xi = \scrr[e^{\xi t}, t]$,
we obtain:

\begin{corollary}[Binomial row-generating matrix of an exponential Riordan array]
   \label{cor.exp_riordan.Bx}
We have
\be
   \scrr[F,G] \, B_\xi
   \;=\;
   \scrr[\e^{\xi G} F, G]
   \;.
 \label{eq.cor.exp_riordan.Bx}
\ee
\end{corollary}

We can now determine the production matrix of an exponential
Riordan array $\scrr[F,G]$.
Let $\ba = (a_n)_{n \ge 0}$ and $\bz = (z_n)_{n \ge 0}$
be sequences in a commutative ring $R$,
with ordinary generating functions
$A(s) = \sum_{n=0}^\infty a_n s^n$
and $Z(s) = \sum_{n=0}^\infty z_n s^n$.
We then define the \textbfit{exponential AZ matrix}
associated to the sequences $\ba$ and $\bz$
to be the lower-Hessenberg matrix with entries
\be
   \EAZ(\ba,\bz)_{nk}
   \;=\;
   {n! \over k!} \: (z_{n-k} \,+\, k \, a_{n-k+1})
 \label{def.EAZ.1}
\ee
(where $z_{-1} \eqdef 0$),
or equivalently (if $R$ contains the rationals)
\be
   \EAZ(\ba,\bz)
   \;=\;
   D \, T_\infty(\bz) \, D^{-1} \:+\: D \, T_\infty(\ba) \, D^{-1} \, \Delta
 \label{def.EAZ.2}
\ee
where $D = \diag\bigl( (n!)_{n \ge 0} \bigr)$.
We also write $\EAZ(A,Z)$ as a synonym for $\EAZ(\ba,\bz)$.

\smallskip

\begin{theorem}[Production matrices of exponential Riordan arrays]
   \label{thm.riordan.exponential.production}
Let $L$ be a lower-triangular matrix
(with entries in a commutative ring $R$ containing the rationals)
with invertible diagonal entries and $L_{00} = 1$,
and let $P = L^{-1} \Delta L$ be its production matrix.
Then $L$ is an exponential Riordan array
if and only~if $P$ is an exponential AZ matrix.

More precisely, $L = \scrr[F,G]$ if and only~if $P = \EAZ(A,Z)$,
where the generating functions $\big( F(t), G(t) \big)$
and $\big( A(s), Z(s) \big)$ are connected by
\be
   G'(t) \;=\; A(G(t))  \;,\qquad
   {F'(t) \over F(t)} \;=\; Z(G(t))
 \label{eq.prop.riordan.exponential.production.1}
\ee
or equivalently
\be
   A(s)  \;=\;  G'(\bar{G}(s))  \;,\qquad
   Z(s)  \;=\;  {F'(\bar{G}(s)) \over F(\bar{G}(s))}
 \label{eq.prop.riordan.exponential.production.2}
\ee
where $\bar{G}(s)$ is the compositional inverse of $G(t)$.
\end{theorem}

\par\bigskip\noindent{\sc Proof}
(mostly contained in \cite[pp.~217--218]{Barry_16}).
Suppose that $L = \scrr[F,G]$.
The hypotheses on $L$ imply that $f_0 = 1$
and that $g_1$ is invertible in $R$;
so $G(t)$ has a compositional inverse.
Now let $P = (p_{nk})_{n,k \ge 0}$ be a matrix;
its column exponential generating functions are, by definition,
$P_k(t) = \sum_{n=0}^\infty p_{nk} \, t^n/n!$.
Applying the FTERA to each column of $P$,
we see that $\scrr[F,G] P$ is a matrix
whose column exponential generating functions
are $\big( F(t) \, P_k(G(t)) \big)_{k \ge 0}$.
On~the other hand, $\Delta \, \scrr[F,G]$
is the matrix $\scrr[F,G]$ with its zeroth row removed
and all other rows shifted upwards,
so it has column exponential generating functions
\be
   {d \over dt} \, \big( F(t) \, G(t)^k/k! \big)
   \;=\;
   {1 \over k!} \: \Big[ F'(t) \, G(t)^k
                         \:+\: k \, F(t) \, G(t)^{k-1} \, G'(t) \Big]
   \;.
\ee
Comparing these two results, we see that
$\Delta \, \scrr[F,G] = \scrr[F,G] \, P$
if and only~if
\be
   P_k(G(t))
   \;=\;
   {1 \over k!} \:
   {F'(t) \, G(t)^k \:+\: k \, F(t) \, G(t)^{k-1} \, G'(t)
    \over
    F(t)}
   \;,
\ee
or in other words
\be
   P_k(t)
   \;=\;
   {1 \over k!}  \:
      \biggl[ {F'(\bar{G}(t)) \over F(\bar{G}(t))} \, t^k
              \:+\: k \, t^{k-1} \, G'(\bar{G}(t))
      \biggr]
   \;.
\ee
Therefore
\begin{subeqnarray}
   p_{nk}
   & = &
   {n! \over k!} \: [t^n] \,
      \biggl[ {F'(\bar{G}(t)) \over F(\bar{G}(t))} \, t^k
              \:+\: k \, t^{k-1} \, G'(\bar{G}(t))
      \biggr]
    \\[2mm]
   & = &
   {n! \over k!} \:
      \biggl[ [t^{n-k}] \: {F'(\bar{G}(t)) \over F(\bar{G}(t))}
              \:+\: k \, [t^{n-k+1}] \: G'(\bar{G}(t))
      \biggr]
    \\[2mm]
   & = &
   {n! \over k!} \: (z_{n-k} \,+\, k \, a_{n-k+1})
\end{subeqnarray}
where $\ba = (a_n)_{n \ge 0}$ and $\bz = (z_n)_{n \ge 0}$
are given by \reff{eq.prop.riordan.exponential.production.2}.

Conversely, suppose that $P = \EAZ(A,Z)$.
Define $F(t)$ and $G(t)$
as the unique solutions (in the formal-power-series ring $R[[t]]$)
of the differential equations \reff{eq.prop.riordan.exponential.production.1}
with initial conditions $F(0) = 1$ and $G(0) = 0$.
Then running the foregoing computation backwards
shows that $\Delta \, \scrr[F,G] = \scrr[F,G] \, P$.
\qed

\medskip

A central role will be played later in this paper
by a simple but remarkable identity for $B_x^{-1} \EAZ(\ba,\bz) \, B_x$,
where $B_x$ is the $x$-binomial matrix defined in \reff{def.Bx}:

\begin{lemma}[Identity for $B_x^{-1} \EAZ(\ba,\bz) \, B_x$]
   \label{lemma.BxinvEAZBx}
Let $\ba = (a_n)_{n \ge 0}$, $\bz = (z_n)_{n \ge 0}$ and $x$
be indeterminates.
Then
\be
   B_x^{-1} \EAZ(\ba,\bz) \, B_x
   \;=\;
   \EAZ(\ba,\bz+x\ba)
 \label{eq.lemma.BxinvEAZBx}
\ee
as an identity in $\Z[\ba,\bz]$.
\end{lemma}

The special case $\bz = \bzero$ of this lemma was proven in
\cite[Lemma~3.6]{latpath_lah};
a simpler proof was given in \cite[Lemma~2.16]{forests_totalpos}.
Here we give the easy generalization to include~$\bz$,
taken from \cite{Chen-Sokal_trees_totalpos}:

\proof
We work temporarily in the ring $\Q[\ba,\bz]$
and use the matrix definition \reff{def.EAZ.2}:
\be
   \EAZ(\ba,\bz)
   \;=\;
   D \, T_\infty(\bz) \, D^{-1} \:+\: D \, T_\infty(\ba) \, D^{-1} \, \Delta
\ee
where $D = \diag\bigl( (n!)_{n \ge 0} \bigr)$.
Since $\EAZ(\ba,\bz) = \EAZ(\ba,\bzero) + \EAZ(\bzero,\bz)$,
it suffices to consider separately the two contributions.

The key observation is that
$B_x = D \, T_\infty\big( (x^n/n!)_{n \ge 0} \big) \, D^{-1}$.
Now two Toeplitz matrices always commute:
$T_\infty(\ba) \, T_\infty(\bb) =
 T_\infty(\ba \star \bb) =
 T_\infty(\bb) \, T_\infty(\ba)$.
It follows that $D T_\infty(\bz) D^{-1}$ and $D T_\infty(\ba) D^{-1}$
commute with $B_x$.
Therefore
\be
   B_x^{-1} \EAZ(\bzero,\bz) \, B_x
   \;=\;
   \EAZ(\bzero,\bz)
   \;.
 \label{eq.EAZ.bz}
\ee
On the other hand, the classic recurrence for binomial coefficients implies
\be
   \Delta B_x  \;=\; B_x \, (xI + \Delta)
\ee
(cf.\ Example~\ref{exam.binomial.matrix.TP}).  Therefore
\begin{subeqnarray}
   B_x^{-1} \EAZ(\ba,\bzero) B_x
   & = &
   B_x^{-1} \, D T_\infty(\ba) D^{-1} \, \Delta \, B_x
      \\[2mm]
   & = &
   B_x^{-1} \, D T_\infty(\ba) D^{-1} \, B_x \, (xI + \Delta)
      \\[2mm]
   & = &
   D T_\infty(\ba) D^{-1} \, (xI + \Delta)
      \\[2mm]
   & = &
   \EAZ(\ba,x\ba)
   \;.
 \label{eq.EAZ.ba}
\end{subeqnarray}
Adding \reff{eq.EAZ.bz} and \reff{eq.EAZ.ba} yields \reff{eq.lemma.BxinvEAZBx}.
\qed

The identity \reff{eq.lemma.BxinvEAZBx} can alternatively be proven
by combining Lemma~\ref{lemma.production.AB}
and Corollary~\ref{cor.exp_riordan.Bx}
with Theorem~\ref{thm.riordan.exponential.production}:
see \cite{Chen-Sokal_trees_totalpos} for details.

\section{Production matrices: Univariate case} \label{sec.prodmat.1var}

In this section we will prove Proposition~\ref{prop.prodmat.univariate.NEW},
which gives the production matrices
for the Laguerre coefficient matrix $\sfL^{(\alpha)}$
[defined in \reff{def.coeffmat}]
and for the binomial row-generating matrix $\sfL^{(\alpha)} B_x$.
The proofs are in fact quite easy.
We also give a simple direct proof of the coefficientwise total positivity
of the Laguerre coefficient matrix $\sfL^{(-1+\lambda)}$,

\subsection[Coefficient matrix $\sfL^{(\alpha)}$: Proof of Proposition~\ref{prop.prodmat.univariate.NEW}(a)]{Coefficient matrix $\sfL^{\bm{(\alpha)}}$: Proof of Proposition~\ref{prop.prodmat.univariate.NEW}(a)}
   \label{subsec.prodmat.1var.L}

We begin by proving Proposition~\ref{prop.prodmat.univariate.NEW}(a),
which asserts that the production matrix
for the Laguerre coefficient matrix $\sfL^{(\alpha)}$
is the tridiagonal unit-lower-Hessenberg matrix
$P^\circ = (p^\circ_{ij})_{i,j \ge 0}$ defined by
\begin{subeqnarray}
   p^\circ_{n,n+1}  & = &   1   \\
   p^\circ_{n,n}    & = &   2n+1+\alpha   \\
   p^\circ_{n,n-1}  & = &   n(n+\alpha)   \\
   p^\circ_{n,k}    & = &   0 \qquad\textrm{if $k < n-1$ or $k > n+1$}
 \label{eq.prop.prodmat.univariate.NEW.1.BIS}
\end{subeqnarray}
We give two proofs:
one by direct computation,
and one using the theory of exponential Riordan arrays.

\firstproof
Recall that the Laguerre coefficient matrix $\sfL^{(\alpha)}$
has entries
\be
   \ell^{(\alpha)}_{n,k}
   \;\eqdef\;
   \binom{n}{k} \, (\alpha+1+k)^{\overline{n-k}}
   \;.
\ee
We need to verify that
$\sum_j \ell^{(\alpha)}_{n,j} \, p^\circ_{j,k} = \ell^{(\alpha)}_{n+1,k}$,
or in other words that
\be
   \ell^{(\alpha)}_{n,k-1}  \:+\:  (2k+1+\alpha) \ell^{(\alpha)}_{n,k}
       \:+\: (k+1)(k+1+\alpha) \ell^{(\alpha)}_{n,k+1}
   \;=\;
   \ell^{(\alpha)}_{n+1,k}
   \;.
\ee
This is a straightforward computation.
\qed

\secondproof
It follows from \reff{eq.egf.scrlna} and \reff{eq.exp_riordan.bivariate_egf}
that the matrix $\sfL^{(\alpha)}$ is an exponential Riordan array $\scrr[F,G]$
with $F(t) = (1-t)^{-(1+\alpha)}$ and $G(t) = t/(1-t)$.
Since
\be
   G'(t)  \;=\;  [1 + G(t)]^2  \;,
\ee
from \reff{eq.prop.riordan.exponential.production.1} we have
\be
   A(s)  \;=\;  1 + 2s + s^2  \;.
 \label{eq.prodmat1.A}
\ee
Then
\be
   {F'(t) \over F(t)}  \;=\;  {1 + \alpha \over 1-t}
                       \;=\;  (1+\alpha) \, [1 + G(t)]
   \;,
\ee
so from \reff{eq.prop.riordan.exponential.production.1} we have
\be
   Z(s)  \;=\;  (1+\alpha) \, (1+s)  \;.
 \label{eq.prodmat1.Z}
\ee
By Theorem~\ref{thm.riordan.exponential.production},
the production matrix for $\sfL^{(\alpha)}$
is the exponential AZ matrix $\EAZ(A,Z)$.
Inserting \reff{eq.prodmat1.A}/\reff{eq.prodmat1.Z} into \reff{def.EAZ.1}
yields \reff{eq.prop.prodmat.univariate.NEW.1.BIS}.
\qed

We now wish to make an observation about the connection
of the Laguerre coefficient matrix $\sfL^{(\alpha)}$ with continued fractions.
The zeroth column of $\sfL^{(\alpha)}$ is the sequence
of rising powers (= Stirling cycle polynomials)
$(\lambda^{\overline{n}})_{n \ge 0}$, where $\lambda \eqdef \alpha + 1$.
And this sequence
has a well-known classical S-fraction for its ordinary generating function,
which was found more than two-and-a-half centuries ago by Euler
\cite[section~26]{Euler_1760} \cite{Euler_1788}\footnote{
   The paper \cite{Euler_1760},
   which is E247 in Enestr\"om's \cite{Enestrom_13} catalogue,
   was probably written circa 1746;
   it~was presented to the St.~Petersburg Academy in 1753,
   and published in 1760.
   The paper \cite{Euler_1788},
   which is E616 in Enestr\"om's \cite{Enestrom_13} catalogue,
   was apparently presented to the St.~Petersburg Academy in 1776,
   and published posthumously in 1788.
}:
\be
   \sum_{n=0}^\infty \lambda^{\overline{n}} \: t^n
   \;=\;
   \cfrac{1}{1 - \cfrac{\lambda t}{1 - \cfrac{t}{1 - \cfrac{(\lambda+1)t}{1- \cfrac{2t}{1- \cdots}}}}}
  \label{eq.stirlingcycle.Sfraction.NEW}
\ee
with coefficients $\alpha_{2k-1} = \lambda+k-1$ and $\alpha_{2k} = k$.
Now, any S-fraction has an associated production matrix
\be
   P  \;=\;
   \begin{bmatrix}
      \alpha_1          & 1                   &     &   &      \\
      \alpha_1 \alpha_2 & \alpha_2 + \alpha_3 & 1   &   &      \\
                        & \alpha_3 \alpha_4   & \alpha_4 + \alpha_5 & 1   &  \\
                        &                     & \ddots & \ddots & \ddots
   \end{bmatrix}
 \label{def.Snl.production.0}
\ee
(see eq.~\reff{def.Snl.production} in Appendix~\ref{subsec.SR}),
which has as its output matrix the triangular array of
generalized Stieltjes--Rogers polynomials of the first kind.
And the production matrix corresponding to the
coefficients $\alpha_{2k-1} = \lambda+k-1$ and $\alpha_{2k} = k$
is precisely \reff{eq.prop.prodmat.univariate.NEW.1.BIS}.
In view of Proposition~\ref{prop.prodmat.univariate.NEW}(a), this shows:

\begin{proposition}
   \label{prop.L.SR}
The Laguerre coefficient matrix $\sfL^{(\alpha)}$
is the matrix $\sfS = (S_{n,\ell}(\balpha))_{n,\ell \ge 0}$
of generalized Stieltjes--Rogers polynomials of the first kind
corresponding to the coefficients
$\alpha_{2k-1} = k+\alpha$ and $\alpha_{2k} = k$.
\end{proposition}

{\bf Remark.}
Proposition~\ref{prop.L.SR} looks at first sight a bit bizarre,
because the coefficient matrix for the monic orthogonal polynomials
associated to a Stieltjes moment sequence $\ba = (a_n)_{n \ge 0}$
is in~general the {\em inverse matrix}\/ of the corresponding matrix $\sfS$
of generalized Stieltjes--Rogers polynomials\footnote{
   This is a well-known fact:  see, for instance
   \cite[p.~III-2, Th\'eor\`eme~1]{Viennot_83}.
   And for a generalization to unit-lower-Hessenberg production matrices
   that are not necessarily tridiagonal,
   see \cite[Proposition~3.2 and Remark~1 following it]{Sokal_multiple_OP}.
};
therefore, the coefficient matrix for the monic unsigned Laguerre polynomials
should be the {\em unsigned inverse matrix}\/ of $\sfS$, i.e.
\be
   \sfS^\sharp  \;\eqdef\; Q \, \sfS^{-1} \, Q
\ee
where $Q = ((-1)^i \delta_{ij})_{i,j \ge 0}$
is the diagonal matrix of alternating 1's and $-1$'s.
This apparent paradox is explained by the curious fact that the
matrix $\sfL^{(\alpha)}$ equals its own unsigned inverse:
that is,
\be
   \sfL^{(\alpha)}  \;=\;  Q \, (\sfL^{(\alpha)})^{-1} \, Q
   \;.
\ee
{\sc Proof.}  From \reff{def.coeffmat} we see that
$\sfL^{(\alpha)} = D \, B^{(\alpha)} \, D^{-1}$
where $B^{(\alpha)}$ is the generalized binomial matrix
\be
   (B^{(\alpha)})_{nk}  \;=\;  \binom{n+\alpha}{n-k}
\ee
and $D = \diag\big( (n!)_{n \ge 0} \big)$.
So it suffices to prove that $B^{(\alpha)}$ equals its own unsigned inverse:
\begin{subeqnarray}
   \sum_j \binom{n+\alpha}{n-j} \, \binom{j+\alpha}{j-k} \, (-1)^{j-k}
   & = &
   \binom{n+\alpha}{n-k} \, \sum_j \binom{n-k}{j-k} \, (-1)^{j-k}  \\[2mm]
   & = &
   \delta_{nk}
\end{subeqnarray}
(cf.~\cite[eq.~(5.21)]{Graham_94}).
\myendremark

\subsection[Binomial row-generating matrix $\sfL^{(\alpha)} B_x$: Proof of Proposition~\ref{prop.prodmat.univariate.NEW}(b)]{Binomial row-generating matrix $\sfL^{\bm{(\alpha)}} \bm{B_x}$: Proof of Proposition~\ref{prop.prodmat.univariate.NEW}(b)}
   \label{subsec.prodmat.1var.LBx}

We now prove Proposition~\ref{prop.prodmat.univariate.NEW}(b),
which asserts that the production matrix
for the binomial row-generating matrix $\sfL^{(\alpha)} B_x$
is the quadridiagonal unit-lower-Hessenberg matrix
$P = (p_{ij})_{i,j \ge 0}$ defined by
\begin{subeqnarray}
   p_{n,n+1}  & = &   1   \\
   p_{n,n}    & = &   (2n+1+\alpha) \,+\, x   \\
   p_{n,n-1}  & = &   n(n+\alpha) \,+\, 2nx \\
   p_{n,n-2}  & = &   n(n-1) x  \\
   p_{n,k}    & = &   0 \qquad\textrm{if $k < n-2$ or $k > n+1$}
 \label{eq.prop.prodmat.univariate.NEW.2.BIS}
\end{subeqnarray} 
This is an easy consequence of the theory of exponential Riordan arrays.
In the preceding subsection we saw that the production matrix $P^\circ$
for the Laguerre coefficient matrix $\sfL^{(\alpha)}$
is the exponential AZ matrix $\EAZ(A,Z)$
where $A$ and $Z$ are given by \reff{eq.prodmat1.A}/\reff{eq.prodmat1.Z}.
By Lemma~\ref{lemma.production.AB}
the production matrix for $\sfL^{(\alpha)} B_x$ is $P = B_x^{-1} P^\circ B_x$;
and by Lemma~\ref{lemma.BxinvEAZBx} we have
$B_x^{-1} P^\circ B_x = \EAZ(A,Z+xA)$.
Using this together with \reff{eq.prodmat1.A}/\reff{eq.prodmat1.Z}
and \reff{def.EAZ.1}
yields \reff{eq.prop.prodmat.univariate.NEW.2.BIS}.

In Appendix~\ref{app.banded} we will explain why $P = B_x^{-1} P^\circ B_x$
is quadridiagonal, by answering the more general question:
Which lower-Hessenberg matrices $P$ have the property that
$B_\xi^{-1} P B_\xi$ is $(r,1)$-banded?

\bigskip

We also have an explicit formula for the binomial row-generating matrix:

\begin{proposition}[Binomial row-generating matrix of the Laguerre coefficient array]
   \label{prop.lagcoeff.Bx.NEW}
The matrix elements of the binomial row-generating matrix
$\sfL^{(\alpha)} B_x$ are
\be
   (\sfL^{(\alpha)} B_x)_{n,k}
   \;=\;
   {1 \over k!} \, {d^k \over dx^k} \, \scrlna(x)
   \;=\;
   \binom{n}{k} \, \scrl_{n-k}^{(\alpha+k)}(x)
   \;.
 \label{eq.prop.lagcoeff.Bx.NEW}
\ee
\end{proposition}

\proof
The first equality is a special case of the general property \reff{def.Ank}
of binomial row-generating polynomials.
The second equality follows by repeated use of the differentiation formula
\be
   {d \over dx} \scrl_n^{(\alpha)}(x)
   \;=\;
   n \scrl_{n-1}^{(\alpha+1)}(x)
   \;.
\ee
\qed
   
It can also be shown, by direct computation using classical identities
for the Laguerre polynomials,
that the matrix \reff{eq.prop.lagcoeff.Bx.NEW}
is indeed the output matrix for the production matrix
\reff{eq.prop.prodmat.univariate.NEW.2.BIS}.
But since this proof is rather lengthy, we refrain from showing it here.

\subsection[Total positivity of the coefficient matrix $\sfL^{(-1+\lambda)}$]{Total positivity of the coefficient matrix $\bm{\sfL^{(-1+\lambda)}}$}
   \label{subsec.univariate_coeff_TP}

Finally, let us give a simple direct proof of the
coefficientwise total positivity of the Laguerre
coefficient matrix $\sfL^{(-1+\lambda)}$.
Start from the binomial matrix $a_{nk} = \binom{n}{k}$,
which is totally positive by Lemma~\ref{lemma.binomialmatrix.TP}.
Now use Lemma~\ref{lemma.diagmult.TP}(b) with
$x_i = \lambda + i-1$ for $i \ge 1$,
which are nonnegative elements of $\Z[\lambda]$.
Then
\begin{subeqnarray}
   b_{nk}  & = & x_{k+1} x_{k+2} \,\cdots\, x_n \, a_{nk}
      \\[2mm]
   & = &  (\lambda+k)^{\overline{n-k}} \, \binom{n}{k}
   \;,
\end{subeqnarray}
which is \reff{def.coeffmat}  with $\alpha = -1+\lambda$.

\section{Multivariate Laguerre polynomials: \\
         Exponential generating functions and production matrices}
    \label{sec.egf}

In this section we compute the
column and bivariate exponential generating functions
for the multivariate Laguerre coefficient matrix \reff{def.coeffmat.gen.second};
we then use these exponential generating functions,
together with the theory of exponential Riordan arrays,
to deduce the corresponding production matrices.

\subsection{Exponential generating functions}  \label{subsec.egf.egf}

In this subsection we shall compute the
column and bivariate exponential generating functions
for the coefficient matrix
$\sfLhat^{(\alpha)}(\yp,\yv,\yda,\ydd,\yfp)$
that was defined in \reff{def.coeffmat.gen.second},
which enumerates Laguerre digraphs according to the status of vertices
as peaks, valleys, double ascents, double descents or fixed points.
Indeed, we shall go farther, by assigning different weights
for the vertices belonging to a cycle or to a path.
For a Laguerre digraph $G$,
let us write $\pkcyc(G)$, $\valcyc(G)$, $\dacyc(G)$, $\ddcyc(G)$, $\fp(G)$
for the number of peaks, valleys, double ascents, double descents
and fixed points that belong to a cycle of $G$,
and $\pkpa(G)$, $\valpa(G)$, $\dapa(G)$, $\ddpa(G)$
for the number of peaks, valleys, double ascents and double descents
that belong to a path of $G$
(of course fixed points can only belong to a cycle).
We then assign weights $\yp,\yv,\yda,\ydd,\yfp$
to the vertices belonging to a cycle,
and weights $\zp,\zv,\zda,\zdd$
to the vertices belonging to a path.
We therefore define the
\textbfit{generalized second multivariate Laguerre coefficient matrix}
\begin{eqnarray}
   & &
   \!\!\!\!
   \sfLtilde^{(\alpha)}(\yp,\yv,\yda,\ydd,\yfp, \zp,\zv,\zda,\zdd)_{n,k}
         \nonumber \\[3mm]
   & &
   \quad\eqdef
   \sum_{G \in \LD_{n,k}}
       y_{\rm p}^{\pkcyc(G)} y_{\rm v}^{\valcyc(G)}
       y_{\rm da}^{\dacyc(G)} y_{\rm dd}^{\ddcyc(G)} y_{\rm fp}^{\fp(G)}
       \,
       z_{\rm p}^{\pkpa(G)} z_{\rm v}^{\valpa(G)}
       z_{\rm da}^{\dapa(G)} z_{\rm dd}^{\ddpa(G)}
         \:\times
         \nonumber \\[-1mm]
   & & \qquad\qquad\quad\;\,
            \, (1+\alpha)^{\cyc(G)}
   \:.
   \qquad
 \label{def.coeffmat.gen.second.gen}
\end{eqnarray}
This polynomial is homogeneous of degree $n$ in
$\yp,\yv,\yda,\ydd,\yfp,\zp,\zv,\zda,\zdd$.

Because of the 0--0 boundary conditions,
each path contains at least one peak;
so~we can, if we wish, remove a factor $z_{\rm p}^k$
and define a unit-lower-triangular matrix by
\be
   \sfLtilde^{(\alpha)\flat}(\yp,\yv,\yda,\ydd,\yfp, \zp,\zv,\zda,\zdd)_{n,k}
   \;\eqdef\;
   \sfLtilde^{(\alpha)}(\yp,\yv,\yda,\ydd,\yfp, \zp,\zv,\zda,\zdd)_{n,k} / z_{\rm p}^k
   \;.
 \label{def.coeffmat.gen.second.gen.flat}
\ee
This polynomial is homogeneous of degree $n-k$ in
$\yp,\yv,\yda,\ydd,\yfp,\zp,\zv,\zda,\zdd$.

We now proceed to compute the exponential generating functions
for the matrices \reff{def.coeffmat.gen.second.gen}
and \reff{def.coeffmat.gen.second.gen.flat}.
We do this by combining the known exponential generating functions
for permutations with cyclic statistics \cite[Th\'eor\`eme~1]{Zeng_93}
and permutations with linear statistics \cite[Proposition~4]{Zeng_93}.
These formulae are as follows:

\medskip

1) Permutations with cyclic statistics are enumerated by the polynomials
\be
   P_n^{\rm cyc}(\yp,\yv,\yda,\ydd,\yfp,\lambda)
   \;\eqdef\;
   \sum_{\sigma \in \Sym_n}
      y_{\rm p}^{\pkcyc(\sigma)} y_{\rm v}^{\valcyc(\sigma)}
       y_{\rm da}^{\dacyc(\sigma)} y_{\rm dd}^{\ddcyc(\sigma)}
       y_{\rm fp}^{\fp(\sigma)}
      \, \lambda^{\cyc(\sigma)}
\ee
where $\pkcyc(\sigma)$, $\valcyc(\sigma)$, $\dacyc(\sigma)$, $\ddcyc(\sigma)$,
$\fp(\sigma)$
denote the number of cycle peaks, cycle valleys, cycle double rises,
cycle double falls and fixed points in $\sigma$,
and $\cyc(\sigma)$ denotes the number of cycles in $\sigma$.
By convention we set $P_0^{\rm cyc} = 1$.
The polynomial $P_n^{\rm cyc}$ is homogeneous of degree $n$ in
$\yp,\yv,\yda,\ydd,\yfp$.
We write
\be
   F(t;\yp,\yv,\yda,\ydd,\yfp,\lambda)
   \;\eqdef\;
   \sum_{n=0}^\infty P_n^{\rm cyc}(\yp,\yv,\yda,\ydd,\yfp,\lambda)
       \: {t^n \over n!}
\ee
for the corresponding exponential generating function.

\begin{lemma} {\bf\cite[Th\'eor\`eme~1]{Zeng_93}}
   \label{lemma.zeng.cyclic}
We have
\begin{subeqnarray}
   F(t;\yp,\yv,\yda,\ydd,\yfp,\lambda)
   & = &
   e^{\lambda \yfp t}
   \,
   \biggl( {r_1 \,-\, r_2
            \over
            r_1 e^{r_2 t} \,-\, r_2 e^{r_1 t}
           }
   \biggr)^{\! \lambda}
       \\[2mm]
   & = &
   F(t;\yp,\yv,\yda,\ydd,\yfp,1)^\lambda
 \label{eq.lemma.zeng.cyclic.1}
\end{subeqnarray}
where $r_1 r_2 = \yp \yv$ and $r_1 + r_2 = \yda + \ydd$.
Otherwise put, $r_1$ and $r_2$ are the roots (in~either order)
of the quadratic equation $\rho^2 - (\yda + \ydd) \rho + \yp \yv = 0$.
Concretely,
\be
   r_{1,2}
   \;=\;
   {\yda + \ydd \pm \sqrt{(\yda + \ydd)^2 - 4\yp \yv}
    \over
    2
   }
   \;.
 \label{eq.lemma.zeng.cyclic.2}
\ee
\end{lemma}

\medskip

2) Permutations with linear statistics and 0--0 boundary conditions
are enumerated by the polynomials
\be
   P_n^{\rm lin(00)}(\zp,\zv,\zda,\zdd)
   \;\eqdef\;
   \sum_{\sigma \in \Sym_n}
      z_{\rm p}^{\pkpa(\sigma)} z_{\rm v}^{\valpa(\sigma)}
       z_{\rm da}^{\dapa(\sigma)} z_{\rm dd}^{\ddpa(\sigma)}
\ee
where $\pkpa(\sigma)$, $\valpa(\sigma)$, $\dapa(\sigma)$, $\ddpa(\sigma)$
denote the number of peaks, valleys, double ascents and double descents
in the permutation $\sigma$ written as a word $\sigma_1 \cdots \sigma_n$,
where we impose the boundary conditions $\sigma_0 = \sigma_{n+1} = 0$.
By convention we restrict attention to $n \ge 1$.
The polynomial $P_n^{\rm lin(00)}$ is homogeneous of degree $n$ in
$\zp,\zv,\zda,\zdd$.
We write
\be
   G(t;\zp,\zv,\zda,\zdd)
   \;\eqdef\;
   \sum_{n=1}^\infty P_n^{\rm lin(00)}(\zp,\zv,\zda,\zdd)
       \: {t^n \over n!}
\ee
for the corresponding exponential generating function
(note that the sum starts at $n=1$).

\begin{lemma}  {\bf\cite[Proposition~4]{Zeng_93}}
   \label{lemma.zeng.linear}
We have
\be
   G(t;\zp,\zv,\zda,\zdd)
   \;=\;
   \zp
   \:
   \biggl( {e^{r_1 t} \,-\, e^{r_2 t}
            \over
            r_1 e^{r_2 t} \,-\, r_2 e^{r_1 t}
           }
   \biggr)
 \label{eq.lemma.zeng.linear.1}
\ee
where $r_1 r_2 = \zp \zv$ and $r_1 + r_2 = \zda + \zdd$
analogously to Lemma~\ref{lemma.zeng.cyclic}.
This function satisfies the differential equation
\be
   G'(t)  \;=\; \zp \,+\, (\zda+\zdd) G(t) \,+\, \zv G(t)^2
   \;.
 \label{eq.lemma.zeng.linear.2}
\ee
\end{lemma}

We can now put these ingredients together to determine
the exponential generating functions
for the matrix \reff{def.coeffmat.gen.second.gen}.
A Laguerre digraph $G \in \LD_{n,k}$
consists of a permutation (that is, a collection of disjoint cycles)
on some subset $S \subseteq [n]$
together with $k$ disjoint paths on $[n] \setminus S$.
Each of these paths can be considered as a permutation written in word form.
By the exponential formula, the exponential generating function
for the $k$th column of the matrix \reff{def.coeffmat.gen.second.gen} is then
\begin{eqnarray}
   & &
   \sum_{n=0}^\infty
      \sfLtilde^{(\alpha)}(\yp,\yv,\yda,\ydd,\yfp, \zp,\zv,\zda,\zdd)_{n,k}
      \: {t^n \over n!}
             \nonumber \\[3mm]
   & &
   \qquad =\;
   F(t;\yp,\yv,\yda,\ydd,\yfp,1)^{1+\alpha}
   \;
   {G(t;\zp,\zv,\zda,\zdd)^k \over k!}
   \;,
   \qquad\qquad
\end{eqnarray}
where the $1/k!$ comes because the paths are indistinguishable.
The bivariate exponential generating function is therefore
\begin{eqnarray}
   & &
   \sum_{n=0}^\infty \sum_{k=0}^n
      \sfLtilde^{(\alpha)}(\yp,\yv,\yda,\ydd,\yfp, \zp,\zv,\zda,\zdd)_{n,k}
      \: {t^n \over n!} \, u^k
             \nonumber \\[3mm]
   & &
   \qquad =\;
   F(t;\yp,\yv,\yda,\ydd,\yfp,1)^{1+\alpha}
   \;
   \exp\bigl[ u \, G(t;\zp,\zv,\zda,\zdd) \bigr]
   \;.
   \qquad\qquad
\end{eqnarray}
Comparing this with \reff{def.exp_riordan}/\reff{eq.exp_riordan.bivariate_egf}
shows:

\begin{proposition}[Generalized second multivariate Laguerre coefficient matrix as exponential Riordan array]
The matrix
$\sfLtilde^{(\alpha)}(\yp,\yv,\yda,\ydd,\yfp, \zp,\zv,\zda,\zdd)$,
which was defined in \reff{def.coeffmat.gen.second.gen},
is an exponential Riordan array $\scrr[F,G]$
where $F$ is given by \reff{eq.lemma.zeng.cyclic.1} with $\lambda = 1+\alpha$
and $G$ is given by \reff{eq.lemma.zeng.linear.1}.
\end{proposition}

For the matrix $\sfLtilde^{(\alpha)\flat}$
defined in \reff{def.coeffmat.gen.second.gen.flat},
the formulae are the same except that $G$ is replaced by
\be
   G^\flat(t;\zp,\zv,\zda,\zdd)
   \;\eqdef\;
   G(t;\zp,\zv,\zda,\zdd) /\zp
   \;=\;
    {e^{r_1 t} \,-\, e^{r_2 t}
            \over
            r_1 e^{r_2 t} \,-\, r_2 e^{r_1 t}
           }
   \;,
 \label{eq.lemma.zeng.linear.1.flat}
\ee
which satisfies the differential equation
\be
   (G^\flat)'(t)  \;=\; 1 \,+\, (\zda+\zdd) G^\flat(t) \,+\, \zp\zv G^\flat(t)^2
   \;.
 \label{eq.lemma.zeng.linear.2.flat}
\ee

\bigskip

{\bf Remark.}
For the first multivariate coefficient matrix \reff{def.coeffmat.gen},
the formulae simplify:
taking $\yp = \ydd = \zp = \zdd = v_-$,
$\yv = \yda = \zv = \zda = v_-$ and $\yfp = v_0$,
we~get $r_1 = v_+$, $r_2 = v_-$ and hence,
using \reff{eq.second-to-first.1},
\begin{eqnarray}
   & &
   \sum_{n=0}^\infty \sum_{k=0}^n
      \sfL^{(\alpha)}(v_-, v_0, v_+)_{n,k}
      \: {t^n \over n!} \, u^k
             \nonumber \\[3mm]
   & &
   \qquad =\;
   e^{(1+\alpha) v_0 t}
   \:
   \biggl( {v_+ \,-\, v_-
            \over
            v_+ e^{v_- t} \,-\, v_- e^{v_+ t}
           }
   \biggr)^{\! 1 + \alpha}
   \:
   \exp\Biggl[ u \,
   \biggl( {e^{v_+ t} \,-\, e^{v_- t}
            \over
            v_+ e^{v_- t} \,-\, v_- e^{v_+ t}
           }
   \biggr)
   \Biggr]
   \;.
   \qquad\qquad
  \label{eq.first_multivariate.egf}
\end{eqnarray}

\subsection{Production matrices: Proof of Proposition~\ref{prop.prodmat.multivariate.NEW}(a,b)}  \label{subsec.egf.prodmat}

We now restrict attention to the case in which the vertices belonging
to paths and cycles are given the same weights,
i.e.\ $(\zp,\zv,\zda,\zdd) = (\yp,\yv,\yda,\ydd)$.
In this situation it is easy to find the series $A(s)$ and $Z(s)$
that satisfy the relations \reff{eq.prop.riordan.exponential.production.1}
leading to the production matrix of an exponential Riordan array.
Namely, from \reff{eq.lemma.zeng.linear.2} we have immediately
\be
   A(s)  \;=\; \yp \,+\, (\yda+\ydd) s \,+\, \yv s^2
   \;.
 \label{eq.multivariate_second.A}
\ee
And a straightforward computation shows that
the logarithmic derivative of $F$ equals $\lambda \yfp + \lambda \yv G$,
so we have
\be
   Z(s)  \;=\;  \lambda \yfp + \lambda \yv s
   \;.
 \label{eq.multivariate_second.Z}
\ee
Theorem~\ref{thm.riordan.exponential.production} then yields:

\begin{proposition}[Production matrix of the second multivariate Laguerre coefficient matrix]
   \label{prop.multivariate_second.prodmat}
The second multivariate Laguerre coefficient matrix
$\sfLhat^{(\alpha)}(\yp,\yv,\yda,\ydd,\yfp)$,
which was defined in \reff{def.coeffmat.gen.second},
has production matrix $P^\circ = \EAZ(A,Z)$ where $A$ and $Z$ are given by
\reff{eq.multivariate_second.A}/\reff{eq.multivariate_second.Z}
with $\lambda = 1 + \alpha$.
In detail, $P^\circ$ is the tridiagonal matrix
\begin{subeqnarray}
   p^\circ_{n,n+1}  & = &   \yp   \\[1mm]
   p^\circ_{n,n}    & = &   (1+\alpha)\yfp \,+\, n(\yda+\ydd)   \\[1mm]
   p^\circ_{n,n-1}  & = &   n(n+\alpha)\yv   \\[1mm]
   p^\circ_{n,k}    & = &   0 \qquad\textrm{if $k < n-1$ or $k > n+1$}
 \label{eq.prop.multivariate_second.prodmat}
\end{subeqnarray}
\end{proposition}

{}From this it is straightforward to deduce the production matrix
for the binomial row-generating matrix
$\sfLhat^{(\alpha)}(\yp,\yv,\yda,\ydd,\yfp) \, B_x$.
Namely, by Lemma~\ref{lemma.production.AB}
this production matrix is $P = B_x^{-1} P^\circ B_x$;
and by Lemma~\ref{lemma.BxinvEAZBx} we have
$B_x^{-1} P^\circ B_x = \EAZ(A,Z+xA)$.
We have therefore proven:

\begin{proposition}[Production matrix of the second multivariate Laguerre binomial row-generating matrix]
   \label{prop.multivariate_second.rowgen.prodmat}
The binomial row-generating matrix
$\sfLhat^{(\alpha)}(\yp,\yv,\yda,\ydd,\yfp) \, B_x$
has production matrix $P = \EAZ(A,Z+xA)$ where $A$ and $Z$ are given by
\reff{eq.multivariate_second.A}/\reff{eq.multivariate_second.Z}
with $\lambda = 1 + \alpha$.
In detail, $P$ is the quadridiagonal lower-Hessenberg matrix
\begin{subeqnarray}
   p_{n,n+1}  & = &   \yp   \\[1mm]
   p_{n,n}    & = &   (1+\alpha)\yfp \,+\, n(\yda+\ydd) \,+\, \yp x   \\[1mm]
   p_{n,n-1}  & = &   n(n+\alpha)\yv \,+\, n(\yda+\ydd) x   \\[1mm]
   p_{n,n-2}  & = &   n(n-1)\yv x   \\[1mm]
   p_{n,k}    & = &   0 \qquad\textrm{if $k < n-2$ or $k > n+1$}
 \label{eq.prop.multivariate_second.rowgen.prodmat}
\end{subeqnarray}
\end{proposition}

For the matrix $\sfLtilde^{(\alpha)\flat}$,
the formulae are the same except that $A$ and $Z$ are replaced by
\begin{eqnarray}
   A^\flat(s)  & = &  1 \,+\, (\yda+\ydd) s \,+\, \yp\yv s^2
       \label{eq.multivariate_second.A.flat} \\[2mm]
   Z^\flat(s)  & = &  \lambda \yfp + \lambda \yp\yv s
       \label{eq.multivariate_second.Z.flat}
\end{eqnarray}
This leads to the production matrices
\begin{subeqnarray}
   p^{\circ\flat}_{n,n+1}  & = &   1   \\[1mm]
   p^{\circ\flat}_{n,n}    & = &   (1+\alpha)\yfp \,+\, n(\yda+\ydd)   \\[1mm]
   p^{\circ\flat}_{n,n-1}  & = &   n(n+\alpha)\yp\yv   \\[1mm]
   p^{\circ\flat}_{n,k}    & = &   0 \qquad\textrm{if $k < n-1$ or $k > n+1$}
 \label{eq.prop.multivariate_second.prodmat.flat}
\end{subeqnarray}
and
\begin{subeqnarray}
   p^\flat_{n,n+1}  & = &   1   \\[1mm]
   p^\flat_{n,n}    & = &   (1+\alpha)\yfp \,+\, n(\yda+\ydd) \,+\, x \\[1mm]
   p^\flat_{n,n-1}  & = &   n(n+\alpha)\yp\yv \,+\, n(\yda+\ydd) x   \\[1mm]
   p^\flat_{n,n-2}  & = &   n(n-1) \yp\yv x   \\[1mm]
   p^\flat_{n,k}    & = &   0 \qquad\textrm{if $k < n-2$ or $k > n+1$}
 \label{eq.prop.multivariate_second.rowgen.prodmat.flat}
\end{subeqnarray}
This completes the proof of Proposition~\ref{prop.prodmat.multivariate.NEW}.

\section[Total positivity of the production matrices: Univariate case]{Total positivity of the production matrices: \\ Univariate case}
                   \label{sec.prodmat.TP.univariate}

In this section we will prove the total positivity
of the production matrices for the univariate Laguerre polynomials:
namely, the tridiagonal production matrix $P^\circ$
[defined in \reff{eq.prop.prodmat.univariate.NEW.1}]
for the Laguerre coefficient matrix $\sfL^{(\alpha)}$,
and the quadridiagonal production matrix $P$
[defined in \reff{eq.prop.prodmat.univariate.NEW.2}]
for the binomial row-generating matrix $\sfL^{(\alpha)} \, B_x$.
Throughout, we use the parameter $\lambda = 1 + \alpha$.

The proofs are easy factorizations;
in the quadridiagonal case we will also need to invoke the
tridiagonal comparison theorem (Proposition~\ref{prop.comparison}).

\subsection[The tridiagonal matrix $P^\circ$: Proof of Proposition~\ref{prop.prodmat.TP.univariate.NEW}(a)]{The tridiagonal matrix $\bm{P^\circ}$: Proof of Proposition~\ref{prop.prodmat.TP.univariate.NEW}(a)}

The tridiagonal matrix $P^\circ$ is given by
\begin{subeqnarray}
   p^\circ_{n,n+1}  & = &   1   \\
   p^\circ_{n,n}    & = &   2n+\lambda   \\
   p^\circ_{n,n-1}  & = &   n(n-1+\lambda)   \\
   p^\circ_{n,k}    & = &   0 \qquad\textrm{if $k < n-1$ or $k > n+1$}
 \label{eq.prop.prodmat.univariate.NEW.1.BIS2}
\end{subeqnarray} 
As observed preceding Proposition~\ref{prop.L.SR},
this is the production matrix \reff{def.Snl.production}
corresponding to the classical S-fraction with coefficients
$\alpha_{2k-1} = k-1+\lambda$ and $\alpha_{2k} = k$
(see Appendix~\ref{subsec.SR}).
And as observed in \reff{def.Snl.production.factorization},
the production matrix \reff{def.Snl.production}
can be factorized in the form $LU$ where
$L$ is the lower-bidiagonal matrix with 1 on the diagonal
and $\alpha_2,\alpha_4,\ldots$ on the subdiagonal,
and $U$ is the upper-bidiagonal matrix with 1 on the superdiagonal
and $\alpha_1,\alpha_3,\ldots$ on the diagonal.
And finally, as also observed there,
this factorization shows, by Lemma~\ref{lemma.bidiagonal},
that the matrix \reff{def.Snl.production} is coefficientwise totally positive
in the indeterminates $\balpha$.
Putting this all together, we conclude:

\begin{proposition}[Factorization and total positivity of $P^\circ$]
   \label{prop.factorization.Pcirc}
\hfill\break\noindent
\vspace*{-6mm}
\begin{itemize}
   \item[(a)]
The tridiagonal matrix $P^\circ$ defined by
\reff{eq.prop.prodmat.univariate.NEW.1.BIS2}
has the factorization $P^\circ = LU$ where
\begin{itemize}
   \item[$\bullet$] $L$ is the lower-bidiagonal matrix with 1 on the diagonal
      and $1,2,3,\ldots$ on the subdiagonal, and
   \item[$\bullet$] $U$ is the upper-bidiagonal matrix with 1 on the
      superdiagonal and $\lambda, {\lambda+1},$ ${\lambda+2}, \ldots$ on the diagonal. 
\end{itemize}
   \item[(b)] $P^\circ$ is totally positive in the ring $\Z[\lambda]$
equipped with the coefficientwise order.
\end{itemize}
\end{proposition}

\noindent
In particular, this proves Proposition~\ref{prop.prodmat.TP.univariate.NEW}(a).

Note that we can also write $U = L\Delta + \lambda I$,
since $L\Delta$ is the upper-bidiagonal matrix with 1~on the superdiagonal
and $0,1,2,\ldots$ on the diagonal.

\subsection[The quadridiagonal matrix $P$: Proof of Proposition~\ref{prop.prodmat.TP.univariate.NEW}(b)]{The quadridiagonal matrix $\bm{P}$: Proof of Proposition~\ref{prop.prodmat.TP.univariate.NEW}(b)}

The quadridiagonal matrix $P$ is given by
\begin{subeqnarray}
   p_{n,n+1}  & = &   1   \\
   p_{n,n}    & = &   (2n+\lambda) \,+\, x   \\
   p_{n,n-1}  & = &   n(n-1+\lambda) \,+\, 2nx \\
   p_{n,n-2}  & = &   n(n-1) x  \\
   p_{n,k}    & = &   0 \qquad\textrm{if $k < n-2$ or $k > n+1$}
 \label{eq.prop.prodmat.univariate.NEW.2.BIS2}
\end{subeqnarray} 

\begin{proposition}[Factorization and total positivity of $P$]
   \label{prop.factorization.P}
\hfill\break\noindent
\vspace*{-6mm}
\begin{itemize}
   \item[(a)]
The quadridiagonal matrix $P$ defined by
\reff{eq.prop.prodmat.univariate.NEW.2.BIS2}
has the factorization $P = {L(L U_x + \lambda I)}$ where
\begin{itemize}
   \item[$\bullet$] $L$ is the lower-bidiagonal matrix with 1 on the diagonal
      and $1,2,3,\ldots$ on the subdiagonal, and
   \item[$\bullet$] $U_x = \Delta + xI$ is the upper-bidiagonal matrix
      with 1 on the superdiagonal and $x$ on the diagonal.
\end{itemize}
   \item[(b)] $P$ is totally positive in the ring $\Z[x,\lambda]$
equipped with the coefficientwise order.
\end{itemize}
\end{proposition}

\proof
(a) It is easy to see that $P = P^\circ + x LL$.
Since $P^\circ = L(L\Delta + \lambda I)$
by Proposition~\ref{prop.factorization.Pcirc}(a),
we have $P = L[L(\Delta + xI) + \lambda I]$,
which proves part~(a).

(b) By Lemma~\ref{lemma.bidiagonal}, the bidiagonal matrices $L$ and $U_x$
are totally positive in the ring $\Z[x]$
equipped with the coefficientwise order.
It follows that the tridiagonal matrix $L U_x$
is also totally positive in the ring $\Z[x]$
equipped with the coefficientwise order.
Then the tridiagonal comparison theorem (Proposition~\ref{prop.comparison})
shows that the tridiagonal matrix $L U_x + \lambda I$
is totally positive in the ring $\Z[x,\lambda]$
equipped with the coefficientwise order.
Therefore the same holds for $P = L(L U_x + \lambda I)$.
\qed

This completes the proof of Proposition~\ref{prop.prodmat.TP.univariate.NEW}(b).

\bigskip

{\bf Remarks.}
1.  For $\lambda=0$ ($\alpha=-1$), we have $P = L L U_x$,
which is the production matrix for the known 2-S-fraction
for the Lah polynomials \cite[Theorem~1.5 with $r=2$]{latpath_lah},
which has coefficients
\begin{subeqnarray}
   \alpha_{3n-1}  & = & x   \\
   \alpha_{3n}    & = & n   \\
   \alpha_{3n+1}  & = & n
 \label{eq.lah.alphas}
\end{subeqnarray}
(For the production matrix associated to a 2-S-fraction,
see \cite[Propositions~7.2 and 8.2 and eqn.~(7.8)]{latpath_SRTR}
or the case $j=0$ of Proposition~\ref{prop.Smj.prodmat} below.)
However, it turns out that this same production matrix
has a one-parameter family of distinct factorizations $P = L_1 L_2 U_x$,
arising from the 2-S-fractions with
\begin{subeqnarray}
   \alpha_{3n-1}  & = & x   \\
   \alpha_{3n}    & = & c_n \, n   \\
   \alpha_{3n+1}  & = & (2 - c_n) \, n
\end{subeqnarray}
where
\be
   c_n  \;\eqdef\;  {(n-1) \,-\, (n-2)\kappa  \over  n \,-\, (n-1)\kappa}
   \quad\hbox{for } n \ge 1
\ee
and $\kappa \in [0,1]$:
see Proposition~\ref{prop.Sfrac.j=0.solutions} below.
The case $\kappa=1$ is the known 2-S-fraction;
all the others appear to be new.

2. It is easy to see that $L \Delta + I = \Delta L$
(this is a special case of \reff{eq.delta.shift} below)
and hence $L U_x + I = U_x L$.
This implies that for $\lambda = 1$ ($\alpha=0$), we have
$P = L(L U_x + I) = L U_x L$.
This is the production matrix for the
generalized 2-Stieltjes--Rogers polynomials of type $j=1$
with the coefficients \reff{eq.lah.alphas}:
see Proposition~\ref{prop.Smj.prodmat} for the general theory,
and Proposition~\ref{prop.Sfrac.j=1.solutions.alpha=0} for this application.
\myendremark

\section[Total positivity of the production matrices: Multivariate case]{Total positivity of the production matrices: \\ Multivariate case}
     \label{sec.prodmat.TP.multivariate}

In this section we will prove
Proposition~\ref{prop.prodmat.TP.multivariate.NEW}
on the total positivity of the production matrices
for the multivariate Laguerre polynomials.
The proof of Proposition~\ref{prop.prodmat.TP.multivariate.NEW}(a),
dealing with the tridiagonal matrix $P^{\circ\flat}$,
is an easy factorization
combined with an appeal to the tridiagonal comparison theorem.
By contrast,
the proof of Proposition~\ref{prop.prodmat.TP.multivariate.NEW}(b),
dealing with the quadridiagonal matrix $P^\flat$,
is decidedly nontrivial.

\subsection[The tridiagonal matrix $P^{\circ\flat}$: Proof of Proposition~\ref{prop.prodmat.TP.multivariate.NEW}(a)]{The tridiagonal matrix $\bm{P^{\circ\flat}}$: Proof of Proposition~\ref{prop.prodmat.TP.multivariate.NEW}(a)}

\proofof{Proposition~\ref{prop.prodmat.TP.multivariate.NEW}{\rm (a)}}
Let $Q$ the matrix \reff{def.Snl.production.factorization}
with $\alpha_{2k-1} = (k+\alpha)\yp$ and $\alpha_{2k} = k\yv$ for $k \ge 1$;
note that these coefficients are nonnegative in the ring~$R$
because of the hypotheses $\alpha \ge -1$, $\yp \ge 0$ and $\yv \ge 0$.
It then follows from Lemma~\ref{lemma.bidiagonal}
that $Q$ is totally positive in the ring $R$.
The matrix $Q$ has elements
\begin{subeqnarray}
   q_{n,n+1}  & = &   1   \\[1mm]
   q_{n,n}    & = &   (1+\alpha)\yp \,+\, n(\yp+\yv)   \\[1mm]
   q_{n,n-1}  & = &   n(n+\alpha)\yp\yv   \\[1mm]
   q_{n,k}    & = &   0 \qquad\textrm{if $k < n-1$ or $k > n+1$}
\end{subeqnarray}
Then $P^{\circ\flat} = Q + D$,
where $D$ is the diagonal matrix with entries
$d_n = {(1+\alpha)(\yfp-\yp)} + {n(\yda+\ydd-\yp-\yv)}$.
By hypothesis these entries are nonnegative in the ring $R$,
so the tridiagonal comparison theorem (Proposition~\ref{prop.comparison})
implies that $P^{\circ\flat}$ is totally positive in the ring $R$.
\qed

\subsection[The quadridiagonal matrix $P^\flat$: Proof of Proposition~\ref{prop.prodmat.TP.multivariate.NEW}(b)]{The quadridiagonal matrix $\bm{P^\flat}$: Proof of Proposition~\ref{prop.prodmat.TP.multivariate.NEW}(b)}
   \label{subsec.prodmat.TP.multivariate.quadridiagonal}


We will now prove Proposition~\ref{prop.prodmat.TP.multivariate.NEW}(b),
which asserts the total positivity
of the production matrix 
\reff{eq.prop.multivariate_second.rowgen.prodmat.flat.BIS0}
when the variables are substituted to elements
of a partially ordered commutative ring $R$
that satisfy $\lambda\geq 0$,
$\lambda \yfp = \lambda \yp$, $\yp \ge 0$, $\yv \ge 0$,
$\yda+\ydd \ge \yp + \yv$ and $x \ge 0$.
We will do this by proving the {\em coefficientwise}\/ total positivity
of a {\em much more general}\/ quadridiagonal lower-Hessenberg matrix $P$,
defined by
\begin{subeqnarray}
   P  & \eqdef &  L_1 U L_2 \, +\, L_1 D_1 \,+\, D_2  L_2
           \slabel{eq.thm.prodmat.TP.bis.gen.new2.a} \\[2mm] 
      & = &    L_1(U L_2 \, +\, D_1) \,+\, D_2  L_2
           \slabel{eq.thm.prodmat.TP.bis.gen.new2.b} \\[2mm] 
      & = &   (L_1  U \,+\, D_2)L_2 \,+\, L_1 D_1 
           \slabel{eq.thm.prodmat.TP.bis.gen.new2.c}
   \label{eq.thm.prodmat.TP.bis.gen.new2}
\end{subeqnarray}
where
\begin{itemize}
	\item $L_1$ is the lower-bidiagonal matrix with 
		the sequence $a_0,a_1,\ldots$ on the diagonal, 
		the sequence $b_1, b_2, \ldots$ on the subdiagonal,
		and zeroes elsewhere;
	\item $U$ is the upper-bidiagonal matrix with 
		the sequence $c_1, c_2,\ldots$ on the superdiagonal,
		the sequence $d_0, d_1, \ldots$ on the diagonal,
		and zeroes elsewhere;
	\item $L_2$ is the lower-bidiagonal matrix with 
		the sequence $e_0,e_1,\ldots$ on the diagonal,
                the sequence $f_1, f_2, \ldots$ on the subdiagonal,
                and zeroes elsewhere;
	\item $D_1$ is the diagonal matrix with entries $g_0, g_1,\ldots\,$;
	\item $D_2$ is the diagonal matrix with entries $h_0, h_1,\ldots\,$;
\end{itemize}
and $\bfa = (a_n)_{n\geq 0}$, $\bfb = (b_n)_{n\geq 1}$, 
              $\bfc = (c_n)_{n\geq 1}$, $\bfd = (d_n)_{n\geq 0}$,
	      $\bfe = (e_n)_{n\geq 0}$, $\bfff = (f_n)_{n\geq 1}$,
	      $\bfg = (g_n)_{n\geq 0}$, and $\bfh = (h_n)_{n\geq 0}$
	      are all indeterminates.
The entries in the $n$th row of $P = (p_{n,k})_{k\geq 0}$
are given by
\begin{subeqnarray}
   p_{n,n+1}  & = &   a_n c_{n+1} e_{n+1}   \\[1mm]
   p_{n,n}    & = &  
       a_n d_n e_n + b_n c_n e_n + a_n c_{n+1} f_{n+1} + a_n g_n + h_n e_n \\[1mm]
   p_{n,n-1}  & = &
	a_n d_n f_n + b_n c_n f_n + b_n d_{n-1} e_{n-1} + b_n g_{n-1} + h_n f_n \\[1mm]                                                                 
   p_{n,n-2}  & = &   b_{n} d_{n-1} f_{n-1}   \\[1mm]
   p_{n,k}    & = &   0 \qquad\textrm{if $k < n-2$ or $k > n+1$}
  \label{eq.thm.prodmat.TP.bis.gen.new2.2}
\end{subeqnarray}
where by definition $b_0 = c_0 = f_0 = 0$
and $a_n = b_n = c_n = d_n = e_n = f_n = g_n = h_n = 0$ whenever $n < 0$.
Our main result is:

\begin{theorem}[Total positivity of the generalized production matrix]
\label{thm.prodmat.TP.bis.gen.new2} 
The matrix $P$ defined by 
\reff{eq.thm.prodmat.TP.bis.gen.new2}/\reff{eq.thm.prodmat.TP.bis.gen.new2.2}
is totally positive, coefficientwise in the indeterminates
$\bfa, \bfb, \bfc, \bfd, \bfe,\bfff,\bfg,\bfh$.
\end{theorem}

\begin{samepage}
\proofof{Proposition~\ref{prop.prodmat.TP.multivariate.NEW}(b)
   assuming Theorem~\ref{thm.prodmat.TP.bis.gen.new2}}
Specialize the matrix \reff{eq.thm.prodmat.TP.bis.gen.new2.2} by setting
\begin{subeqnarray}
   & & a_n \;=\; c_n \;=\; e_n \;=\; 1   \\
   & & b_n \;=\; n \yv \\
   & & d_n \;=\; n \yp \\
   & & f_n \;=\; x  \\ 
   & & g_n \;=\; \lambda\yp \: (= \lambda \yfp)  \\
   & & h_n \;=\; n (\yda + \ydd - \yp - \yv)
\end{subeqnarray}
\qed
\end{samepage}

Theorem~\ref{thm.prodmat.TP.bis.gen.new2} will be proven as follows.
Define the matrix $Q = (q_{n,k})_{n,k\geq 0}$ by
\be
   Q  \;\eqdef\; P \big|_{\bfh = \bzero}  \;=\; L_1(U L_2 \, +\, D_1)
   \;,
 \label{def.prodmat.Q.new2}
\ee
or in other words
\begin{subeqnarray}
   q_{n,n+1}  & = &   a_n c_{n+1} e_{n+1}   \\[1mm]
   q_{n,n}    & = &  
       a_n d_n e_n + b_n c_n e_n + a_n c_{n+1} f_{n+1} + a_n g_n \\[1mm]
   q_{n,n-1}  & = &
	a_n d_n f_n + b_n c_n f_n + b_n d_{n-1} e_{n-1} + b_n g_{n-1} \\[1mm]                                                                 
   q_{n,n-2}  & = &   b_{n} d_{n-1} f_{n-1}   \\[1mm]
   q_{n,k}    & = &   0 \qquad\textrm{if $k < n-2$ or $k > n+1$}
   \label{def.Qmat2}
\end{subeqnarray}
Then
\be
P \;=\; Q \, +\, D_2L_2.
\label{eq.Q.P.difference}
\ee
We will begin by proving (Lemma~\ref{lemma.TP.Q.new2})
that $Q$ is coefficientwise totally positive;
this proof uses the factorization $Q = L_1(U L_2 +D_1)$
together with the tridiagonal comparison theorem,
in close analogy with Propositions~\ref{prop.prodmat.TP.multivariate.NEW}(a)
and \ref{prop.factorization.P}.
It follows that for every integer $m \ge 0$,
the matrix consisting of first $m$ rows of $Q$ is also
coefficientwise totally positive.

The rest of the proof shows how to restore the terms in $P$ involving~$\bfh$.
In terms of the row vectors
$(\bp_n)_{n\geq 0}$, $(\bq_n)_{n\geq 0}$, $(\bll_n)_{n\geq 0}$
associated to the matrices $P,Q,L_2$,
equation~\reff{eq.Q.P.difference} can be rewritten as
\be
   \bp_n \;=\; \bq_n \,+\, h_n \bll_n
   \;,
\ee
where
\begin{subeqnarray}
  \bll_0
	& = &
   \begin{bmatrix}
	  e_0 & \bzero_{1\times \infty}
   \end{bmatrix}  \\[1mm]
   \bll_n
	&=&
   \begin{bmatrix}
	   \bzero_{1 \times n-1} & f_n  & e_n  & \bzero_{1\times \infty}  
   \end{bmatrix}
   \quad\textrm{for $n \ge 1$}
\end{subeqnarray}

\medskip

{\bf Remark.}
Here we have chosen to write
$Q \eqdef P |_{\bfh = \bzero} = L_1(U L_2 \, +\, D_1)$
and $P = Q + D_2 L_2$ and to argue by rows
(because $D_2$ acts on the left).
We could equally well have started from
$Q' \eqdef P |_{\bfg = \bzero} = (L_1  U \,+\, D_2)L_2$
and $P = Q' + L_1 D_1$ and to argue by columns
(because $D_1$ acts on the right).
The reader is invited to work out this variant proof.
\myendremark

\medskip

Our proof is based on considering the matrix
$
\displaystyle
\begin{bmatrix}
	\bq_0\\
	\vdots\\
	\bq_{n-1}\\
	\bp_{n}\\
	\vdots\\
	\bp_m
\end{bmatrix}
$,
which we shall sometimes write henceforth for typographical simplicity
as $[[\bq_0,\ldots,\bq_{n-1},\bp_n,\ldots,\bp_m]]$
(and likewise for other matrices written by rows).
We will show (Lemma~\ref{lemma.TP.QPfinal.new2})
that for every pair of integers $0 \le n \le m+1$,
this matrix
is totally positive;
and we will do this, for each fixed $m \ge 0$,
by induction on $n = m+1,m,m-1,\ldots,0$.
The base case $n=m+1$ of this induction is thus Lemma~\ref{lemma.TP.Q.new2},
and the final case $n=0$ is Theorem~\ref{thm.prodmat.TP.bis.gen.new2}.
The proof of Lemma~\ref{lemma.TP.QPfinal.new2} will involve the following steps:

\bigskip

   {\bf Lemma~\ref{lemma.TP.QR.new2}:}
   The matrix $[[\bq_0, \ldots,\bq_{n-1}, \bll_n]]$
       is totally positive.

\bigskip

   \hangindent=1.5cm
   {\bf Lemma~\ref{lemma.TP.QPR.if.new2}:}
       If the matrix $[[\bp_{n+1},\ldots,\bp_m]]$
       is totally positive,\\
       then so is 
       $[[\bq_0,\ldots,\bq_{n-1},\bll_n,\bp_{n+1},\ldots,\bp_m]]$.

\bigskip

   \hangindent=1.5cm
   {\bf The induction step (Lemma~\ref{lemma.TP.QP.inductive.new2}):}
       If the matrix 
       $[[\bq_0,\ldots,\bq_{n},\bp_{n+1},\ldots,\bp_m]]$
       is totally positive,
       then so is 
       $[[\bq_0,\ldots,\bq_{n-1},\bp_{n},\ldots,\bp_m]]$.
\bigskip

\noindent
Putting this all together will prove Lemma~\ref{lemma.TP.QPfinal.new2}
and hence Theorem~\ref{thm.prodmat.TP.bis.gen.new2}.

As preparation, we need some simple lemmas
concerning operations that preserve total positivity.
All of these hold in an arbitrary partially ordered commutative ring,
and apply to rectangular (i.e.\ not necessarily square)
matrices of suitably conformable dimensions.

\begin{lemma}[Increasing the upper-left entry]
   \label{lemma.TP.0}
Let $A$ be a (finite or infinite) matrix
with entries in a partially ordered commutative ring $R$;
let $c$ be a nonnegative element of $R$;
and let $A'$ be the matrix obtained from $A$
by adding $c$ to the upper-left entry
and leaving all other entries unchanged.
If $A$ is totally positive of order~$r$, then so is $A'$.
\end{lemma}

\proof
Consider a minor on rows $I$ and columns $J$.
If $I$ and $J$ do not both contain~1,
then obviously $\det (A')_{IJ} = \det A_{IJ}$.
If $I$ and $J$ both contain~1,
then $\det (A')_{IJ} = \det A_{IJ} + c \det A_{I\setminus 1, J\setminus 1}$.
\qed

\begin{lemma}[Matrices that coincide except in one row]
   \label{lemma.TP.1}
\hfill\break
\vspace*{-2mm}

\noindent
Let $\displaystyle A = \begin{bmatrix} M \\ \balpha \\ N \end{bmatrix}$
and $\displaystyle B = \begin{bmatrix} M \\ \bbeta \\ N \end{bmatrix}$
be matrices that coincide except in one row;
let $a,b \ge 0$;
and let
$\displaystyle C = \begin{bmatrix} M \\ a\balpha + b\bbeta \\ N \end{bmatrix}$.
If $A$ and $B$ are totally positive of order~$r$,
then so is $C$.
\end{lemma}

\proof
This is an immediate consequence of the row linearity of determinants.
\qed

\noindent
Of course, an analogous result holds for matrices that coincide
except in one column.

\begin{lemma}[Block-diagonal matrices]
   \label{lemma.TP.2}
If the matrices $A$ and $B$ are totally positive of order~$r$,
then so is the matrix
$\displaystyle C = \left[
                      \begin{array}{c|c}
                            A   &    \\
                            \hline
                                &   B
                       \end{array}
                    \right]$.
\end{lemma}

\proof
This is trivial: every nonzero minor of $C$
is a product of minors of $A$ and $B$.
\qed

A slightly less trivial result \cite[p.~398]{Karlin_68}
concerns matrices that are not quite block-diagonal,
being shifted so as to overlap in one row:

\begin{lemma}[Almost-block-diagonal matrices]
   \label{lemma.TP.3}
Let $A = \displaystyle \begin{bmatrix} A' \\ \balpha \end{bmatrix}$
and $B = \displaystyle \begin{bmatrix} \bbeta \\ B' \end{bmatrix}$
be two matrices,
where the row vector $\balpha$ is the last row of $A$
and the row vector $\bbeta$ is the first row of $B$.
Now form the matrix
\be
   M
   \;=\;
   \left[
      \begin{array}{c|c}
         A'  & \bzero    \\
         \hline
         \balpha & \bbeta \\
         \hline
         \bzero & B'
      \end{array}
   \right]
   \;,
\ee
in which the blocks $A$ and $B$ overlap in one row.
If $A$ and $B$ are totally positive of order~$r$, then so is $M$.
\end{lemma}

\proof
By Lemma~\ref{lemma.TP.2}, the block-diagonal matrices
\be
   M_1
   \;\eqdef\;
   \left[
       \begin{array}{c|c}
              \\[-3mm]
           A   &    \\[2mm]
           \hline
               &   B'
       \end{array}
   \right]
   \;=\;
   \left[
      \begin{array}{c|c}
         A'  &   \\
         \hline
         \balpha & \bzero  \\
         \hline
             & B'
      \end{array}
   \right]
   \quad\hbox{and}\quad
   M_2
   \;\eqdef\;
   \left[
       \begin{array}{c|c}
           A'   &    \\
           \hline
              \\[-2mm]
               &   B  \\[2mm]
       \end{array}
   \right]
   \;=\;
   \left[
      \begin{array}{c|c}
         A'  &   \\
         \hline
         \bzero   & \bbeta \\
         \hline
             & B'
      \end{array}
   \right]
\ee
are totally positive of order~$r$.
But then, by Lemma~\ref{lemma.TP.1}, so is $M$.
\qed

\noindent
Of course, an analogous result holds for matrices that overlap in one column.

\bigskip

We now begin the proof of Theorem~\ref{thm.prodmat.TP.bis.gen.new2}.

\begin{lemma}[Total positivity of $Q$]
   \label{lemma.TP.Q.new2}
The matrix $Q$ defined by \reff{def.prodmat.Q.new2}/\reff{def.Qmat2}
is totally positive, coefficientwise in the indeterminates
$\bfa,\bfb,\bfc,\bfd,\bfe,\bfff,\bfg$.

In particular, for every integer $m \ge 0$,
the matrix
	$[[\bq_0,\ldots,\bq_m]]$ 
is coefficientwise totally positive.
\end{lemma}

\proof
{}From \reff{eq.thm.prodmat.TP.bis.gen.new2}/\reff{def.prodmat.Q.new2}
we have
\be
   Q \;=\; L_1 (U L_2\, +\, D_1)  \;.
\ee
By Lemma~\ref{lemma.bidiagonal},
the matrices $L_1$, $L_2$, $U$ and $D_1$ are all totally positive.
Therefore, $U L_2$ is totally positive and tridiagonal,
while $D_1$ is totally positive and diagonal.
The tridiagonal comparison theorem (Proposition~\ref{prop.comparison})
then implies that $U L_2 + D_1$ is totally positive.
The result then follows.
\qed

\medskip

%

\begin{lemma}
   \label{lemma.TP.QR.new2}
For each integer $n \ge 0$,
the matrix
	$[[\bq_0,\ldots,\bq_{n-1},\bll_n]]$ 
is coefficientwise totally positive.
\end{lemma}

\proof
The case when $n=0$ is trivial.
So we assume that $n\geq 1$.

Let $\widetilde{\bq}_{n-1} = (\widetilde{q}_{n-1,k})_{k\geq 0}$
be the row vector  with entries
\begin{subeqnarray}
	\widetilde{q}_{n-1,n}  & = &   a_{n-1} c_{n} e_{n}   \\[1mm]
   \widetilde{q}_{n-1,n-1}    & = &
	a_{n-1} d_{n-1} e_{n-1} + b_{n-1} c_{n-1} e_{n-1} +   a_{n-1} g_{n-1} \\[1mm]
   \widetilde{q}_{n-1,n-2}  & = &
	a_{n-1} d_{n-1} f_{n-1} + b_{n-1} c_{n-1} f_{n-1} +  b_{n-1} d_{n-2} e_{n-2} + b_{n-1} g_{n-2}  \qquad \\[1mm]
   \widetilde{q}_{n-1,n-3}  & = &   b_{n-1} d_{n-2} f_{n-2}   \\[1mm]
   \widetilde{q}_{n-1,k}    & = &   0 \qquad\textrm{if $k < n-3$ or $k > n$}
\end{subeqnarray}
Thus, $\widetilde{\bq}_{n-1}$ is identical to $\bq_{n-1}$ except
that the entry $\widetilde{q}_{n-1,n-1}$ does not contain 
the summand $a_{n-1} c_n f_n$.
Also, let $\widetilde{\bll}_{n} = (\widetilde{l}_{n,k})_{k \ge 0}$
be the row vector
\be
\widetilde{\bll}_{n} \;=\; 
\begin{bmatrix}
	\bzero_{1\times n} & e_{n} & \bzero_{1\times \infty}
\end{bmatrix}
   \;.
\ee
Then
\begin{subeqnarray}
   q_{n-1,n-1}  & = & \widetilde{q}_{n-1,n-1} \,+\,\frac{f_n}{e_{n}} \,   \widetilde{q}_{n-1,n}
	\\[2mm]
   l_{n,n-1}  & = & \widetilde{l}_{n,n-1} \,+\, \frac{f_n}{e_{n}} \,  \widetilde{l}_{n,n}
        \label{eq.bidiag.rightmultiply}
\end{subeqnarray}
(actually $\widetilde{l}_{n,n-1} = 0$ but it is convenient to write it anyway).

By Lemma~\ref{lemma.TP.Q.new2}, the matrix
\be
M_1  \;=\;
[[ \bq_0, \ldots, \bq_{n-2}, \bq_{n-1}]]
\ee
is coefficientwise totally positive.
Therefore the matrix
\be
M_2  \;=\;
[[ \bq_0, \ldots, \bq_{n-2}, \widetilde{\bq}_{n-1}]]
   \;,
\ee
which is obtained from $M_1$ by setting $f_{n} = 0$,
is also coefficientwise totally positive.
(Setting $f_n = 0$ also affects rows $n$ and $n+1$ of $Q$,
 but these are not contained in $M_1$.)
Then
\be
M_3 \;=\;
[[\bq_0, \ldots, \bq_{n-2}, \widetilde{\bq}_{n-1}, \widetilde{\bll}_{n}]]
\ee
is obtained from $M_2$ by adjoining a row that has $e_n$ in the $n$th column
and 0 elsewhere;
this is totally positive because $M_2$ is zero beyond the $n$th column
(this reasoning is a very special case of Lemma~\ref{lemma.TP.3}).
Finally, the matrix 
\be
M \;=\;
[[\bq_0, \ldots, \bq_{n-1}, \bll_{n}]]
\ee
can be obtained by right-multiplying $M_3$
by the lower-bidiagonal matrix that has $1$ on the diagonal,
$f_n/e_{n}$ in position $(n,n-1)$, and 0 elsewhere
(this adds $f_n/e_{n}$ times the $n$th column to the $(n-1)$st column,
which is the content of \reff{eq.bidiag.rightmultiply}).
This completes the proof of Lemma~\ref{lemma.TP.QR.new2}.
\qed

{\bf Remark.}  To justify the use of the quantity $f_n/e_n$,
it suffices to work in a ring of Laurent polynomials
with indeterminates $\bfe = (e_n)_{n \ge 0}$
and $\bfe^{-1} = (e_n^{-1})_{n \ge 0}$,
and then to consider coefficientwise total positivity in this ring.
At the end, all quantities will belong to the subring
consisting of ordinary polynomials with nonnegative exponents.
\myendremark

\begin{lemma}
   \label{lemma.TP.QPR.if.new2}
Fix integers $0 \le n \le m$.
If the matrix 
$
	[[\bp_{n+1},\ldots, \bp_m]]
$
is coefficientwise totally positive,
then so is the matrix
$
M_0 \;\eqdef\;
	[[\bq_0, \ldots, \bq_{n-1}, \bll_{n},\bp_{n+1},\ldots, \bp_m]]
   \;.
$
\end{lemma}

\proof
The case $n=m$ is Lemma~\ref{lemma.TP.QR.new2};  so assume that $n < m$.

Let $\bt_{n+1}$ be obtained from $\bp_{n+1}$ by specializing $d_{n}$ to zero;
then
\be
   \bp_{n+1}  \;=\;  \bt_{n+1} \,+\, b_{n+1} d_n \bll_n
   \;.
   \label{eq.p.t.difference}
\ee
Note that this specialization does not affect $\bp_\ell$ for $\ell > n+1$.
By hypothesis 
\be
M_1
\;\eqdef\;
\begin{bmatrix}
        \bp_{n+1}\\
        \vdots\\
        \bp_{m}\\
\end{bmatrix}
\ee
is coefficientwise totally positive;
this implies, by specialization, that
\be
M_2
\;\eqdef\;
\begin{bmatrix}
        \bt_{n+1}\\
	\bp_{n+2}\\
        \vdots\\
        \bp_{m}\\
\end{bmatrix}
\ee
is coefficientwise totally positive.
Also let
\be
M_3 
\;\eqdef\;
\begin{bmatrix}
        \bq_0\\
        \vdots\\
        \bq_{n-1}\\
        \bll_n
\end{bmatrix}
\;.
\ee
(Here for clarity we have avoided the notation $[[\;\cdots\;]]$.)

Now observe that the matrix $S \eqdef 
\left[
\begin{array}{c}
	M_3\\
	\hline
	M_2
\end{array}
\right]$
consists of two blocks overlapping in a single column
(namely, column $n$ when the columns are numbered starting at 0):
\be
   S
   \;=\;
   \left[
\begin{array}{c}
        M_3\\
        \hline
        M_2
\end{array}
\right]\;=\;
   \left[
   \begin{array}{c}
        \bq_{0}\\
        \vdots\\
        \bq_{n-1}\\
        \bll_n\\
	   \hline
        \bt_{n+1}\\
        \bp_{n+2}\\
        \vdots\\
        \bp_m
   \end{array}
   \right]
   \;=\;
   \left[
       \begin{array}{c|c|c|c}
	       *   &   *  &  * & \bzero_{n\times \infty} \\[1mm]
	       \bzero_{1\times n-1}   & f_n & e_n & \bzero_{1\times \infty} \\[1mm]
           \hline
           0   &  0      & *   &  *                   \\[1mm]
       \end{array}
   \right]
\ee
where the asterisks stand for blocks of unspecified entries
(which may be zero or nonzero).
By Lemma~\ref{lemma.TP.QR.new2}, 
the matrix $M_3$ 
is totally positive;
and we have just shown that the matrix
$M_2$
is totally positive.
So the transpose of Lemma~\ref{lemma.TP.3} 
implies that the matrix $S$ is totally positive.

On the other hand, 
the matrix 
$M_0$
can be obtained from $S$ by left-multiplying it
by the lower-bidiagonal matrix that has
1 on the diagonal, $b_{n+1}d_n$ 
in position $(n+1,n)$ and zeroes elsewhere
(this adds $b_{n+1}d_n$ times row~$n$ to row~$n+1$,
which is the content of \reff{eq.p.t.difference}).
This proves that the desired matrix $M_0$ is totally positive.
\qed

\begin{lemma}[Induction step]
   \label{lemma.TP.QP.inductive.new2}
\hfill\break
Fix integers $0 \le n \le m$.
If the matrix 
$[[\bq_0,\ldots,\bq_{n},\bp_{n+1},\ldots,\bp_m]]$
is coefficientwise totally positive,
then so is 
$[[\bq_0,\ldots,\bq_{n-1},\bp_{n},\ldots,\bp_m]]$.
\end{lemma}

\proof
Notice that $\bp_n = \bq_n + h_n \bll_n$.
The matrix
$
   [[\bq_0,\ldots,\bq_{n-1},\bq_{n},\bp_{n+1},\ldots,\bp_m]]
$
is totally positive by hypothesis;
and the matrix
 $[[\bq_0,\ldots,\bq_{n-1},\bll_n,\bp_{n+1},\ldots,\bp_m]]
$
is totally positive by Lemma~\ref{lemma.TP.QPR.if.new2}.
So the conclusion follows by Lemma~\ref{lemma.TP.1}
with $a=1$ and $b = h_n$.
\qed

\begin{lemma}
   \label{lemma.TP.QPfinal.new2}
\hfill\break
For every pair of integers $0 \le n \le m+1$,
the matrix 
$[[\bq_0,\ldots,\bq_{n-1},\bp_n,\ldots,\bp_m]]$ 
is coefficientwise totally positive.
\end{lemma}

\proof
Start from Lemma~\ref{lemma.TP.Q.new2},
and apply Lemma~\ref{lemma.TP.QP.inductive.new2} for $n=m,m-1,\ldots,0$.
\qed

This completes the proof of Theorem~\ref{thm.prodmat.TP.bis.gen.new2}.

\bigskip

We conclude by posing the following open problem:

\begin{problem}
Find a combinatorial interpretation for the output matrix $A = \scro(P)$
generated by the production matrix
\reff{eq.thm.prodmat.TP.bis.gen.new2}/\reff{eq.thm.prodmat.TP.bis.gen.new2.2},
or by interesting specializations thereof.
\end{problem}

\section{The connection with multiple orthogonal polynomials}  \label{sec.MOP}

In this section we would like to explain how we were led to guess
the quadridiagonal production matrix \reff{eq.prop.prodmat.univariate.NEW.2},
by virtue of an unexpected connection with multiple orthogonal polynomials.

Let us first show that, for each $\alpha \ge -1$ and $x \ge 0$,
the sequence $(\scrlna(x))_{n \ge 0}$ is a Stieltjes moment sequence,
by explicitly exhibiting its Stieltjes moment representation.
Indeed, we have the integral representation
(compare \cite[p.~103, eq.~(5.4.1)]{Szego_75})
\be
   \scrlna(x)
   \;=\;
   e^{-x} x^{-\alpha/2}  \int\limits_0^\infty \! u^{n+\alpha/2} \, e^{-u} \,
                             I_\alpha(2 \sqrt{xu}) \: du
 \label{eq.scrlna.stieltjes}
\ee
valid for $\alpha > -1$, where $I_\alpha$ is a modified Bessel function
of the first kind \cite[p.~77]{Watson_44}
\be
   I_\alpha(z)
   \;=\;
   \sum_{k=0}^\infty {(z/2)^{\alpha+2k} \over k! \, \Gamma(\alpha+k+1)}
   \;.
 \label{def.BesselI}
\ee
For $n=0$, the formula \reff{eq.scrlna.stieltjes} is easily verified
by inserting \reff{def.BesselI} and integrating term-by-term.
The general formula then follows by computing the
exponential generating function of the right-hand side
of \reff{eq.scrlna.stieltjes},
i.e.\ multiplying \reff{eq.scrlna.stieltjes} by $t^n/n!$
and summing over $n \ge 0$;
the result is $(1-t)^{-(\alpha+1)} e^{xt/(1-t)}$,
which is indeed [by \reff{eq.genfn.1}]
the exponential generating function of $\scrlna(x) = n! \, L_n^{(\alpha)}(x)$.
And since $I_\alpha$ is nonnegative on $[0,\infty)$,
the integral representation \reff{eq.scrlna.stieltjes}
is a Stieltjes moment representation, with representing measure
\be
   d\mu_{\alpha,x}(u)
   \;=\;
   e^{-x} x^{-\alpha/2} \, u^{\alpha/2} \, e^{-u} \, I_\alpha(2 \sqrt{xu})
   \: du
   \quad\hbox{ on } [0,\infty)
 \label{def.mu.alpha.x}
\ee
for $\alpha > -1$
(with the corresponding limiting measure when $\alpha = -1$).\footnote{
   These facts were observed many years ago by Karlin
   \cite[p.~440]{Karlin_68} \cite[p.~62]{Karlin_68b}.
}

We now invoke a recently discovered \cite{Sokal_multiple_OP}
connection between multiple orthogonal polynomials
and branched continued fractions.
Recall that multiple orthogonal polynomials \cite[Chapter~23]{Ismail_05}
are a generalization of conventional orthogonal polynomials
\cite{Szego_75,Chihara_78,Ismail_05}
in which the polynomials satisfy orthogonality relations
with respect to several measures $\mu_1,\ldots,\mu_r$ rather than just one.
In particular, the multiple orthogonal polynomials of type~II,
denoted $P_{\bfn}(x)$
and indexed by a multi-index $\bfn = (n_1,\ldots,n_r) \in \N^r$,
are monic polynomials
of degree $|\bfn| = \sum_{i=1}^r n_i$
that satisfy the orthogonality relations
\be
   \int \! x^k \, P_{\bfn}(x) \: d\mu_i(x)
   \;=\;
   0
   \qquad\hbox{for $1 \le i \le r$ and $0 \le k < n_i$}
   \;.
 \label{eq.MOP}
\ee
(Here we restrict attention to the case in which
the system of measures $\mu_1,\ldots,\mu_r$ is {\em perfect}\/:
namely, there exists, for each multi-index $\bfn$,
a {\em unique}\/ monic polynomial $P_{\bfn}(x)$
of degree $|\bfn|$ satisfying \reff{eq.MOP}.
Several general sufficient conditions for a system to be perfect are known
(Angelesco systems, AT systems, Nikishin systems, \ldots):
see \cite[Chapter~23]{Ismail_05} \cite{VanAssche_21}.)

Now recall
\cite[Theorems~50.1 and 51.1]{Wall_48}
\cite[Chapitre~V]{Viennot_83}
\cite[Theorem~2.3 and Corollary~2.5]{Corteel_16}
\cite[Section~5.2.1]{Zeng_21}
that the conventional orthogonal polynomials
associated to a measure~$\mu$ (which we here normalize to total mass~1)
satisfy a three-term recurrence relation of the form
\be
   P_{n+1}(x)
   \;=\;
   (x - \gamma_n) \, P_n(x) \:-\: \beta_n \, P_{n-1}(x)
   \;,
\ee
where the coefficients $\gamma_n$ and $\beta_n$ occurring
in this recurrence relation are identical to those occurring
in the classical J-fraction for the moments of $\mu$,
i.e.
\be
   \sum_{n=0}^\infty a_n \, t^n
   \;=\;
   \cfrac{1}{1 - \gamma_0 t - \cfrac{\beta_1 t^2}{1 - \gamma_1 t - \cfrac{\beta_2 t^2}{1 - \gamma_2 t - \cfrac{\beta_3 t^2}{1- \gamma_3 t - \cdots}}}}
 \label{eq.Jtype.cfrac}
\ee
where $a_n = \int \! x^n \: d\mu(x)$.
Now, Flajolet \cite{Flajolet_80} showed that the Taylor coefficient $a_n$
in \reff{eq.Jtype.cfrac}
is the generating polynomial for Motzkin paths from $(0,0)$ to $(n,0)$
in which each rise gets weight 1,
each level step at height $i$ gets weight $\gamma_i$,
and each fall from height $i$ gets weight $\beta_i$.
In other words, $a_n$ is the $(n,0)$ matrix element
of the output matrix $A = \scro(\Pi)$,
where $\Pi$ is the tridiagonal production matrix
\be
   \Pi
   \;=\;
   \begin{bmatrix}
      \gamma_0   & 1           &            &        &         \\
      \beta_1    & \gamma_1    & 1          &        &         \\
                 & \beta_2     & \gamma_2   & 1      &         \\
                 &             & \ddots     & \ddots & \ddots
   \end{bmatrix}
   \;.
 \label{def.Jnl.production.bis}
\ee

It turns out \cite{Sokal_multiple_OP}
that a generalization of this connection
holds for multiple orthogonal polynomials.
Fix measures $\mu_1,\ldots,\mu_r$,
and let $(P_\bfn(x))_{\bfn \in \N^r}$
be the corresponding multiple orthogonal polynomials of type~II.
Among these, let us consider the multi-indices $\bfn = (n_1,\ldots,n_r)$
lying on the so-called \textbfit{stepline}:
this is the near-diagonal sequence starting at $(0,0,\ldots,0)$
and following the path
$
   (n,n,\ldots,n) \to
   (n+1,n,\ldots,n) \to
   (n+1,n+1,\ldots,n) \to \ldots \to
   (n+1,n+1,\ldots,n+1) \to \ldots
   \;.
$
That is, we define a singly-indexed sequence $(\Ptilde_n(x))_{n \ge 0}$
by
\be
   \Ptilde_n(x)  \;=\;  P_{(n_1,\ldots,n_r)}(x)
   \quad\hbox{where }
   n_i \:=\: \Bigl\lfloor {n+r-i \over r} \Bigr\rfloor
   \;\;\hbox{for } 1 \le i \le r
   \;.
 \label{def.Ptilde.stepline}
\ee
It is well known \cite[Theorem~23.1.7]{Ismail_05}
that the stepline sequence $(\Ptilde_n(x))_{n \ge 0}$
satisfies an $(r+2)$-term recurrence of the form
\be
   \Ptilde_{n+1}(x)
   \;=\;
   (x - \pi_{nn}) \, \Ptilde_n(x) \:-\:
     \sum_{i=1}^{r} \pi_{n,n-i} \, \Ptilde_{n-i}(x)
   \;,
 \label{eq.prop.Ptildenx.recurrence.1}
\ee
or equivalently
\be
   x \Ptilde_n(x)   \;=\;  \sum_{k=n-r}^{n+1} \pi_{nk} \, \Ptilde_k(x)
 \label{eq.prop.Ptildenx.recurrence.2}
\ee
where $\pi_{n,n+1} = 1$.
Here the production matrix $\Pi = (\pi_{nk})_{n,k \ge 0}$
is an $(r,1)$-banded unit-lower-Hessenberg matrix,
i.e.\ $\pi_{nk} = 0$ for $k < n-r$ or $k > n+1$.
Furthermore, it turns out \cite{Sokal_multiple_OP}
--- and this is the key fact ---
that the zeroth column of the output matrix $A = \scro(\Pi)$
is precisely the moment sequence of the measure $\mu_1$.
Therefore, if one knows the recurrence relation
\reff{eq.prop.Ptildenx.recurrence.1}/\reff{eq.prop.Ptildenx.recurrence.2}
for the multiple orthogonal polynomials along the stepline,
one also knows a production matrix for
the moment sequence of the measure $\mu_1$.

\medskip

{\bf Remark.} The subsequent columns of the output matrix $A = \scro(\Pi)$
also have interpretations \cite{Sokal_multiple_OP}
as linear combinations of the moment sequences of $\mu_1,\ldots,\mu_r$.
For instance, for $0 \le i \le r-1$,
the $i$th column of the output matrix
is that particular linear combination of the moment sequences
of $\mu_1,\mu_2,\ldots,\mu_{i+1}$
that annihilates the first $i$ entries
and makes the next entry equal to 1.
But we will not need this fact here.
\myendremark

\medskip

Return now to our Laguerre polynomials
and their Stieltjes moment representation
with the measure $\mu_{\alpha,x}$ defined in \reff{def.mu.alpha.x}.
We then benefit from exceedingly good fortune:
some years ago, Coussement and {Van Assche} \cite{Coussement_03}
studied the multiple orthogonal polynomials with $r=2$
associated to the pair of measures
$(\mu_1,\mu_2) = (\mu_{\alpha,x},\mu_{\alpha+1,x})$.
(They actually used a slightly different normalization,
 so that their moments are a rescaled version of the
 reversed monic unsigned Laguerre polynomials $\scrlbarna(x)$:
 compare \cite[Lemma~1]{Coussement_03} to our~\reff{def.scrlbarna}.)
And they computed explicitly, among other things,
the four-term recurrence relation for the multiple orthogonal polynomials
of type~II along the stepline \cite[Theorem~9]{Coussement_03}.
After translating their normalization to ours,
this four-term recurrence corresponds precisely to the
quadridiagonal production matrix \reff{eq.prop.prodmat.univariate.NEW.2}
for the monic unsigned Laguerre polynomials $\scrlna(x)$.
This is how we first discovered these production matrices.

{}From here it was a small step to discover the production matrix
\reff{eq.prop.multivariate_second.rowgen.prodmat.flat.BIS0}
for our multivariate Laguerre polynomials.
Indeed, once we had exploited the Foata--Strehl \cite{Foata_84}
combinatorial model to define the multivariate Laguerre polynomials
\reff{def.multivariate.rowgen.scrl.second},
it was not difficult to guess
--- helped by a bit of computer experimentation ---
how the production matrix \reff{eq.prop.prodmat.univariate.NEW.2},
with its single variable $x$, should be refined to introduce
the further variables $\yp,\yv,\yda,\ydd,\yfp$.

Let us also observe that there is another way that we could have
discovered the production matrix
\reff{eq.prop.prodmat.univariate.NEW.2}
for the univariate Laguerre polynomials.
As observed preceding Proposition~\ref{prop.L.SR},
we already knew a tridiagonal production matrix for
the Laguerre polynomials at $x=0$,
which are $\scrlna(0) = (1+\alpha)^{\overline{n}}$.
It would then have been natural to try perturbing this tridiagonal matrix
by terms linear in $x$;  and if one also realized that one ought to consider
introducing a second subdiagonal, one could try the Ansatz
\begin{subeqnarray}
   p_{n,n+1}  & = &   1   \\[1mm]
   p_{n,n}    & = &   (2n+1+\alpha) \,+\, A_n x   \\[1mm]
   p_{n,n-1}  & = &   n(n+\alpha) \,+\, B_n x \\[1mm]
   p_{n,n-2}  & = &   C_n x  \\[1mm]
   p_{n,k}    & = &   0 \qquad\textrm{if $k < n-2$ or $k > n+1$}
\end{subeqnarray}
A little computer work
--- setting the zeroth column of the output matrix equal to $\scrlna(x)$ ---
would then lead quickly to the values of $A_n,B_n,C_n$ for small $n$,
from which one could easily guess that $A_n = 1$, $B_n = 2n$ and $C_n = n(n-1)$.
This is not, in fact, the way we first found the production matrix
\reff{eq.prop.prodmat.univariate.NEW.2},
but it easily could have been.

\section*{Acknowledgments}


One of us (A.D.S.)\ wishes to thank the organizers of the
15th International Symposium on Orthogonal Polynomials,
Special Functions and Applications
(Hagenberg, Austria, 22--26 July 2019)
for inviting him to give a talk there;
this allowed him to meet Walter Van Assche and to discover
the unexpected connection \cite{Sokal_multiple_OP}
between branched continued fractions and multiple orthogonal polynomials,
which was the starting point for this work.

We wish to thank Xi Chen, Sylvie Corteel and Walter Van Assche
for helpful discussions.

This research was supported in part by
the U.K.~Engineering and Physical Sciences Research Council grant EP/N025636/1
and by the research fellowship DY 133/1-1
from the Deutsche Forschungsgemeinschaft.
The first author is currently supported by
DIMERS project ANR-18-CE40-0033 funded by
Agence Nationale de la Recherche (ANR, France).

\appendix

\section{Generalized and modified $\bm{m}$-Stieltjes--Rogers polynomials}   \label{app.modified}

In this Appendix we introduce, for each integer $m \ge 1$,
an infinite sequence of triangular arrays, indexed by an integer $j \ge 0$,
whose matrix elements we call
the {\em generalized $m$-Stieltjes--Rogers polynomials of type $j$}\/.
We examine in particular the entries in these matrices' zeroth columns,
which we call the
{\em modified $m$-Stieltjes--Rogers polynomials of type $j$}\/.
When $j=0$ these polynomials
reduce to the (generalized) $m$-Stieltjes--Rogers polynomials
introduced in \cite[sections~5, 7 and 9]{latpath_SRTR}.
When $m=1$ they reduce to the generalized
classical Stieltjes--Rogers polynomials
of the first ($j=0$) and second ($j=1$) kinds.
The discussion here expands and supersedes the treatment
in \cite{latpath_SRTR}.

We begin (Section~\ref{subsec.SR})
by reviewing the the theory of classical Stieltjes--Rogers polynomials,
and then (Section~\ref{subsec.mSR})
the basic ideas from \cite{latpath_SRTR}
concerning the $m$-Stieltjes--Rogers polynomials (of type 0).
Then (Section~\ref{subsec.mSR.gen}) we introduce the
generalized $m$-Stieltjes--Rogers polynomials of type $j$.
Finally (Section~\ref{subsec.mSR.univariate}),
we apply this theory to the univariate Laguerre production matrix
\reff{eq.prop.prodmat.univariate.NEW.2}.

\subsection{Classical Stieltjes--Rogers polynomials}  \label{subsec.SR}

A \textbfit{Dyck path} is a path in the upper half-plane $\Z \times \N$,
starting and ending on the horizontal axis,
using steps $(1,1)$ [``rise'' or ``up step'']
and $(1,-1)$ [``fall'' or ``down step''].
More generally, a \textbfit{Dyck path at level $k$}
is a path in $\Z \times \N_{\ge k}$,
starting and ending at height $k$,
using steps $(1,1)$ and $(1,-1)$.
Clearly a Dyck path must be of even length;
we denote by $\scrd_{2n}$ the set of Dyck paths from $(0,0)$ to $(2n,0)$.

Now let $\balpha = (\alpha_i)_{i \ge 1}$ be an infinite set of indeterminates,
and let $S_n(\balpha)$ be the generating polynomial
for Dyck paths from $(0,0)$ to $(2n,0)$ in which each rise gets weight~1
and each fall from height~$i$ gets weight $\alpha_i$.
Clearly $S_n(\balpha)$ is a homogeneous polynomial
of degree~$n$ with nonnegative integer coefficients;
following Flajolet \cite{Flajolet_80},
we call it the \textbfit{Stieltjes--Rogers polynomial} of order~$n$.

Let $f_0(t) = \sum_{n=0}^\infty S_n(\balpha) \, t^n$
be the ordinary generating function for Dyck paths with these weights
(considered as a formal power series in $t$);
and more generally, let $f_k(t)$ be the ordinary generating function
for Dyck paths at level $k$ with these same weights.
(Obviously $f_k$ is just $f_0$ with each $\alpha_i$ replaced by $\alpha_{i+k}$;
 but we shall not explicitly use this fact.)
Then a straightforward ``renewal'' argument \cite{Flajolet_80}
gives the functional equation
\be
   f_k(t)  \;=\;  1 \:+\: \alpha_{k+1} t \, f_k(t) \, f_{k+1}(t)
 \label{eq.SRfk.1}
\ee
or equivalently
\be
   f_k(t)  \;=\;  {1 \over 1 \:-\: \alpha_{k+1} t \, f_{k+1}(t)}
   \;.
 \label{eq.SRfk.2}
\ee
Iterating \reff{eq.SRfk.2}, we see immediately that $f_k$
is given by the continued fraction
\be
   f_k(t)
   \;=\;
   \cfrac{1}{1 - \cfrac{\alpha_{k+1} t}{1 - \cfrac{\alpha_{k+2} t}{1- \cfrac{\alpha_{k+3} t}{1- \cdots}}}}
 \label{eq.fk.Sfrac}
\ee
and in particular that $f_0$ is given by
\be
   f_0(t)
   \;=\;
   \cfrac{1}{1 - \cfrac{\alpha_{1} t}{1 - \cfrac{\alpha_{2} t}{1- \cfrac{\alpha_{3} t}{1- \cdots}}}}
   \;\,.
 \label{eq.f0.Sfrac}
\ee
The right-hand sides of \reff{eq.fk.Sfrac}/\reff{eq.f0.Sfrac}
are called (classical) \textbfit{Stieltjes-type continued fractions},
or (classical) \textbfit{S-fractions} for short.
This combinatorial interpretation of S-fractions
in terms of weighted Dyck paths is due to Flajolet \cite{Flajolet_80}.

We now generalize the Stieltjes--Rogers polynomials to a triangular array,
as follows:
We use the term \textbfit{partial Dyck path}
to denote a path in the upper half-plane $\Z \times \N$,
using steps $(1,1)$ and $(1,-1)$,
that starts on the horizontal axis
but is allowed to end anywhere in the upper half-plane.
For $n,\ell \ge 0$, we define the
\textbfit{generalized Stieltjes--Rogers polynomial of the first kind}
$S_{n,\ell}(\balpha)$
to be the generating polynomial for partial Dyck paths
starting at $(0,0)$ and ending at $(2n,2\ell)$,
in which each rise gets weight 1
and each fall from height~$i$ gets weight $\alpha_i$.
Obviously $S_{n,\ell}$ is nonvanishing only for $0 \le \ell \le n$,
so we have an infinite lower-triangular array
$\sfS = (S_{n,\ell}(\balpha))_{n,\ell \ge 0}$
in which the zeroth column displays
the ordinary Stieltjes--Rogers polynomials $S_{n,0} = S_n$.
In particular we have $S_{n,n} = 1$
and $S_{n,n-1} = \sum_{i=1}^{2n-1} \alpha_i$.

Analogously, we define
the \textbfit{generalized Stieltjes--Rogers polynomial of the second kind}
$S'_{n,\ell}(\balpha)$
to be the generating polynomial for partial Dyck paths
starting at $(0,0)$ and ending at $(2n+1,2\ell+1)$,
in which again each rise gets weight 1
and each fall from height~$i$ gets weight $\alpha_i$.
Since $S'_{n,\ell}$ is nonvanishing only for $0 \le \ell \le n$,
we obtain a second infinite lower-triangular array
$\sfS' = (S'_{n,\ell}(\balpha))_{n,\ell \ge 0}$.
In particular we have $S'_{n,n} = 1$
and $S'_{n,n-1} = \sum_{i=1}^{2n} \alpha_i$.

The polynomials $S_{n,\ell}(\balpha)$ and $S'_{n,\ell}(\balpha)$
manifestly satisfy the joint recurrence
\begin{subeqnarray}
   S'_{n,\ell}   & = &  S_{n,\ell} \:+\: \alpha_{2\ell+2} \, S_{n,\ell+1}
       \\[2mm]
   S_{n+1,\ell}  & = &  S'_{n,\ell-1} \:+\: \alpha_{2\ell+1} \, S'_{n,\ell}
 \slabel{eq.SnlSprimenl.recurrence.b}
 \label{eq.SnlSprimenl.recurrence}
\end{subeqnarray}
for $n,\ell \ge 0$,
with the initial condition $S_{0,\ell} = \delta_{\ell 0}$
(where of course we also set ${S'_{n,-1} = 0}$).
It follows that the $S_{n,\ell}$ satisfy the recurrence
\be
   S_{n+1,\ell}
   \;=\;
   S_{n,\ell-1}
      \:+\:  (\alpha_{2\ell} + \alpha_{2\ell+1}) \, S_{n,\ell}
      \:+\:  \alpha_{2\ell+1} \alpha_{2\ell+2} \, S_{n,\ell+1}
 \label{eq.Snl.recursion}
\ee
with the initial condition $S_{0,\ell} = \delta_{\ell 0}$
(where we set $S_{n,-1} = 0$ and $\alpha_0 = 0$).
In other words, the unit-lower-triangular array $\sfS$
has the tridiagonal production matrix
\be
   P  \;=\;
   \begin{bmatrix}
      \alpha_1          & 1                   &     &   &      \\
      \alpha_1 \alpha_2 & \alpha_2 + \alpha_3 & 1   &   &      \\
                        & \alpha_3 \alpha_4   & \alpha_4 + \alpha_5 & 1   &  \\
                        &                     & \ddots & \ddots & \ddots
   \end{bmatrix}
   \;.
 \label{def.Snl.production}
\ee

Please observe now that the production matrix \reff{def.Snl.production}
can be factorized as a product of two bidiagonal matrices:
\be
\Scale[0.95]{
   \begin{bmatrix}
      \alpha_1          & 1                   &     &   &      \\
      \alpha_1 \alpha_2 & \alpha_2 + \alpha_3 & 1   &   &      \\
                        & \alpha_3 \alpha_4   & \alpha_4 + \alpha_5 & 1   &  \\
                        &                     & \ddots & \ddots & \ddots
   \end{bmatrix}
   \;=\;
   \begin{bmatrix}
      1                 &                     &     &   &      \\
      \alpha_2          & 1                   &     &   &      \\
                        & \alpha_4            & 1   &   &  \\
                        &                     & \ddots & \ddots &
   \end{bmatrix}
   \begin{bmatrix}
      \alpha_1          & 1                   &     &   &      \\
                        & \alpha_3            & 1   &   &      \\
                        &                     & \alpha_5 & 1   &  \\
                        &                     &        & \ddots & \ddots
   \end{bmatrix}
}
 \label{def.Snl.production.factorization}
\ee
By Lemma~\ref{lemma.bidiagonal} this shows \cite[Lemma~3.3]{Wang_16}
that the production matrix \reff{def.Snl.production}
is totally positive in the ring $Z[\balpha]$ equipped with the
coefficientwise order.
It then follows \cite{Sokal_totalpos}
from Theorems~\ref{thm.iteration.homo} and \ref{thm.iteration2bis}
that the unit-lower-triangular array $\sfS$ of
generalized Stieltjes--Rogers polynomials
is coefficientwise totally positive,
and that the sequence $(S_n(\balpha))_{n \ge 0}$
of Stieltjes--Rogers polynomials is coefficientwise Hankel-totally positive.

Similar considerations apply to the matrix $\sfS'$
of generalized Stieltjes--Rogers polynomials of the second kind
\cite[eq.~(7.12)]{latpath_SRTR};
we leave the details as an exercise for the reader.

\subsection[$m$-Stieltjes--Rogers polynomials (branched S-fractions)]{$\bm{m}$-Stieltjes--Rogers polynomials (branched S-fractions)}  \label{subsec.mSR}


In what follows we fix an integer $m \ge 1$.
We recall \cite{Aval_08,Cameron_16,Prodinger_16,latpath_SRTR}
that an \textbfit{$\bm{m}$-Dyck path}
is a path in the upper half-plane $\Z \times \N$,
starting and ending on the horizontal axis,
using steps $(1,1)$ [``rise'' or ``up step'']
and $(1,-m)$ [``$m$-fall'' or ``down step''].
More generally, an \textbfit{$\bm{m}$-Dyck path at level $\bm{k}$}
is a path in $\Z \times \N_{\ge k}$,
starting and ending at height $k$,
using steps $(1,1)$ and $(1,-m)$.
Since the number of up steps must equal $m$ times the number of down steps,
the length of an $m$-Dyck path must be a multiple of $m+1$.

Now let $\balpha = (\alpha_i)_{i \ge m}$ be an infinite set of indeterminates.
Then \cite{latpath_SRTR}
the \textbfit{$\bm{m}$-Stieltjes--Rogers polynomial} of order~$n$,
denoted $S^{(m)}_n(\balpha)$, is the generating polynomial
for $m$-Dyck paths from $(0,0)$ to $((m+1)n,0)$
in which each rise gets weight~1
and each $m$-fall from height~$i$ gets weight $\alpha_i$.
Clearly $S_n^{(m)}(\balpha)$ is a homogeneous polynomial
of degree~$n$ with nonnegative integer coefficients.

Let $f_0(t) = \sum_{n=0}^\infty S^{(m)}_n(\balpha) \, t^n$
be the ordinary generating function for $m$-Dyck paths with these weights;
and more generally, let $f_k(t)$ be the ordinary generating function
for $m$-Dyck paths at level $k$ with these same weights.
(Obviously $f_k$ is just $f_0$ with each $\alpha_i$ replaced by $\alpha_{i+k}$;
 but we shall not explicitly use this fact.)
Then straightforward combinatorial arguments \cite{latpath_SRTR}
lead to the functional equation
\be
   f_k(t)  \;=\;  1 \:+\: \alpha_{k+m} t \, f_k(t) \, f_{k+1}(t) \,\cdots\, f_{k+m}(t)
 \label{eq.mSRfk.1}
\ee
or equivalently
\be
   f_k(t)  \;=\;  {1 \over 1 \:-\: \alpha_{k+m} t \, f_{k+1}(t) \,\cdots\, f_{k+m}(t)}
   \;.
 \label{eq.mSRfk.2}
\ee
Iterating \reff{eq.mSRfk.2}, we see immediately that $f_k$
is given by the branched continued fraction
\begin{subeqnarray}
   f_k(t)
   & = &
   \cfrac{1}
         {1 \,-\, \alpha_{k+m} t
            \prod\limits_{i_1=1}^{m}
                 \cfrac{1}
            {1 \,-\, \alpha_{k+m+i_1} t
               \prod\limits_{i_2=1}^{m}
               \cfrac{1}
            {1 \,-\, \alpha_{k+m+i_1+i_2} t
               \prod\limits_{i_3=1}^{m}
               \cfrac{1}{1 - \cdots}
            }
           }
         }
%
      \slabel{eq.fk.mSfrac.a} \\[2mm]
   & = &
\Scale[0.6]{
   \cfrac{1}{1 - \cfrac{\alpha_{k+m} t}{
     \Biggl( 1 - \cfrac{\alpha_{k+m+1} t}{
        \Bigl( 1  - \cfrac{\alpha_{k+m+2} t}{(\cdots) \,\cdots\, (\cdots)} \Bigr)
        \,\cdots\,
        \Bigl( 1  - \cfrac{\alpha_{k+2m+1} t}{(\cdots) \,\cdots\, (\cdots)} \Bigr)
       }
     \Biggr)
     \,\cdots\,
     \Biggl( 1 - \cfrac{\alpha_{k+2m} t}{
        \Bigl( 1  - \cfrac{\alpha_{k+2m+1} t}{(\cdots) \,\cdots\, (\cdots)} \Bigr)
        \,\cdots\,
        \Bigl( 1  - \cfrac{\alpha_{k+3m} t}{(\cdots) \,\cdots\, (\cdots)} \Bigr)
       }
     \Biggr)
    }
   }
}
     \nonumber \\
 \slabel{eq.fk.mSfrac.b}
 \label{eq.fk.mSfrac}
\end{subeqnarray}
and in particular that $f_0$ is given by
the specialization of \reff{eq.fk.mSfrac} to $k=0$.
We shall call the right-hand side of \reff{eq.fk.mSfrac}
an \textbfit{$\bm{m}$-branched Stieltjes-type continued fraction},
or \textbfit{$\bm{m}$-S-fraction} for short.

\medskip

{\bf Remark.}
In truth, we hardly ever use the branched continued fraction
\reff{eq.fk.mSfrac};
instead, we work directly with the $m$-Dyck paths
and/or with the recurrence \reff{eq.mSRfk.1}/\reff{eq.mSRfk.2}
that their generating functions satisfy.
\myendremark

\medskip


\subsection[Generalized $m$-Stieltjes--Rogers polynomials of type $j$]{Generalized $\bm{m}$-Stieltjes--Rogers polynomials of type $\bm{j}$}
    \label{subsec.mSR.gen}

Fix again an integer $m \ge 1$.
A \textbfit{partial $\bm{m}$-Dyck path}
is a path in the upper half-plane $\Z \times \N$,
starting on the horizontal axis but ending anywhere in the upper half-plane,
using steps $(1,1)$ [``rise''] and $(1,-m)$ [``$m$-fall''].
A partial $m$-Dyck path starting at $(0,0)$
must stay always within the set
$V_m = \{ (x,y) \in \Z \times \N \colon\: x=y \bmod (m+1) \}$.
An \textbfit{$\bm{m}$-Dyck path} is simply a partial $m$-Dyck path
that ends on the horizontal axis.

Now let $\balpha = (\alpha_i)_{i \ge m}$ be an infinite set of indeterminates.
We recall (Section~\ref{subsec.mSR})
that the $m$-Stieltjes--Rogers polynomial $S^{(m)}_n(\balpha)$
is defined to be the generating polynomial
for $m$-Dyck paths from $(0,0)$ to ${((m+1)n,0)}$,
in~which each rise gets weight~1
and each $m$-fall from height~$i$ gets weight $\alpha_i$.
Now, more generally, for any integer $j \ge 0$,
we define the
\textbfit{modified $\bm{m}$-Stieltjes--Rogers polynomial of type $\bm{j}$},
denoted $S^{(m;j)}_{n}(\balpha)$,
to be the generating polynomial
for partial $m$-Dyck paths from $(0,0)$ to ${((m+1)n+j,j)}$,
in~which each rise gets weight~1
and each $m$-fall from height~$i$ gets weight $\alpha_i$.\footnote{
     In \cite{latpath_SRTR} we used the notation $S^{(m)}_{n|j}(\balpha)$
     for what we are now calling $S^{(m;j)}_{n}(\balpha)$.
}

For $j=0$, the modified $m$-Stieltjes--Rogers polynomials
are of course the usual $m$-Stieltjes--Rogers polynomials
$S^{(m;0)}_{n}(\balpha) = S^{(m)}_{n}(\balpha)$,
with ordinary generating function $f_0(t)$.
For general $j \ge 0$, the ordinary generating function
of the modified $m$-Stieltjes--Rogers polynomials of type $j$
was found in \cite[section~2.3]{latpath_SRTR}
(though this terminology was not used there).
Recall that $f_k(t)$ is the ordinary generating function
for $m$-Dyck paths at level $k$ with the weights $\balpha$.
We then have:

\begin{proposition}[Generating function of modified $m$-Stieltjes--Rogers polynomials]
   \label{prop.modified_SR.genfn}
For each integer $j \ge 0$, we have
\be
   \sum_{n=0}^\infty S^{(m;j)}_{n}(\balpha) \: t^n
   \;=\;
   f_0(t) \, f_1(t) \,\cdots\, f_j(t)
   \;.
 \label{eq.prop.modified_SR.genfn}
\ee
\end{proposition}

\noindent
The proof is a simple combinatorial argument
based on splitting the partial $m$-Dyck path
at its last return to level 0, then its last return to level 1, etc.\ 
\cite[sections~2.1 and 2.3]{latpath_SRTR}.

For $0 \le j \le m$,
the modified $m$-Stieltjes--Rogers polynomials $S^{(m;j)}_{n}(\balpha)$
also have interpretations in terms of the ordinary
$m$-Stieltjes--Rogers polynomials $S^{(m)}_{n}(\balpha)$.
Besides the trivial case $j=0$, the simplest case is $j=m$:
since an $m$-Dyck path of nonzero length must always end
with an $m$-fall from height $m$, we have
\be
   S^{(m;m)}_{n}(\balpha)  \;=\;  {S^{(m)}_{n+1}(\balpha) \over \alpha_m}
   \;.
 \label{eq.Smm}
\ee
Observe also that, multiplying both sides of \reff{eq.Smm} by $t^n$
and summing over $n \ge 0$ and then using \reff{eq.prop.modified_SR.genfn},
we obtain
\be
   f_0(t) \, f_1(t) \,\cdots\, f_m(t)
   \;=\;
   {f_0(t) \,-\, 1  \over \alpha_m t}
   \;,
\ee
which is precisely the fundamental functional equation \reff{eq.mSRfk.1}
at $k=0$.

Now consider \reff{eq.Smm} specialized to $\alpha_m = 0$:
it generates partial $m$-Dyck paths from $(0,0)$ to $((m+1)n+m,m)$
with the constraint that they must stay always at height $\ge 1$
(except of course for the starting point).
Removing the first step, which is necessarily a rise,
we have paths from $(1,1)$ to $((m+1)n+m,m)$
with the constraint that must stay always at height $\ge 1$;
or translating everything by $(-1,-1)$,
we have partial $m$-Dyck paths from $(0,0)$ to $((m+1)n+m-1,m-1)$.
It follows that
\be
   S^{(m;m-1)}_{n}(\balpha)
   \;=\;
   \left.  {S^{(m)}_{n+1}(\balpha) \over \alpha_m}
   \right| _{\alpha_m = 0 ,\: \alpha_i \to \alpha_{i-1}}
   \;.
 \label{eq.Smm-1}
\ee
Continuing in the same way, we have for $0 \le \ell \le m$,
\be
   S^{(m;m-\ell)}_{n}(\balpha)
   \;=\;
   \left.  {S^{(m)}_{n+1}(\balpha) \over \alpha_m}
   \right| _{\alpha_m = \alpha_{m+1} = \ldots = \alpha_{m+\ell-1} = 0 ,\;
             \alpha_i \to \alpha_{i-\ell}}
   \;.
 \label{eq.Smm-ell}
\ee
So all the modified $m$-Stieltjes--Rogers polynomials $S^{(m;j)}_{n}(\balpha)$
for $0 \le j \le m$ can be interpreted as $S^{(m)}_{n+1}(\balpha) / \alpha_m$
with some initial coefficients $\alpha_i$ set to zero
(and the variables renamed).
Multiplying both sides of \reff{eq.Smm-ell} by $t^n$
and summing over $n \ge 0$,
we obtain an alternate form for the generating function
of Proposition~\ref{prop.modified_SR.genfn} when $j \le m$:

\begin{proposition}[Alternate generating function of modified $m$-Stieltjes--Rogers polynomials]
   \label{prop.modified_SR.genfn.alternate}
For each integer $j$ satisfying $0 \le j \le m$, we have
\be
   \sum_{n=0}^\infty S^{(m;j)}_{n}(\balpha) \: t^n
   \;=\;
   \left.  {f_0(t) \,-\, 1  \over \alpha_m t}
   \right|_{\alpha_m = \alpha_{m+1} = \ldots = \alpha_{2m-j-1} = 0 ,\;
            \alpha_i \to \alpha_{i+j-m}}
   \;.
 \label{eq.prop.modified_SR.genfn.alternate}
\ee
\end{proposition}


\bigskip

We now generalize this construction by introducing,
for each integer $j \ge 0$, a triangular array.
For each triplet of integers $j,n,k \ge 0$,
let $S^{(m;j)}_{n,k}(\balpha)$ be the generating polynomial
for partial $m$-Dyck paths from $(0,0)$ to ${((m+1)n+j,(m+1)k+j)}$,
in~which each rise gets weight~1
and each $m$-fall from height~$i$ gets weight $\alpha_i$.
We call the $S^{(m;j)}_{n,k}$ the
\textbfit{generalized $\bm{m}$-Stieltjes--Rogers polynomials of type $\bm{j}$}.
Obviously $S^{(m;j)}_{n,k}$ is nonvanishing only for $0 \le k \le n$,
and $S^{(m;j)}_{n,n} = 1$.
We therefore have, for each integer $j \ge 0$,
an infinite unit-lower-triangular array
$\sfS^{(m;j)} = \big( S^{(m;j)}_{n,k}(\balpha) \big)_{n,k \ge 0}$.
The zeroth column of this array ($k=0$)
consists of the modified $m$-Stieltjes--Rogers polynomials of type $j$
defined earlier:
$S^{(m;j)}_{n,0}(\balpha) = S^{(m;j)}_{n}(\balpha)$,

Let us now take a closer look at the triangular arrays $\sfS^{(m;j)}$.
Note first that the arrays for $0 \le j \le m$ are the most fundamental,
in the sense that any array $\sfS^{(m;j)}$ for $j \ge m+1$
is a submatrix of the array $\sfS^{(m;j')}$ with $j' \eqdef j \bmod (m+1)$:
namely, if $j = j' + (m+1)\ell$,
then $S^{(m;j)}_{n,k} = S^{(m;j')}_{n+\ell,k+\ell}$.
In terms of the matrix $\Delta$ with 1 on the superdiagonal and 0 elsewhere,
we can write
\be
   \sfS^{(m;j' + (m+1)\ell)}
   \;=\;
   \Delta^\ell \sfS^{(m;j')} (\Delta^{\rm T})^\ell
   \;.
 \label{eq.S.Deltas}
\ee
So we shall concentrate, wherever necessary, on the cases $0 \le j \le m$.

The polynomials $S^{(m;j)}_{n,k}$ satisfy a simple pair of recurrences,
based on examining the two possibilities (rise or $m$-fall)
for the final step of a partial $m$-Dyck path:
\begin{eqnarray}
   S^{(m;j+1)}_{n,k}
   & = &
   S^{(m;j)}_{n,k} \:+\: \alpha_{(m+1)(k+1)+j} \, S^{(m;j)}_{n,k+1}
      \label{eq.genSR.1}  \\[2mm]
   S^{(m;j)}_{n+1,k}
   & = &
   S^{(m;j+m)}_{n,k-1} \:+\: \alpha_{(m+1)k+j+m} \, S^{(m;j+m)}_{n,k}
      \label{eq.genSR.2}
\end{eqnarray}
In each case, the first (resp.~second) term on the right-hand side
corresponds to the final step being a rise (resp.~an $m$-fall).
For the classical case $m=1$ and $j=0$,
these recurrences are \reff{eq.SnlSprimenl.recurrence}.

The production matrix for the triangle $\sfS^{(m;0)}$
was found in \cite[sections~7.1 and 8.2]{latpath_SRTR};
and the production matrix for the triangle $\sfS^{(m;m)}$
was in essence found in \cite[section~7.2]{latpath_SRTR}
(though this terminology was not used there).
Here we will generalize this analysis,
and find the production matrices for all the triangles $\sfS^{(m;j)}$
with $0 \le j \le m$.

We begin by defining some special matrices $M = (m_{ij})_{i,j \ge 0}$:
\begin{itemize}
   \item $L(s_1,s_2,\ldots)$ is the lower-bidiagonal matrix
       with 1 on the diagonal and $s_1,s_2,\ldots$ on the subdiagonal:
\be
   L(s_1,s_2,\ldots)
   \;=\;
   \begin{bmatrix}
      1  &     &     &     &    \\
      s_1 & 1  &     &     &    \\
          & s_2 & 1  &     &    \\
          &     & s_3 & 1  &    \\
          &     &     & \ddots & \ddots
   \end{bmatrix}
   \;.
 \label{def.L}
\ee
   \item $U^\star(s_1,s_2,\ldots)$ is the upper-bidiagonal matrix
       with 1 on the superdiagonal and $s_1,s_2,\ldots$ on the diagonal:
\be
   U^\star(s_1,s_2,\ldots)
   \;=\;
   \begin{bmatrix}
      s_1 & 1   &     &     &     &    \\
          & s_2 & 1   &     &     &    \\
          &     & s_3 & 1   &     &    \\
          &     &     & s_4 & 1   &    \\
          &     &     &     & \ddots & \ddots
   \end{bmatrix}
   \;.
 \label{def.Ustar}
\ee
\end{itemize}
Note that
\begin{subeqnarray}
    \Delta \, L(s_1,s_2,\ldots)   & = & U^\star(s_1,s_2,\ldots)  \\[2mm]
    L(s_1,s_2,\ldots) \, \Delta   & = & U^\star(0,s_1,s_2,\ldots)
 \label{eq.delta.shift}
\end{subeqnarray}

Let us now define, for each $r \ge 0$, the matrices
\be
   L_r  \;\eqdef\; L(\alpha_{m+r}, \alpha_{2m+r+1}, \alpha_{3m+r+2}, \ldots)
\ee
and
\be
   U_r  \;\eqdef\; \Delta L_r \;=\;
       U^\star(\alpha_{m+r}, \alpha_{2m+r+1}, \alpha_{3m+r+2}, \ldots)
   \;.
\ee
We then have:

\begin{proposition}[Production matrices for the generalized $m$-Stieltjes--Rogers polynomials of type $j$]
   \label{prop.Smj.prodmat}
For $0 \le j \le m$, the generalized $m$-Stieltjes--Rogers matrix of type $j$,
namely $\sfS^{(m;j)} = \big( S^{(m;j)}_{n,k}(\balpha) \big)_{n,k \ge 0}$,
has production matrix
$P^{(m;j)} \eqdef (\sfS^{(m;j)})^{-1} \Delta \sfS^{(m;j)}$
given by
\be
   P^{(m;j)}  \;=\;
   L_{j+1} L_{j+2} \,\cdots\, L_m \, U_0 \, L_1 L_2 \,\cdots\, L_j
   \;.
 \label{eq.prop.Smj.prodmat}
\ee
(Here $L_{j+1} L_{j+2} \,\cdots\, L_m$ is the identity matrix when $j=m$,
 and $L_1 L_2 \,\cdots L_j$ is the identity matrix when $j=0$.)
\end{proposition}

This result was also found independently,
and almost simultaneously, by H\'elder Lima \cite{Lima_preprint}.

\medskip

{\bf Remark.}
Any production matrix $P^{(m;j)}$ as in \reff{eq.prop.Smj.prodmat}
can be rewritten as $P^{(m;j')}$ for any $j' \in [j,m]$
with different parameters $\balpha$,
by using \reff{eq.delta.shift} to move the matrix $\Delta$ to the left:
$L(s_1,s_2,\ldots) \, \Delta = \Delta \, L(0,s_1,s_2,\ldots)$.
A little thought then shows that this result can be stated as
\be
   P^{(m;j)}(\alpha_m,\alpha_{m+1},\alpha_{m+2},\ldots)
   \;=\;
   P^{(m;j+1)}(0,\alpha_m,\alpha_{m+1},\alpha_{m+2},\ldots)
   \quad\hbox{for } 0 \le j \le m-1
   \;.
 \label{eq.j_to_j+1}
\ee
There is also a simple combinatorial argument for this identity:
as explained below,
$(P^{(m;j)})_{k,k'}$ is the sum over $(m+1)$-step walks in $\N$
going from height $(m+1)k+j$ to height $(m+1)k'+j$,
using steps $(1,1)$ and $(1,-m)$,
with weight $\alpha_i$ for each $m$-fall from height $i$;
and these quantities manifestly satisfy \reff{eq.j_to_j+1}.

However, the converse is not true:
it is not in general possible to {\em reduce}\/ $j$.
In particular, the modified $m$-Stieltjes--Rogers polynomials of type $j>0$
are a genuine generalization of those of type 0.
\myendremark

\medskip

The proof of Proposition~\ref{prop.Smj.prodmat}
is in fact a very minor modification
of the proof of \cite[Proposition~7.2]{latpath_SRTR},
but for clarity we include a few more details:

\proof
Let us first clarify what the matrices $L_r$ and $U_0$ do.
The $i$th row of the matrix $L_r$ corresponds to steps
starting at height $(m+1)i+r-1$:
more precisely,
the diagonal ($i \to i$) entries in the matrix $L_r$
give the weights for rises from height $(m+1)i+r-1$ to height $(m+1)i+r$,
while the subdiagonal ($i \to i-1$) entries in the matrix $L_r$
give the weights for $m$-falls from height $(m+1)i+r-1$
to height $(m+1)(i-1) + r$.
Similarly, the $i$th row of the matrix $U_0$ corresponds to steps
starting at height $(m+1)i+m$:
more precisely,
the superdiagonal ($i \to i+1$) entries in the matrix $U_0$
give the weights for rises from height $(m+1)i+m$ to height $(m+1)(i+1)$,
while the diagonal ($i \to i$) entries in the matrix $U_0$
give the weights for $m$-falls from height $(m+1)i+m$ to height $(m+1)i$.

The production matrix of $\sfS^{(m;j)}$
tells us how to get from the $n$th row of $\sfS^{(m;j)}$
to the $(n+1)$st row.
That is, to get the $(k,k')$ matrix element of $P^{(m;j)}$,
we need to enumerate the $(m+1)$-step walks in $\N$
going from height $(m+1)k+j$ to height $(m+1)k'+j$,
using steps $(1,1)$ and $(1,-m)$,
with weights $\balpha$;
we will show that the right-hand side of \reff{eq.prop.Smj.prodmat}
does the job.
The first step starts at height $(m+1)k+j$
and uses the matrix $L_{j+1}$:
it is either a rise, using the diagonal ($i \to i$) with weight~1,
or an $m$-fall, using the subdiagonal ($i \to i-1$)
with weight $\alpha_{(m+1)k+j}$.
We are now at height $(m+1)k_1 + j+1$,
where $k_1$ is either $k$ (if we made a rise)
or $k-1$ (if we made an $m$-fall).
The next step uses the matrix $L_{j+2}$ in an analogous way,
and we end up at height $(m+1)k_2 + j+2$,
where $k_2$ is either $k_1$ (if we made a rise)
or $k_1 -1$ (if we made an $m$-fall).
And so forth, through the matrix $L_m$.
After using the matrix $L_m$,
we are at height $(m+1)k_{m-j} + m$,
and we use the matrix $U_0$:
this step is either a rise, using the superdiagonal ($i \to i+1$)
with weight~1,
or else an $m$-fall, using the diagonal ($i \to i$)
with weight $\alpha_{(m+1)k_{m-j} + m}$.
We are now at height $(m+1)\widehat{k}$,
where $\widehat{k}$ equals $k+1$ minus the number of $m$-falls
that have occurred so far.
Now we use the matrix $L_1$:
it is either a rise, using the diagonal ($i \to i$) with weight~1,
or an $m$-fall, using the subdiagonal ($i \to i-1$)
with weight $\alpha_{(m+1)\widehat{k}}$.
We are now at height $(m+1)\widehat{k}_1 +1$,
where $\widehat{k}_1$ is either $\widehat{k}$ (if we made a rise)
or $\widehat{k}-1$ (if we made an $m$-fall).
The next step uses the matrix $L_2$ in an analogous way,
and so forth through $L_j$.
We end at height $(m+1)\widehat{k}_j + j$,
where $\widehat{k}_j$ equals $k+1$ minus the number of $m$-falls
that have occurred throughout the whole process.
But $k+1$ minus the number of $m$-falls is exactly the index
we have arrived at in the matrix product
$L_{j+1} L_{j+2} \,\cdots\, L_m \, U_0 \, L_1 L_2 \,\cdots\, L_j$,
because each rise (resp.~$m$-fall)
corresponds to a step $i \to i$ (resp.~$i \to i-1$),
except for the case of $U_0$,
for which a rise (resp.~$m$-fall)
corresponds to a step $i \to i+1$ (resp.~$i \to i$).
That is, $\widehat{k}_j = k'$, and the proof is complete.
\qed

It will be useful to write out explicitly the production matrices
for the case $m=2$, which is the one arising in this paper.
For $j=0$ we have \cite[eq.~(7.11)]{Sokal_multiple_OP}
\begin{subeqnarray}
   P_{n,n}    & = &  \alpha_{3n} \,+\, \alpha_{3n+1} \,+\, \alpha_{3n+2}
      \\[1mm]
   P_{n,n-1}  & = &  \alpha_{3n-2} \alpha_{3n} \,+\,
                       \alpha_{3n-1} \alpha_{3n} \,+\,
                       \alpha_{3n-1} \alpha_{3n+1}
      \\[1mm]
   P_{n,n-2}  & = &  \alpha_{3n-4} \alpha_{3n-2} \alpha_{3n}
 \label{eq.P20}
\end{subeqnarray}
{\em provided that we make the convention $\alpha_0 = \alpha_1 = 0$.}\/
(If $\alpha_{3n}$ and $\alpha_{3n+1}$ are given by polynomial expressions
in $n$ that do {\em not}\/ vanish when $n=0$,
then $P_{n,n}$ and $P_{n,n-1}$ are given by the corresponding
polynomial expressions {\em plus}\/ correction terms
proportional to $\delta_{n,0}$ and $\delta_{n,1}$.)
For $j=1$ we have
\begin{subeqnarray}
   P_{n,n}    & = & \alpha_{3n+1} \,+\, \alpha_{3n+2} \,+\, \alpha_{3n+3}
      \\[1mm]
   P_{n,n-1}  & = & \alpha_{3n-1} \alpha_{3n+1} \,+\,
                    \alpha_{3n} \alpha_{3n+1} \,+\,
                    \alpha_{3n} \alpha_{3n+2} 
      \\[1mm]
   P_{n,n-2}  & = &  \alpha_{3n-3} \alpha_{3n-1} \alpha_{3n+1}
 \label{eq.P21}
\end{subeqnarray}
{\em provided that we make the convention $\alpha_1 = 0$.}\/
For $j=2$ we have
\begin{subeqnarray}
   P_{n,n}    & = &  \alpha_{3n+2} \,+\, \alpha_{3n+3} \,+\, \alpha_{3n+4}
      \\[1mm]
   P_{n,n-1}  & = &  \alpha_{3n} \alpha_{3n+2} \,+\,
                     \alpha_{3n+1} \alpha_{3n+2} \,+\,
                     \alpha_{3n+1} \alpha_{3n+3} 
      \\[1mm]
   P_{n,n-2}  & = & \alpha_{3n-2} \alpha_{3n} \alpha_{3n+2} 
 \label{eq.P22}
\end{subeqnarray}
Note that these expressions arise from \reff{eq.P20} simply by
incrementing all the indices by $j$;
this fact is equivalent to \reff{eq.j_to_j+1}.

%


\bigskip

Let us now prove some theorems on total positivity.
We begin with the total positivity of the production matrices:

\begin{proposition}[Total positivity of the production matrices]
   \label{prop.Pmj.TP}
For $0 \le j \le m$, the production matrix
$P^{(m;j)} \eqdef (\sfS^{(m;j)})^{-1} \Delta \sfS^{(m;j)}$
is totally positive in the ring $\Z[\balpha]$
equipped with the coefficientwise order.
\end{proposition}

\proof
This is an immediate consquence of the factorization \reff{eq.prop.Smj.prodmat}
combined with Lemma~\ref{lemma.bidiagonal}.
\qed

Next we prove the total positivity of the triangular arrays $\sfS^{(m;j)}$.
As in \cite{latpath_SRTR},
we will give two proofs of each result:
a graphical proof, based on the Lindstr\"om--Gessel--Viennot lemma;
and an algebraic proof, based on production matrices.

\begin{theorem}[Total positivity for generalized $m$-Stieltjes--Rogers polynomials]
   \label{thm.Smj.TP}
\hfill\break\noindent
For each integer $j \ge 0$, the lower-triangular matrix $\sfS^{(m;j)}$
is totally positive in the ring $\Z[\balpha]$
equipped with the coefficientwise order.
\end{theorem}

\graphicalproof
We apply the Lindstr\"om--Gessel--Viennot lemma
(see \cite[section~9.4]{latpath_SRTR} for a summary)
to the directed graph $G_m = (V_m, E_m)$ with vertex set
\be
   V_m  \;=\;  \{ (x,y) \in \Z \times \N \colon\:
                  x=y \bmod (m+1)
               \}
 \label{def.Vm}
\ee
and edge set
\be
   E_m  \;=\;  \bigl\{ \bigl( (x_1,y_1),\, (x_2,y_2) \bigr) \in V_m \times V_m
                       \colon\:
                       x_2 - x_1 = 1 \hbox{ and } y_2 - y_1 \in \{1,-m\}
               \bigr\}
   \;.
 \label{def.Em}
\ee
It is easy to see that $G_m$ is planar and acyclic.

Now consider an $\ell \times \ell$ minor of $\sfS^{(m;j)}$:
that is, we choose sets of integers
$I = \{i_1,i_2,\ldots,i_\ell\}$ with $0 \le i_1 < i_2 < \ldots < i_\ell$
and $J = \{j_1,j_2,\ldots,j_\ell\}$ with $0 \le j_1 < j_2 < \ldots < j_\ell$,
and we consider the $\ell \times \ell$ submatrix
\be
   \sfS^{(m;j)}_{IJ}
   \;=\;
   \bigl( S^{(m;j)}_{i_r,j_s}(\balpha) \bigr)_{1 \le r,s \le \ell}
\ee
and the corresponding minor
\be
   \Delta_{IJ}(\sfS^{(m;j)})  \;=\;  \det \sfS^{(m;j)}_{IJ}  \;.
\ee
We can write the elements of the submatrix $\sfS^{(m;j)}_{IJ}$
as sums over walks in the directed graph $G_m$:
it is easy to see that $S^{(m;j)}_{i_r,j_s}(\balpha)$
is the sum over walks from $(-(m+1)i_r,0)$ to $(j,(m+1)j_s +j)$,
with a weight~1 on each rising directed edge
and a weight~$\alpha_i$ on each $m$-falling directed edge
starting at height~$i$.
Note that the sets of vertices
$I^\star = \{ (0,0)=i^\star_0 < (-(m+1),0)=i^\star_1 < (-2(m+1),0)=i^\star_2 < \ldots \}$
and $J^\star = \{ (j,j)=j^\star_0 < (j,m+1+j)=j^\star_1 < (j,2(m+1)+j)=j^\star_2 < \ldots \}$
lie on the boundary of $G_m \restrict ((-\infty,j] \times \N)$
in the order ``first $I^\star$ in reverse order, then $J^\star$ in order'';
therefore, by \cite[Lemma~9.18]{latpath_SRTR}
they form a fully nonpermutable pair.
Applying the Lindstr\"om--Gessel--Viennot lemma, we conclude that
the minor $\Delta_{IJ}(\sfS^{(m;j)})$ is the generating polynomial
for families of vertex-disjoint paths $P_1,\ldots,P_\ell$ in the digraph $G_m$,
where path $P_r$ starts at $(-(m+1)i_r,0)$ and ends at $(j,(m+1)j_r+j)$,
in which each rise gets weight 1
each $m$-fall from height~$i$ gets weight $\alpha_i$.
Since every such polynomial manifestly has nonnegative integer coefficients,
Theorem~\ref{thm.Smj.TP} is proven.
\qed

\algebraicproof
Since the array $\sfS^{(m;j)}$ for $j \ge m+1$
is a submatrix of the array $\sfS^{(m;j')}$ with $j' \eqdef j \bmod (m+1)$,
it suffices to prove the theorem for $0 \le j \le m$.
And the latter is an immediate consequence of the total positivity
of the production matrix $P^{(m;j)}$ (Proposition~\ref{prop.Pmj.TP})
combined with Theorem~\ref{thm.iteration.homo}.
\qed

{\bf Remark.}
Theorem~\ref{thm.Smj.TP} was proven for the special case $j=0$
in \cite[Theorem~9.10]{latpath_SRTR},
with both graphical and algebraic proofs.
The graphical proof given here for general $j \ge 0$
is the obvious generalization of the graphical proof given there.
The algebraic proof given here reduces to the one given there when $j=0$,
but for $1 \le j \le m$ it requires the new production matrices
\reff{eq.prop.Smj.prodmat},
and for $j \ge m+1$ it requires the obvious submatrix argument.
\myendremark

\medskip

The second result concerns the Hankel-total positivity of the
modified $m$-Stieltjes--Rogers polynomials $S^{(m;j)}_n(\balpha)$.
It was already proven in \cite[Theorem~9.12]{latpath_SRTR},
so we simply state the result here:

\begin{theorem}[Hankel-total positivity for modified $m$-Stieltjes--Rogers polynomials]
   \label{thm.Stype.minors.extended}
For $0 \le j \le m$,
the sequence
$\bS^{(m;j)} = ( S^{(m;j)}_{n}(\balpha) )_{n \ge 0}$
of modified $m$-Stieltjes--Rogers polynomials of type $j$
is a Hankel-totally positive sequence in the polynomial ring
$\Z[\balpha]$ equipped with the coefficientwise partial order.
\end{theorem}

The proof of Theorem~\ref{thm.Stype.minors.extended}
using the Lindstr\"om--Gessel--Viennot lemma
is similar to that of Theorem~\ref{thm.Smj.TP},
but using sink vertices
$J^{\star\star} = \{ (j,j)=j^{\star\star}_0 < (m+1+j,j)=j^{\star\star}_1 < (2(m+1)+j,j)=j^\star_2 < \ldots \}$
in place of $J^\star$.
For $j \le m$ (but {\em only}\/ then),
the source and sink vertices lie on the boundary of the graph $G_m$,
and the pair $(I^\star, J^{\star\star})$ is fully nonpermutable.
The algebraic proof of Theorem~\ref{thm.Stype.minors.extended}
uses the production matrix \reff{eq.prop.Smj.prodmat}
together with Theorem~\ref{thm.iteration2bis};
once again, it works only for $j \le m$,
since the ``submatrix argument'' does not apply to the zeroth column.
Indeed, the restriction of these proofs to $j \le m$ is no accident:
as explained in \cite[end of section~9.4]{latpath_SRTR},
the conclusion of Theorem~\ref{thm.Stype.minors.extended}
is {\em false}\/ for $j > m$;
in fact, it fails even for the $2 \times 2$ minors.

\bigskip

{\bf Final remark.}
For each fixed $m \ge 1$, the quantities $S^{(m;j)}_{n,k}$
depend on the three parameters $n,k,j$.
We have chosen here to fix $j$ and then assemble these quantities
into a lower-triangular matrix
$\sfS^{(m;j)} = \big( S^{(m;j)}_{n,k}(\balpha) \big)_{n,k \ge 0}$.
Lima \cite{Lima_23} has pointed out that one can also fix $k=0$
and assemble them into a (full) matrix
$\widehat{\sfS}^{(m)} = \big( S^{(m;j)}_{n,0}(\balpha) \big)_{n,j \ge 0}$.
He then shows \cite[Proposition~3.1]{Lima_23} the wonderful $LU$ factorization
\be
   \widehat{\sfS}^{(m)}  \;=\;  \sfS^{(m;0)} \, \Lambda^{(m)}
\ee
for an explicit upper-triangular matrix $\Lambda^{(m)}$.
Furthermore, he shows \cite[Proposition~3.2]{Lima_23}
that $\Lambda^{(m)}$ is coefficientwise totally positive.
Since $\sfS^{(m;0)}$ is also coefficientwise totally positive
by \cite[Theorem~9.10]{latpath_SRTR}
(i.e., the case $j=0$ of Theorem~\ref{thm.Smj.TP}),
this shows that $\widehat{\sfS}^{(m)}$ is coefficientwise totally positive
\cite[Theorem~3.3]{Lima_23}.
\myendremark

\subsection{Application to the univariate \hbox{Laguerre production matrix}}
        \label{subsec.mSR.univariate}

Let us now apply this theory to determine the conditions under which
the univariate Laguerre production matrix
\reff{eq.prop.prodmat.univariate.NEW.2}, namely
\begin{subeqnarray}
   p_{n,n+1}  & = &   1   \\[1mm]
   p_{n,n}    & = &   (2n+1+\alpha) \,+\, x   \\[1mm]
   p_{n,n-1}  & = &   n(n+\alpha) \,+\, 2nx \\[1mm]
   p_{n,n-2}  & = &   n(n-1) x
 \label{eq.app4.prodmat}
\end{subeqnarray}
can be written as the production matrix $P^{(m;j)}$
[cf.~\reff{eq.prop.Smj.prodmat}]
for the generalized $m$-Stieltjes--Rogers polynomials of type $j$
with $m=2$ and $j=0$, 1 or 2 and some $\balpha$.
To do this, we will use the formulae \reff{eq.P20}--\reff{eq.P22}
for $m=2$ and $j=0,1,2$,
remembering that $\alpha_0 = \alpha_1 \eqdef 0$.
In this analysis, we will interpret the Laguerre parameter $\alpha$
as a fixed real number $\ge -1$.
We will require that the Stieltjes parameters $\balpha$
be polynomials with real coefficients in the indeterminate $x$,
and we will most often require that these be polynomials
with {\em nonnegative}\/ real coefficients.

Since some of the computations are fairly lengthy,
it may be helpful to the reader to give now a summary of our results.
We will find that
\begin{itemize}
   \item $j=0$ works for $\alpha=-1$, and only that value;   \\[-7mm]
   \item $j=1$ works for $\alpha=-1$, 0, and only those values;   \\[-7mm]
   \item $j=2$ works for $\alpha=-1$, 0, 1, and only those values.
\end{itemize}
The allowed parameter sets $\balpha$ for each pair $(j,\alpha)$
are summarized in Table~\ref{table.j.univariate}.
We stress that the different columns in each row
correspond to different factorizations of the {\em same}\/ production matrix.

\begin{table}[t]
\begin{center}
\begin{tabular}{|c|c|c|c|}
\cline{2-4}
   \multicolumn{1}{c|}{\vphantom{$\displaystyle\sum$}} & $j=0$ & $j=1$ & $j=2$ \\[1mm]
\hline
$\alpha=-1$ & 
   \parbox{4cm}{\quad \\[-8mm]
                   \begin{center} Proposition~\ref{prop.Sfrac.j=0.solutions}
                   \end{center}
                \begin{eqnarray*}
                                  &       &     \\[-1mm]
                                  &       &     \\[-12mm]
                   \alpha_{3n-1}  & \!=\! & x   \\[-1mm]
                   \alpha_{3n}    & \!=\! & c_n \, n   \\[-1mm]
                   \alpha_{3n+1}  & \!=\! & (2 - c_n) \, n  \\[-7mm]
                \end{eqnarray*}}    &
   \parbox{4cm}{\quad \\[-8mm]
                   \begin{center} Proposition~\ref{prop.Sfrac.j=1.solutions.alpha=-1}
                   \end{center}
                \begin{eqnarray*}
                                  &       &     \\[-12mm]
                   \alpha_2       & \!=\! & 0 \\[-1mm]
                   \alpha_{3n}  & \!=\! & x   \\[-1mm]
                   \alpha_{3n+1}    & \!=\! & c_n \, n   \\[-1mm]
                   \alpha_{3n+2}  & \!=\! & (2 - c_n) \, n  \\[-7mm]
                \end{eqnarray*}}    &
   \parbox{4cm}{\quad \\[-8mm]
                   \begin{center} Proposition~\ref{prop.Sfrac.j=2.solutions.alpha=-1}
                   \end{center}
                \begin{eqnarray*}
                                  &       &     \\[-12mm]
                   \alpha_2 = \alpha_3  & \!=\! & 0 \\[-1mm]
                   \alpha_{3n+1}  & \!=\! & x   \\[-1mm]
                   \alpha_{3n+2}    & \!=\! & c_n \, n   \\[-1mm]
                   \alpha_{3n+3}  & \!=\! & (2 - c_n) \, n  \\[-7mm]
                \end{eqnarray*}} 
     \\
\hline
$\alpha=0$ &     &
   \parbox{4cm}{\quad \\[-8mm]
                   \begin{center} Proposition~\ref{prop.Sfrac.j=1.solutions.alpha=0}
                   \end{center}
                \begin{eqnarray*}
                                  &       &     \\[-12mm]
                                  &       &     \\[-1mm]
                   \alpha_{3n-1}  & \!=\! & x   \\[-1mm]
                   \alpha_{3n}    & \!=\! & n   \\[-1mm]
                   \alpha_{3n+1}  & \!=\! & n  \\[-7mm]
                \end{eqnarray*}}    &
   \parbox{4cm}{\quad \\[-8mm]
                   \begin{center} Proposition~\ref{prop.Sfrac.j=2.solutions.alpha=0}
                   \end{center}
                \begin{eqnarray*}
                                  &       &     \\[-12mm]
                   \alpha_2       & \!=\! & 0 \\[-1mm]
                   \alpha_{3n}  & \!=\! & x   \\[-1mm]
                   \alpha_{3n+1}    & \!=\! & n   \\[-1mm]
                   \alpha_{3n+2}  & \!=\! &  n  \\[-7mm]
                \end{eqnarray*}} 
     \\
\cline{1-1}\cline{3-4}
$\alpha=1$ & \multicolumn{1}{c}{\quad}     &      &
   \parbox{4cm}{\quad \\[-8mm]
                   \begin{center} Proposition~\ref{prop.Sfrac.j=2.solutions.alpha=1}
                   \end{center}
                \begin{eqnarray*}
                                  &       &     \\[-12mm]
                   \alpha_{3n-1}  & \!=\! & x   \\[-1mm]
                   \alpha_{3n}    & \!=\! & c_n \, n   \\[-1mm]
                   \alpha_{3n+1}  & \!=\! & (2 - c_n) \, n  \\[-7mm]
                \end{eqnarray*}} 
     \\
\cline{1-1}\cline{4-4}
\end{tabular}
\end{center}
\vspace*{-2mm}
\caption{
   Cases where the univariate Laguerre production matrix
   \reff{eq.prop.prodmat.univariate.NEW.2}/\reff{eq.app4.prodmat}
   can be written as the production matrix $P^{(m;j)}$
   for the generalized $m$-Stieltjes--Rogers polynomials of type $j$
   with $m=2$ and $j=0$, 1 or 2.
   Here $c_n$ is given by \reff{def.cn.j=0}, with $\kappa \in [0,1]$.
   Note that each passage to the right (at fixed $\alpha$) in this table
   is given by the transformation \reff{eq.j_to_j+1}.
}
\vspace*{5mm}
   \label{table.j.univariate}
\end{table}

\subsubsection[Case $j=0$]{Case $\bm{j=0}$}

\begin{proposition}[$j=0$ implies $\alpha=-1$]
  \label{prop.Sfrac.j=0}
The univariate Laguerre production matrix
\reff{eq.prop.prodmat.univariate.NEW.2}/\reff{eq.app4.prodmat}
can be written as the production matrix $P^{(2;0)}$,
where $\alpha_2,\alpha_4,\alpha_6$ are polynomials in $x$,
only if $\alpha = -1$.
\end{proposition}

\proof
From \reff{eq.app4.prodmat} and \reff{eq.P20} we have
\begin{subequations}
\begin{alignat}{2}
   p_{0,0}  &\;=\;  1 + \alpha + x     &&\;=\;  \alpha_2    \\[1mm]
   p_{2,0}  &\;=\;  2x                 &&\;=\;  \alpha_2 \alpha_4 \alpha_6
\end{alignat}
\end{subequations}
Since $\alpha_4$ and $\alpha_6$ are polynomials in $x$,
the polynomial $p_{0,0} = 1 + \alpha + x$ must divide $p_{2,0} = 2x$,
which occurs only if $\alpha = -1$.
\qed

\begin{proposition}[Solutions for $j=0$ and $\alpha=-1$]
  \label{prop.Sfrac.j=0.solutions}
The univariate Laguerre production matrix
\reff{eq.prop.prodmat.univariate.NEW.2}/\reff{eq.app4.prodmat}
with $\alpha=-1$
can be written as the production matrix $P^{(2;0)}$,
where the $\balpha$ are polynomials in $x$
with \emph{nonnegative} coefficients,
in the following ways (and \emph{only} the following ways):

Let $\kappa$ be any real number in the interval $[0,1]$, and set
\begin{subeqnarray}
   \alpha_{3n-1}  & = & x   \\
   \alpha_{3n}    & = & c_n \, n   \\
   \alpha_{3n+1}  & = & (2 - c_n) \, n
 \label{def.alphas.j=0}
\end{subeqnarray}
where
\be
   c_n  \;\eqdef\;  {(n-1) \,-\, (n-2)\kappa  \over  n \,-\, (n-1)\kappa}
   \quad\hbox{for } n \ge 1
   \;.
 \label{def.cn.j=0}
\ee
\end{proposition}

Please note that
$
   0 \:\le\: \kappa = c_1 \:\le\: c_2 \:\le\: c_3 \:\le\: \ldots \:\le\: 1
$
and $\lim\limits_{n \to \infty} c_n = 1$.

\medskip

\proof
We first prove the necessity, then the sufficiency.

\medskip

{\bf Necessity.}
For $j=0$ and $\alpha = -1$ we have
\begin{subequations}
 \label{eq.proof.prop.Sfrac.j=0.solutions.0}
\begin{alignat}{2}
   p_{0,0}  &\;=\; x            &&\;=\;  \alpha_2    \\[1mm]
   p_{n,n}  &\;=\; 2n+x         &&\;=\;
        \alpha_{3n} + \alpha_{3n+1} + \alpha_{3n+2}
        \quad\hbox{for } n \ge 1                     \\[1mm]
   p_{1,0}  &\;=\; 2x           &&\;=\;  \alpha_2 (\alpha_3 + \alpha_4) \\[1mm]
   p_{n,n-1} &\;=\; n(n-1) + 2nx   &&\;=\;
                       \alpha_{3n-2} \alpha_{3n} \,+\,
                       \alpha_{3n-1} (\alpha_{3n} + \alpha_{3n+1})
        \quad\hbox{for } n \ge 2                     \\[1mm]
   p_{n,n-2}  &\;=\;  n(n-1)x       &&\;=\;
                       \alpha_{3n-4} \alpha_{3n-2} \alpha_{3n}
\end{alignat}
\end{subequations}
From $p_{0,0}$ we get $\alpha_2 = x$;
from $p_{1,0}$ we then get $\alpha_3 + \alpha_4 = 2$;
and from $p_{2,0}$ we then get $\alpha_4 \alpha_6 = 2$.
We will now prove by induction that for all $n \ge 1$, we have
\begin{subeqnarray}
   \alpha_{3n-1}  & = &  x  \\[1mm]
   \alpha_{3n} + \alpha_{3n+1}  & = & 2n  \\[1mm]
   \alpha_{3n+1} \alpha_{3n+3}  & = & n(n+1)
 \label{eq.proof.prop.Sfrac.j=0.solutions}
\end{subeqnarray}
The base case $n=1$ has just been proven.
Using the inductive hypothesis ${\alpha_{3n} + \alpha_{3n+1}} = 2n$
and the formula for $p_{n,n}$, we deduce $\alpha_{3n+2} = x$.
Then using this together with the inductive hypothesis
$\alpha_{3n+1} \alpha_{3n+3} = n(n+1)$
and the formula for $p_{n+1,n}$, we deduce
$\alpha_{3n+3} + \alpha_{3n+4} = 2n+2$.
And finally, from $\alpha_{3n+2} = x$ and the formula for $p_{n+2,n}$,
we deduce $\alpha_{3n+4} \alpha_{3n+6} = (n+1)(n+2)$.
This completes the induction.

Writing $\alpha_{3n} = c_n \, n$,
(\ref{eq.proof.prop.Sfrac.j=0.solutions}b,c) say that
$\alpha_{3n+1} = (2 - c_n) \, n$ and $c_{n+1} = 1/(2 - c_n)$.
It is then easily proven by induction that the solution to this recurrence,
with initial condition $c_1 = \kappa$, is given by \reff{def.cn.j=0}.

By hypothesis we must have $c_1 \ge 0$, i.e.~$\kappa \ge 0$.
On the other hand, if $\kappa > 1$,
then there exists a positive integer $n$ satisfying
\be
   {n-1 \over n-2}  \;<\;  \kappa  \;\le\;  {n \over n-1}
   \;,
\ee
which makes either $c_n < 0$ or $c_n = \infty$,
both of which are forbidden.
So we must have $\kappa \in [0,1]$.

\medskip

{\bf Sufficiency.}
It is easy to verify that, for any $\kappa \in [0,1]$,
the parameters $\balpha$ defined by
\reff{def.alphas.j=0}/\reff{def.cn.j=0}
satisfy \reff{eq.proof.prop.Sfrac.j=0.solutions}
and thence \reff{eq.proof.prop.Sfrac.j=0.solutions.0}.
\qed

\medskip

{\bf Remarks.}
1.  The reasoning in this proof can be abstracted as a {\em general}\/ method
for writing a quadridiagonal unit-lower-Hessenberg matrix $P$
as the production matrix $P^{(2;0)}$ of a 2-S-fraction of type $j=0$
[cf.~\reff{eq.P20}]
whenever this is possible in the chosen ring, as follows:
For $n \ge 1$, define
\begin{subeqnarray}
   a_n  & = &  \alpha_{3n-1}  \\
   b_n  & = &  \alpha_{3n} + \alpha_{3n+1}  \\
   c_n  & = &  \alpha_{3n+1} \alpha_{3n+3}
\end{subeqnarray}
Then the equations for $p_{0,0}$, $p_{1,0}$ and $p_{2,0}$ give
\begin{subeqnarray}
   a_1  & = &  p_{0,0}  \\[1mm]
   b_1  & = &  {p_{1,0} \over p_{0,0}}  \\[1mm]
   c_1  & = &  {p_{2,0} \over p_{0,0}}
\end{subeqnarray}
if these divisions make sense in the chosen ring.
Then the equations for $p_{n,n}$, $p_{n+1,n}$ and $p_{n+2,n}$ give
successively
\begin{subeqnarray}
   a_{n+1}  & = &  p_{n,n} - b_n   \\[1mm]
   b_{n+1}  & = &  {p_{n+1,n} - c_n  \over a_{n+1}}
                   \;=\; {p_{n+1,n} - c_n  \over  p_{n,n} - b_n}   \\[1mm]
   c_{n+1}  & = &  {p_{n+2,n} \over a_{n+1}}
                   \;=\; {p_{n+2,n}  \over  p_{n,n} - b_n}
\end{subeqnarray}
if these divisions make sense.
Now define, for $n \ge 1$,
\begin{subeqnarray}
   d_n  & = &  \alpha_{3n}  \\
   e_n  & = &  \alpha_{3n+1}
\end{subeqnarray}
Then $d_1  = \alpha_3$ can be a freely chosen ring element, call it $\kappa$;
and we have the recurrences
\begin{subeqnarray}
   e_n      & = &  b_n - d_n  \\
   d_{n+1}  & = &  {c_n \over e_n}  \;=\;  {c_n \over b_n - d_n}
\end{subeqnarray}
if these divisions make sense.

2.  There are two extreme cases in Proposition~\ref{prop.Sfrac.j=0.solutions}.
On the one hand, $\kappa = 1$ implies $c_n = 1$ and hence
\be
   \balpha \;=\; (\alpha_i)_{i \ge 2}
      \;=\; x,1,1,x,2,2,x,3,3,x,4,4,\ldots
   \;.
 \label{eq.lah.alphas.1a}
\ee
This is the known 2-S-fraction for the Lah polynomials
\cite[Theorem~1.5 with $r=2$]{latpath_lah}.
On the other hand, $\kappa = 0$ implies $c_n = (n-1)/n$ and hence
\be
   \balpha \;=\; (\alpha_i)_{i \ge 2}
      \;=\; x,0,2,x,1,3,x,2,4,x,3,5,\ldots
   \;,
 \label{eq.lah.alphas.2a}
\ee
which is apparently new.
But there is in fact a one-parameter family, indexed by $\kappa \in [0,1]$,
that interpolates between these two extreme cases.
\myendremark

%

\subsubsection[Case $j=1$]{Case $\bm{j=1}$}

\begin{proposition}[$j=1$ implies $\alpha=-1$ or 0]
  \label{prop.Sfrac.j=1}
The univariate Laguerre production matrix
\reff{eq.prop.prodmat.univariate.NEW.2}/\reff{eq.app4.prodmat}
can be written as the production matrix $P^{(2;1)}$,
where $\alpha_2,\alpha_3,\ldots,\alpha_{13}$ are polynomials in $x$
with \emph{nonnegative} real coefficients,
only if $\alpha = -1$ or $0$.
\end{proposition}

\proof
From \reff{eq.app4.prodmat} and \reff{eq.P21} we have
\begin{subequations}
\begin{alignat}{2}
   p_{0,0}  &\;=\;  1 + \alpha + x     &&\;=\;  \alpha_2 + \alpha_3    \\[1mm]
   p_{1,1}  &\;=\;  3 + \alpha + x     &&\;=\;  \alpha_4 + \alpha_5 + \alpha_6
                                                                       \\[1mm]
   p_{2,2}  &\;=\;  5 + \alpha + x     &&\;=\;  \alpha_7 + \alpha_8 + \alpha_9  
                                                                       \\[1mm]
   p_{1,0}  &\;=\;  1 + \alpha + 2x    &&\;=\;
                  (\alpha_2 + \alpha_3) \alpha_4 + \alpha_3 \alpha_5   \\[1mm]
   p_{2,1}  &\;=\;  4 + 2\alpha + 4x   &&\;=\;
                  (\alpha_5 + \alpha_6) \alpha_7 + \alpha_6 \alpha_8   \\[1mm]
   p_{2,0}  &\;=\;  2x                 &&\;=\;  \alpha_3 \alpha_5 \alpha_7      
                                                                       \\[1mm]
   p_{3,1}  &\;=\;  6x                 &&\;=\;  \alpha_6 \alpha_8 \alpha_{10}   
                                                                       \\[1mm]
   p_{4,2}  &\;=\;  12x                &&\;=\;  \alpha_9 \alpha_{11} \alpha_{13}
\end{alignat}
\end{subequations}
Since $\alpha_2$ and $\alpha_3$ have nonnegative coefficients,
from $p_{0,0}$ we see that $\alpha \ge -1$.
If~$\alpha = -1$ we are done;
so we henceforth assume that $\alpha > -1$.

Comparing the equations for $p_{0,0}$ and $p_{1,0}$ gives
\be
   1 + \alpha + 2x  \;=\;  (1 + \alpha + x) \alpha_4 + \alpha_3 \alpha_5
   \;.
 \label{eq.j=1.p10}
\ee
Since $\alpha_3$ and $\alpha_5$ have nonnegative coefficients,
we can avoid an $x^2$ (or higher-order) term on the right-hand side only if
$\alpha_4$ is some real constant $c$;
moreover, by comparing the constant terms in \reff{eq.j=1.p10}
and using $\alpha > -1$, we see that  $0 \le c \le 1$.

Next we observe that $\alpha_3$ must divide $p_{2,0} = 2x$,
so either $\alpha_3 = d$ or $\alpha_3 = dx$,
for some real constant $d > 0$.

\medskip

{\bf Case $\bm{\alpha_3 = d}$.}
From $p_{0,0}$ with $\alpha_3 = d$ we get $d \le 1+\alpha$.
From \reff{eq.j=1.p10} with $\alpha_4 = c$ and $\alpha_3 = d$ we get
\be
   \alpha_5  \;=\;  {(1+\alpha)(1-c) \over d} \:+\: {2-c \over d} x
   \;.
\ee
Since $\alpha_5$ divides $p_{2,0} = 2x$,
we must either have $(1+\alpha)(1-c) = 0$ or $c=2$.
Since $\alpha > -1$ and $0 \le c \le 1$, we conclude that $c=1$.
Then $\alpha_5 = x/d$, and from $p_{1,1}$ we get
\be
   \alpha_6  \;=\; (2 + \alpha) \,+\, \Big( 1 \,-\, {1 \over d} \Big) x
   \;.
 \label{eq.starj1a}
\ee
Since $\alpha_6$ must divide $p_{3,1} = 6x$,
we must either have $\alpha = -2$ (which is not allowed) or $d=1$;
so $d=1$ and $\alpha_6 = 2+\alpha$.
{}From $p_{2,1}$ we then get
\be
   4 + 2\alpha + 4x
   \;=\;
   (2 + \alpha + x) \alpha_7 \:+\: (2+\alpha) \alpha_8
   \;.
\ee
{}From $p_{3,1}$ we see that $\alpha_8$ is either a constant
or a constant times $x$:
the first possibility gives
$\alpha_7 = 4$ and $\alpha_8 = -2$, which is forbidden;
the second possibility gives $\alpha_7 = 2$ and
\be
   \alpha_8  \;=\;  {2 \over 2+\alpha} \, x
   \;.
\ee
Then using $p_{2,2}$ we get
\be
   \alpha_9  \;=\;  (3+\alpha) \:+\: {\alpha \over 2+\alpha} \, x
   \;.
\ee
But $\alpha_9$ must divide $p_{4,2} = 12x$,
which is possible only if $\alpha=-3$ (which is not allowed) or $\alpha = 0$.
This completes the proof in the case $\alpha_3 = d$.

\medskip

{\bf Case $\bm{\alpha_3 = dx}$.}
From \reff{eq.j=1.p10} with $\alpha_4 = c$ and $\alpha_3 = dx$,
comparing constant terms and using $\alpha \neq -1$ we conclude that $c=1$;
then $\alpha_5 = 1/d$.  From $p_{1,1}$ we then get
\be
   \alpha_6  \;=\; \Big( 2+\alpha \,-\, {1 \over d} \Big) \:+\: x
   \;.
 \label{eq.starj1c}
\ee
Since $\alpha_6$ must divide $p_{3,1} = 6x$,
we must have $d = 1/(2+\alpha)$ and $\alpha_6 = x$.
{}From $p_{2,0}$ we get $\alpha_7 = 2$,
and then $p_{2,1}$ gives $\alpha_8 = 2$.
Then using $p_{2,2}$ we get $\alpha_9 = 1 + \alpha + x$.
But $\alpha_9$ must divide $p_{4,2} = 12x$, so $\alpha=-1$,
contrary to hypothesis.
We conclude that the case $\alpha_3 = dx$ cannot occur when $\alpha \neq -1$.
\qed

\begin{question}
\rm
Does this result hold if we omit the condition that the
coefficients of the polynomials $\balpha$ are nonnegative?
\end{question}

\medskip

\begin{proposition}[Solutions for $j=1$ and $\alpha=-1$]
  \label{prop.Sfrac.j=1.solutions.alpha=-1}
The univariate Laguerre production matrix
\reff{eq.prop.prodmat.univariate.NEW.2}/\reff{eq.app4.prodmat}
with $\alpha=-1$
can be written as the production matrix $P^{(2;1)}$
where the $\balpha$ are polynomials in $x$
with \emph{nonnegative} coefficients,
in the following ways (and \emph{only} the following ways):
\begin{subeqnarray}
   \alpha_2       & = & 0   \\
   \alpha_{3n}    & = & x   \\
   \alpha_{3n+1}  & = & c_n \, n   \\
   \alpha_{3n+2}  & = & (2 - c_n) \, n  \quad\hbox{for } n \ge 1
 \label{def.alphas.j=1.alpha=-1}
\end{subeqnarray}
where $c_n$ is given by \reff{def.cn.j=0}, with $\kappa \in [0,1]$.
\end{proposition}

Note that these coefficients $\balpha$ are precisely those
that are obtained from Proposition~\ref{prop.Sfrac.j=0.solutions}
by the transformation \reff{eq.j_to_j+1}.
So the sufficiency of these coefficients
is an immediate consequence of the sufficiency half
of Proposition~\ref{prop.Sfrac.j=0.solutions}.
The interesting part of the present result is therefore the necessity:
namely, these are the {\em only}\/ solutions for $j=1$ and $\alpha=-1$.

\proof
{\bf Necessity.}
In the proof of Proposition~\ref{prop.Sfrac.j=1}
we saw that either $\alpha_3 = d$ or $\alpha_3 = dx$,
for some real constant $d > 0$;
and that if $\alpha_3 = d$, then $d \le 1 + \alpha$.
But for $\alpha = -1$ this contradicts $d > 0$.
Therefore, for $\alpha = -1$ we must have $\alpha_3 = dx$.

For $j=1$ and $\alpha = -1$ we have
\begin{subequations}
\begin{alignat}{2}
   p_{0,0}  &\;=\; x            &&\;=\;  \alpha_2 + \alpha_3
                                                     \\[1mm]
   p_{n,n}  &\;=\; 2n+x         &&\;=\;
        \alpha_{3n+1} + \alpha_{3n+2} + \alpha_{3n+3}
        \quad\hbox{for } n \ge 1                     \\[1mm]
   p_{n,n-1} &\;=\; n(n-1) + 2nx   &&\;=\;
                       \alpha_{3n-1} \alpha_{3n+1} \,+\,
                       \alpha_{3n} \alpha_{3n+1} \,+\,
                       \alpha_{3n} \alpha_{3n+2}     \\[1mm]
   p_{n,n-2}  &\;=\;  n(n-1)x       &&\;=\;
                       \alpha_{3n-3} \alpha_{3n-1} \alpha_{3n+1}
\end{alignat}
\end{subequations}
From $\alpha_3 = dx$ and $p_{0,0}$ we have $\alpha_2 = (1-d)x$.
Hence $0 < d \le 1$.
From $p_{1,0}$ we deduce that $\alpha_4 + \alpha_5 d = 2$;
so $\alpha_4 = c$ for some real constant $c$,
and $\alpha_5 = (2-c)/d$.
Then from $p_{1,1}$ we get
\be
   \alpha_6
   \;=\;
   \biggl[ 2 \,-\, \Big( c + {2-c \over d} \Big) \biggr]  \:+\: x
   \;,
\ee
Since $\alpha_6$ must divide $p_{3,1} = 6x$,
we conclude that $c + (2-c)/d = 2$ and hence either $c=2$ or $d=1$.
Then either way we have $\alpha_6 = x$, and from $p_{2,1}$ we get
\begin{subeqnarray}
   2 + 4x
   & = &
   \Big( {2-c \over d} \,+\, x \Big) \, \alpha_7
     \:+\: \alpha_8 x
           \\[1mm]
   & = &
   {2-c \over d} \, \alpha_7   \:+\:  (\alpha_7 + \alpha_8) x
      \;.
\end{subeqnarray}
It follows that $\alpha_7 = 2d/(2-c)$ and $\alpha_7 + \alpha_8 = 1$;
therefore, $c=2$ is impossible, and we have $d=1$
and hence $\alpha_2 = 0$.

From here on the proof is identical to that of
Proposition~\ref{prop.Sfrac.j=0.solutions},
with each $\alpha_i$ replaced by $\alpha_{i+1}$.
%
%

\medskip

{\bf Sufficiency.}
It is easily verified that \reff{def.alphas.j=1.alpha=-1},
inserted into \reff{eq.P21},
yields \reff{eq.app4.prodmat} with $\alpha=-1$.
\qed

\begin{proposition}[Solutions for $j=1$ and $\alpha=0$]
  \label{prop.Sfrac.j=1.solutions.alpha=0}
The univariate Laguerre production matrix
\reff{eq.prop.prodmat.univariate.NEW.2}/\reff{eq.app4.prodmat}
with $\alpha=0$
can be written as the production matrix $P^{(2;1)}$
where the $\balpha$ are polynomials in $x$
with \emph{nonnegative} coefficients,
in the following way (and \emph{only} the following way):
\begin{subeqnarray}
   \alpha_{3n-1}  & = & x   \\
   \alpha_{3n}    & = & n   \\
   \alpha_{3n+1}  & = & n
 \label{def.alphas.j=1.alpha=0}
\end{subeqnarray}
\end{proposition}

\proof
{\bf Necessity.}
In the proof of Proposition~\ref{prop.Sfrac.j=1}
we saw that either $\alpha_3 = d$ or $\alpha_3 = dx$,
for some real constant $d > 0$;
and we saw that $\alpha_3 = dx$ leads to a contradiction
whenever $\alpha \neq -1$.
So for $\alpha = 0$ we must have $\alpha_3 = d$;
and in that case we saw that $d=1$ (hence $\alpha_3 = 1$)
and $\alpha_4 = 1$, $\alpha_5 = x$, $\alpha_6 = 2$,
$\alpha_7 = 2$, $\alpha_8 = x$, $\alpha_9 = 3$.
And from $p_{0,0}$ and $\alpha_3 = 1$ we get $\alpha_2 = x$.

For $j=1$ and $\alpha = 0$ we have
\begin{subequations}
\begin{alignat}{2}
   p_{0,0}  &\;=\; 1+x            &&\;=\;  \alpha_2 + \alpha_3
                                                     \\[1mm]
   p_{n,n}  &\;=\; 2n+1+x         &&\;=\;
        \alpha_{3n+1} + \alpha_{3n+2} + \alpha_{3n+3}
        \quad\hbox{for } n \ge 1                     \\[1mm]
   p_{n,n-1} &\;=\; n^2 + 2nx   &&\;=\;
                       \alpha_{3n-1} \alpha_{3n+1} \,+\,
                       \alpha_{3n} \alpha_{3n+1} \,+\,
                       \alpha_{3n} \alpha_{3n+2}     \\[1mm]
   p_{n,n-2}  &\;=\;  n(n-1)x       &&\;=\;
                       \alpha_{3n-3} \alpha_{3n-1} \alpha_{3n+1}
\end{alignat}
\end{subequations}
We now prove by induction that \reff{def.alphas.j=1.alpha=0} holds
for all $n \ge 1$.
The base case $n=1$ has already been proven (as has $n=2$).
Using the inductive hypothesis, from $p_{n,n-1}$ we get $\alpha_{3n+2} = x$.
Then from $p_{n,n}$ we get $\alpha_{3n+3} = n+1$.
And finally from $p_{n+1,n-1}$ we get $\alpha_{3n+4} = n+1$.

\medskip

{\bf Sufficiency.}
It is easily verified that \reff{def.alphas.j=1.alpha=0},
inserted into \reff{eq.P21},
yields \reff{eq.app4.prodmat} with $\alpha=0$.
\qed

Proposition~\ref{prop.Sfrac.j=1.solutions.alpha=0} shows that
the monic unsigned Laguerre polynomials with $\alpha=0$ are given by
the {\em modified}\/ 2-Stieltjes--Rogers polynomials of type $j=1$
with the coefficients
\be
   \balpha \;=\; (\alpha_i)_{i \ge 2}
      \;=\; x,1,1,x,2,2,x,3,3,x,4,4,\ldots
   \;.
 \label{eq.alphas.j=1.alpha=0}
\ee
Equivalently, the rook polynomials of an $n \times n$ chessboard,
which are the reversed monic unsigned Laguerre polynomials with $\alpha=0$,
are given by the modified 2-Stieltjes--Rogers polynomials of type $j=1$
with the coefficients
\be
   \balpha \;=\; (\alpha_i)_{i \ge 2}
      \;=\; 1,x,x,1,2x,2x,1,3x,3x,1,4x,4x,\ldots
   \;.
\ee
It was through this particular example that we first discovered
the importance of the modified $m$-Stieltjes--Rogers polynomials.

Please note that the coefficients \reff{eq.alphas.j=1.alpha=0}
are precisely \reff{eq.lah.alphas.1a},
i.e.~the $\kappa=1$ case of the coefficients in
Proposition~\ref{prop.Sfrac.j=0.solutions}.
Thus, these coefficients with $j=0$ give $\alpha = -1$,
while the same coefficients with $j=1$ give $\alpha = 0$;
and we will see in Proposition~\ref{prop.Sfrac.j=2.solutions.alpha=1}
that the same coefficients with $j=2$ give $\alpha = 1$.
But it is curious that for $j=1$, only $\kappa=1$ is allowed,
while for $j=0$ and $j=2$, we can take any $\kappa \in [0,1]$.

\subsubsection[Case $j=2$]{Case $\bm{j=2}$}

\begin{proposition}[$j=2$ implies $\alpha=-1$ or 0 or 1]
  \label{prop.Sfrac.j=2}
The univariate Laguerre production matrix
\reff{eq.prop.prodmat.univariate.NEW.2}/\reff{eq.app4.prodmat}
can be written as the production matrix $P^{(2;2)}$,
where $\alpha_2,\alpha_3,\ldots$ are polynomials in $x$
with \emph{nonnegative} coefficients,
only if $\alpha = -1$ or $0$ or $1$.
\end{proposition}

\proof
From \reff{eq.app4.prodmat} and \reff{eq.P22} we have
\begin{subequations}
\begin{alignat}{2}
   p_{0,0}  &\;=\;  1 + \alpha + x     &&\;=\;  \alpha_2 + \alpha_3 + \alpha_4
                                                                \\[1mm]
   p_{1,1}  &\;=\;  3 + \alpha + x     &&\;=\;  \alpha_5 + \alpha_6 + \alpha_7
                                                                       \\[1mm]
   p_{2,2}  &\;=\;  5 + \alpha + x     &&\;=\;
                  \alpha_8 + \alpha_9 + \alpha_{10}                    \\[1mm]
   p_{1,0}  &\;=\;  1 + \alpha + 2x    &&\;=\;
                  (\alpha_3 + \alpha_4) \alpha_5 + \alpha_4 \alpha_6   \\[1mm]
   p_{2,1}  &\;=\;  4 + 2\alpha + 4x   &&\;=\;
                  (\alpha_6 + \alpha_7) \alpha_8 + \alpha_7 \alpha_9   \\[1mm]
   p_{3,2}  &\;=\;  9 + 3\alpha + 6x   &&\;=\;
                  (\alpha_9 + \alpha_{10}) \alpha_{11} + \alpha_{10} \alpha_{12}
                                                                       \\[1mm]
   p_{2,0}  &\;=\;  2x                 &&\;=\;  \alpha_4 \alpha_6 \alpha_8      
                                                                       \\[1mm]
   p_{3,1}  &\;=\;  6x                 &&\;=\;  \alpha_7 \alpha_9 \alpha_{11}   
                                                                       \\[1mm]
   p_{4,2}  &\;=\;  12x                &&\;=\;  \alpha_{10} \alpha_{12} \alpha_{14}
\end{alignat}
\end{subequations}
Since $\alpha_2,\alpha_3,\alpha_4$ have nonnegative coefficients,
from $p_{0,0}$ we see that $\alpha \ge -1$.
If~$\alpha = -1$ we are done;
so we henceforth assume that $\alpha > -1$.

We also conclude from $p_{0,0}$ that
\be
   \alpha_3 + \alpha_4  \;=\;  c_1 \,+\, c_2 x
 \label{eq.starj2}
\ee
with $0 \le c_1 \le 1+\alpha$ and $0 \le c_2 \le 1$.
Furthermore, from $p_{2,0}$ we see that either
$\alpha_4 \alpha_6 = d$ or $\alpha_4 \alpha_6 = dx$, for some $d > 0$.
Finally, we observe that $\alpha_3 + \alpha_4$ cannot be identically zero,
because $p_{1,0} = 1 + \alpha + 2x$ has both constant and linear terms nonzero
(since $\alpha > -1$)
but $\alpha_4 \alpha_6$ has only a constant term {\em or}\/ a linear term.
So we must have either $c_1 > 0$ or $c_2 > 0$.

\bigskip

\underline{{\bf Case $\bm{\alpha_4 \alpha_6 = d}$.}}
From $p_{2,0}$ we get $\alpha_8 = (2/d) x$.
Furthermore, from $p_{1,0}$ we get
\begin{itemize}
   \item If $c_2 = 0$, then $c_1 > 0$ and
$\displaystyle \alpha_5 = {(1+\alpha-d) \,+\, 2x  \over  c_1}$.
   \item If $c_2 > 0$, then
$\displaystyle \alpha_5 \,=\, {2 \over c_2}$
and $d = (1+\alpha) - c_1 \alpha_5$.
\end{itemize}

{\bf When $\bm{c_2 = 0}$}, we substitute $\alpha_5$ into $p_{1,1}$ to get
\be
   \alpha_6 + \alpha_7  \;=\; 
   \Big( 3 + \alpha \,-\, {1+\alpha-d \over c_1} \Big)
   \:+\:
   \Big( 1 \,-\, {2 \over c_1} \Big) \, x
   \;.
\ee
From $p_{2,1}$ we have
\be
   4 + 2\alpha + 4x  \;=\;  (\alpha_6 + \alpha_7) \, {2 \over d} x
      \:+\: \alpha_7 \alpha_9
   \;.
 \label{eq.starj2a}
\ee
It follows that $\alpha_6 + \alpha_7$ cannot have a term $x$,
hence $c_1 =2$.
This implies $\alpha_6 + \alpha_7 = (5+\alpha+d)/2$.
Inserting this into \reff{eq.starj2a} we conclude that
\be
   \alpha_7 \alpha_9
   \;=\;
   (4 + 2\alpha)  \:+\: {3d - 5 - \alpha \over d} \, x
   \;.
\ee
Now $\alpha_7 \alpha_9$ must divide $p_{3,1} = 6x$;
but since $4 + 2 \alpha > 0$, we must have $d = (5 + \alpha)/3$.
Then $\alpha_7 \alpha_9 = 4 + 2\alpha$
and $\displaystyle \alpha_{11} = {3 \over 2+\alpha} \, x$.
Also, since $\alpha_8 = (2/d) x$,
we have $\alpha_8 = \displaystyle {6 \over 5+\alpha} \, x$.
Next, from $p_{2,2}$ we get
\be
   \alpha_9 + \alpha_{10}  \;=\;
   (5+\alpha) \:+\: \Big( 1 \,-\, {6 \over 5+\alpha} \Big) \, x
   \;.
 \label{eq.starj2b}
\ee
From $p_{3,2}$ we have
\be
   9 + 3\alpha + 6x  \;=\;  (\alpha_9 + \alpha_{10}) \, {3 \over 2+\alpha} x
      \:+\: \alpha_{10} \alpha_{12}
   \;.
 \label{eq.starj2c}
\ee
Combining \reff{eq.starj2b} and \reff{eq.starj2c},
to avoid an $x^2$ term we must have $6/(5+\alpha) = 1$, i.e.~$\alpha = 1$.

\medskip

{\bf When $\bm{c_2 > 0}$}, we substitute $\alpha_5$ into $p_{1,1}$ to get
\be
   \alpha_6 + \alpha_7  \;=\; 
   \Big( 3 + \alpha \,-\, {2 \over c_2} \Big)  \:+\:  x
   \;.
\ee
But then $\alpha_8 = (2/d)x$ produces an $x^2$ term
on the right-hand side of $p_{2,1}$, which is a contradiction.
So $c_2 > 0$ is impossible.

This completes the case $\alpha_4 \alpha_6 = d$.

\bigskip

\underline{{\bf Case $\bm{\alpha_4 \alpha_6 = dx}$.}}
{}From $p_{2,0}$ we get $\alpha_8 = 2/d$.
Furthermore, from $p_{1,0}$ we get
\begin{itemize}
   \item If $c_2 = 0$, then $c_1 > 0$ and
$\displaystyle \alpha_5 = {1 + \alpha \,+\, (2-d)x  \over c_1}$.
   \item If $c_2 > 0$, then $c_1 > 0$ and
$\displaystyle \alpha_5 \,=\, {1+\alpha \over c_1}
                        \,=\, {2-d \over c_2}$.
\end{itemize}

{\bf When $\bm{c_2 = 0}$}, we substitute $\alpha_5$ into $p_{1,1}$ to get
\be
   \alpha_6 + \alpha_7  \;=\; 
   \Big( 3+\alpha \,-\, {1+\alpha \over c_1} \Big)
   \:+\:
   \Big( 1 \,-\, {2-d \over c_1} \Big) \, x
   \;.
 \label{eq.starj2d}
\ee
From $p_{2,1}$ we have
\be
   4 + 2\alpha + 4x  \;=\;  (\alpha_6 + \alpha_7) \, {2 \over d}
      \:+\: \alpha_7 \alpha_9
   \;,
\ee
so
\be
   \alpha_7 \alpha_9
   \;=\;
   \biggl[ 4 + 2\alpha \,-\,
     \Big( 3+\alpha \,-\, {1+\alpha \over c_1} \Big) \, {2 \over d}
   \biggr]
   \:+\:
   \biggl[ 4 \,-\, \Big( 1 \,-\, {2-d \over c_1} \Big) \, {2 \over d}
   \biggr] \, x
   \;\eqdef\;
   A + Bx
   \;.
\ee
Now $\alpha_7 \alpha_9$ must divide $p_{3,1} = 6x$,
so either $A=0$ or $B=0$.
From $p_{2,2}$ we get
\be
   \alpha_9 + \alpha_{10}  \;=\;
   \Big( 5+\alpha - {2 \over d} \Big) \:+\: x
   \;,
 \label{eq.starj2e}
\ee
and from $p_{3,2}$ we have
\be
   9 + 3\alpha + 6x  \;=\;  (\alpha_9 + \alpha_{10}) \, \alpha_{11}
      \:+\: \alpha_{10} \alpha_{12}
   \;.
 \label{eq.starj2f}
\ee
Combining \reff{eq.starj2e} and \reff{eq.starj2f},
to avoid an $x^2$ term, $\alpha_{11}$ must be a constant;
therefore, from $p_{3,1}$ we conclude that $A=0$,
hence $\alpha_7 \alpha_9 = Bx$ and $\alpha_{11} = 6/B$.
It follows that
\be
   \alpha_{10} \alpha_{12}
   \;=\;
   \biggl[ 9 + 3\alpha \,-\, \Big( 5+\alpha - {2 \over d} \Big) {6 \over B}
   \biggr]
   \:+\:
   \Big( 6 \,-\, {6 \over B} \Big) \, x
   \;.
\ee
Now $\alpha_{10} \alpha_{12}$ must divide $p_{4,2} = 12x$,
so either
$\displaystyle B = {10 + 2\alpha - {4 \over d}  \over  3+\alpha}$ or $B=1$.
Combining this with $\alpha > -1$,
it can be shown that there are no solutions.\footnote{
\begin{samepage}
   {\bf If $\bm{B=1}$,} then $d = 2 (c_1 - 2)/(3c_1 - 2)$;
   substituting this into $A=0$ gives
   $$
     (2+2\alpha)  \:-\: (1+\alpha) c_1  \:+\: (5+\alpha) c_1^2
     \;=\; 0
     \;.
   $$
   But the discriminant of this quadratic is $- (1+\alpha) (39+7\alpha)$,
   which is $< 0$ whenever $\alpha > -1$.
  
   \medskip

   {\bf If $\displaystyle \bm{B = {10 + 2\alpha - {4 \over d}  \over  3+\alpha}}$,}
   then $d = [6+2\alpha - (1+\alpha)c_1] / [3+\alpha - (1+\alpha)c_1]$;
   substituting this into $A=0$ gives (assuming $\alpha \neq -1$)
   $$
      (3+\alpha) \:+\: 2 c_1  \:+\:  c_1^2
      \;=\;  0
     \;.
   $$
   Now the discriminant is $-4(2+\alpha)$,
   which is $< 0$ whenever $\alpha > -2$.
\end{samepage}
}

\bigskip

{\bf When $\bm{c_2 > 0}$}, from $p_{1,1}$ we get
\be
   \alpha_6 + \alpha_7  \;=\; 
   (3 + \alpha - \alpha_5)  \:+\: x
   \;.
 \label{eq.alpha6+7}
\ee
From $p_{2,1}$ we then have
\be
   4 + 2\alpha + 4x  \;=\;
   \Big[ (3 + \alpha - \alpha_5)  \:+\: x \Big] \, {2 \over d}
      \:+\: \alpha_7 \alpha_9
   \;,
\ee
so
\be
   \alpha_7 \alpha_9
   \;=\;
   \biggl[ 4 + 2\alpha \,-\, ( 3+\alpha - \alpha_5) \, {2 \over d} \biggr]
   \:+\:
   \Big( 4 \,-\, {2 \over d} \Big) \, x
   \;\eqdef\;
   A' + B' x
   \;.
 \label{eq.A82}
\ee
Now $\alpha_7 \alpha_9$ must divide $p_{3,1} = 6x$,
so either $A' =0$ or $B' =0$.
From $p_{2,2}$ we get as before
\be
   \alpha_9 + \alpha_{10}  \;=\;
   \Big( 5+\alpha - {2 \over d} \Big) \:+\: x
   \;,
 \label{eq.alpha9+10}
\ee
and from $p_{3,2}$ we get as before
\be
   9 + 3\alpha + 6x  \;=\;  (\alpha_9 + \alpha_{10}) \, \alpha_{11}
      \:+\: \alpha_{10} \alpha_{12}
   \;;
\ee
so to avoid an $x^2$ term, $\alpha_{11}$ must again be a constant.
Therefore, from $p_{3,1}$ we conclude that $A' =0$,
hence $\alpha_7 \alpha_9 = B' x$ and $\alpha_{11} = 6/B'$.
As before we get
\be
   \alpha_{10} \alpha_{12}
   \;=\;
   \biggl[ 9 + 3\alpha \,-\, \Big( 5+\alpha - {2 \over d} \Big) {6 \over B'}
   \biggr]
   \:+\:
   \Big( 6 \,-\, {6 \over B'} \Big) \, x
   \;;
\ee
and since $\alpha_{10} \alpha_{12}$ must divide $p_{4,2} = 12x$,
we have either
$B' = 1$ or $\displaystyle B' = {10 + 2\alpha - {4 \over d}  \over  3+\alpha}$.
In the first case we must have $d = 2/3$ to satisfy \reff{eq.A82},
and hence $\alpha_{10} \alpha_{12} = -3(1+\alpha) < 0$,
contrary to hypothesis.
In the second case we must have $d=1$ and hence $B'=2$,
$\alpha_4 \alpha_6 = x$, $\alpha_8 = 2$,
$\alpha_7 \alpha_9 = 2x$, $\alpha_{11} = 3$
and $\alpha_{10} \alpha_{12} = 3x$,
and then (from $A'=0$)
$\alpha_5 = 1$, $c_2 = 1$ and $c_1 = 1+\alpha$;
here $\alpha > -1$ is arbitrary.
But we now argue as follows:
%
%
%
Since $\alpha_6 + \alpha_7 = 2+\alpha+x$
and $\alpha_6$ divides $p_{2,0} = 2x$
and $\alpha_7$ divides $p_{3,1} = 6x$,
we must have either
\begin{itemize}
   \item $\alpha_6 = 2+\alpha$ and $\alpha_7 = x$,
     which implies $\alpha_9 = 2$ and $\alpha_{10} = 1+\alpha+x$;
     but $\alpha_{10}$ must divide $p_{4,2} = 12x$,
     which is a contradiction since $\alpha \neq -1$.
   \item $\alpha_6 = x$ and $\alpha_7 = 2+\alpha$,
     which implies $\alpha_9 = 2x/(2+\alpha)$ and
     $\alpha_{10} = 3+\alpha + {[1 - 2/(2+\alpha)]x}$;
     but $\alpha_{10}$ must divide $p_{4,2} = 12x$,
     which implies $\alpha=0$.
\end{itemize}
\qed

\begin{proposition}[Solutions for $j=2$ and $\alpha=-1$]
  \label{prop.Sfrac.j=2.solutions.alpha=-1}
The univariate Laguerre production matrix
\reff{eq.prop.prodmat.univariate.NEW.2}/\reff{eq.app4.prodmat}
with $\alpha=-1$
can be written as the production matrix $P^{(2;2)}$
where the $\balpha$ are polynomials in $x$
with \emph{nonnegative} coefficients,
in the following ways (and \emph{only} the following ways):
\begin{subeqnarray}
   \alpha_2 \,=\, \alpha_3  & = & 0   \\
   \alpha_{3n+1}  & = & x   \\
   \alpha_{3n+2}  & = & c_n \, n  \qquad\quad\;\,\hbox{for } n \ge 1  \\
   \alpha_{3n+3}  & = & (2 - c_n) \, n  \quad\hbox{for } n \ge 1
 \label{def.alphas.j=2.alpha=-1}
\end{subeqnarray}
where $c_n$ is given by \reff{def.cn.j=0}, with $\kappa \in [0,1]$.
\end{proposition}

Note that these coefficients $\balpha$ are precisely those
that are obtained from Proposition~\ref{prop.Sfrac.j=0.solutions}
by applying twice the transformation \reff{eq.j_to_j+1}.
So the sufficiency of these coefficients
is an immediate consequence of the sufficiency half
of Proposition~\ref{prop.Sfrac.j=0.solutions}.
The interesting part of the present result is therefore the necessity:
namely, these are the {\em only}\/ solutions for $j=2$ and $\alpha=-1$.

\proof
{\bf Necessity.}
For $j=2$ and $\alpha = -1$ we have
\begin{subequations}
\begin{alignat}{2}
   p_{n,n}  &\;=\; 2n+x         &&\;=\;
        \alpha_{3n+2} + \alpha_{3n+3} + \alpha_{3n+4}    \\[1mm]
   p_{n,n-1} &\;=\; n(n-1) + 2nx   &&\;=\;
                       (\alpha_{3n} + \alpha_{3n+1}) \alpha_{3n+2} \,+\,
                       \alpha_{3n+1} \alpha_{3n+3}     \\[1mm]
   p_{n,n-2}  &\;=\;  n(n-1)x       &&\;=\;
                       \alpha_{3n-2} \alpha_{3n} \alpha_{3n+2}
\end{alignat}
\end{subequations}
From $p_{0,0}$ we have $\alpha_2 + \alpha_3 + \alpha_4 = x$,
hence $\alpha_3 + \alpha_4 = c_2 x$ and $\alpha_4 = \widehat{c} x$
with $0 \le \widehat{c} \le c_2 \le 1$.
From $p_{2,0}$ we have $\alpha_4 \alpha_6 \alpha_8 = 2x$,
hence $\widehat{c} > 0$
and $\alpha_4 \alpha_6 = dx$ for some $d > 0$;
this gives $\alpha_6 = d/\widehat{c}$ and $\alpha_8 = 2/d$.

Next from $p_{1,0}$ we have $2x = c_2 x \alpha_5 + dx$,
hence $\alpha_5 = (2-d)/c_2$.
From $p_{1,1}$ we then get
\be
   \alpha_7  \;=\;  \Bigl( 2 \,-\, {2-d \over c_2} \,-\, {d \over \widehat{c}}
                    \Bigr)
                    \:+\:  x
             \;\eqdef\;  A \,+\, x
   \;.
\ee
Since $\alpha_7$ must divide $p_{3,1} = 6x$, we conclude that $A=0$, hence
\be
   0  \;\;=\;\;     2 \,-\, {2-d \over c_2} \,-\, {d \over \widehat{c}}
      \;\;\le\;\;   2 \,-\, {2-d \over c_2} \,-\, {d \over c_2}
      \;\;=\;\;     2 \,-\, {2 \over c_2}
      \;\;.
\ee
Since $c_2 \le 1$, this implies $c_2 = 1$,
and then equality here implies also $\widehat{c} = c_2 = 1$.
This implies $\alpha_2 = \alpha_3 = 0$.

From here on the proof is identical to that of
Proposition~\ref{prop.Sfrac.j=0.solutions},
with each $\alpha_i$ replaced by $\alpha_{i+2}$.
\qed

\begin{proposition}[Solutions for $j=2$ and $\alpha=0$]
  \label{prop.Sfrac.j=2.solutions.alpha=0}
The univariate Laguerre production matrix
\reff{eq.prop.prodmat.univariate.NEW.2}/\reff{eq.app4.prodmat}
with $\alpha=0$ can be written as the production matrix $P^{(2;2)}$
where the $\balpha$ are polynomials in $x$ with
{\em nonnegative} coefficients, in the following way
(and only the following way):
\begin{subeqnarray}
        \alpha_2 &=& 0\\
        \alpha_{3n} &=& x\\
        \alpha_{3n+1} &=& n\\
        \alpha_{3n+2} &=& n
\end{subeqnarray}
\end{proposition}

\proof
{\bf Necessity.}
In the proof of Proposition~\ref{prop.Sfrac.j=2} we saw that either
$\alpha_4\alpha_6 =d$ or $\alpha_4\alpha_6 =dx$,
for some real constant $d>0$;
and we saw that $\alpha_4\alpha_6 =d$ leads to
$\alpha = 1$ whenever $\alpha\neq -1$.
So, for $\alpha=0$ we must have $\alpha_4\alpha_6 =dx$.
Furthermore, we saw that
$\alpha_3 +\alpha_4 = c_1 + c_2 x$,
and whenever $\alpha\neq -1$ and $\alpha_4\alpha_6 =dx$,
to obtain $\alpha=0$, we need $c_2>0$.
In this case, we saw that we must have $d=1$
with $\alpha_6 = x$
and thus, $\alpha_4=1$.
We also obtained
$c_2 = 1$ and $c_1 =1$,
and thus, $\alpha_3=x$
and from $p_{0,0}$, $\alpha_2 = 0$.

From here on the proof is identical to that of
Proposition~\ref{prop.Sfrac.j=1.solutions.alpha=0},
with each $\alpha_i$ replaced by $\alpha_{i+1}$.

\medskip

{\bf Sufficiency.}  This is easily checked.
\qed

\begin{proposition}[Solutions for $j=2$ and $\alpha=1$]
  \label{prop.Sfrac.j=2.solutions.alpha=1}
The univariate Laguerre production matrix
\reff{eq.prop.prodmat.univariate.NEW.2}/\reff{eq.app4.prodmat}
with $\alpha=1$ can be written as the production matrix $P^{(2;2)}$
where the $\balpha$ are polynomials in $x$ with
{\em nonnegative} coefficients, in the following ways
(and only the following ways):
\begin{subeqnarray}
        \alpha_{3n-1} &=& x\\
        \alpha_{3n} &=& c_n n\\
        \alpha_{3n+1} &=& (2-c_{n}) n
 \label{def.alphas.j=2.alpha=1}
\end{subeqnarray}
where $c_n$ is given by \reff{def.cn.j=0}, with $\kappa \in [0,1]$.
\end{proposition}

\proof
{\bf Necessity.}
In the proof of Proposition~\ref{prop.Sfrac.j=2} we saw that either
$\alpha_4\alpha_6 =d$ or $\alpha_4\alpha_6 =dx$,
for some real constant $d>0$;
and we saw that $\alpha_4\alpha_6 =dx$ leads to
$\alpha = 0$ whenever $\alpha\neq -1$.
So, for $\alpha=1$ we must have $\alpha_4\alpha_6 =d$.
Furthermore, in this case we found that $c_2 = 0$ and $c_1 = 2$,
hence $\alpha_3+\alpha_4 = 2$.
We also found that $d = (5+\alpha)/3 =2$, hence $\alpha_4\alpha_6=2$.

For $j=2$ and $\alpha=1$ we have
\begin{subequations}
 \label{eq.proof.prop.Sfrac.j=2.solutions.alpha=1.0}
\begin{alignat}{2}
   p_{n,n}  &\;=\; 2n+2 + x         &&\;=\;
        \alpha_{3n+2} + \alpha_{3n+3} + \alpha_{3n+4}    \\[1mm]
   p_{n,n-1} &\;=\; n(n+1) + 2nx   &&\;=\;
                       (\alpha_{3n} + \alpha_{3n+1}) \alpha_{3n+2} \,+\,
                        \alpha_{3n+1} \alpha_{3n+3}   \\[1mm]
   p_{n,n-2}  &\;=\;  n(n-1)x       &&\;=\;
                       \alpha_{3n-2} \alpha_{3n} \alpha_{3n+2}
\end{alignat}
\end{subequations}
From $p_{0,0}$ and $\alpha_3+\alpha_4 = 2$, we see that $\alpha_2 = x$.
We will now prove by induction that for all $n \ge 1$, we have
\begin{subeqnarray}
   \alpha_{3n-1}  & = &  x  \\[1mm]
   \alpha_{3n} + \alpha_{3n+1}  & = & 2n  \\[1mm]
   \alpha_{3n+1} \alpha_{3n+3}  & = & n(n+1)
 \label{eq.proof.prop.Sfrac.j=2.solutions.alpha=1}
\end{subeqnarray}
The base case $n=1$ has just been proven.
Using the inductive hypotheses ${\alpha_{3n} + \alpha_{3n+1}} = 2n$
and $\alpha_{3n+1} \alpha_{3n+3} = n(n+1)$
and the formula for $p_{n,n-1}$, we deduce $\alpha_{3n+2} = x$.
Then using this together with the formula for $p_{n,n}$, we deduce
$\alpha_{3n+3} + \alpha_{3n+4} = 2n+2$.
Using $\alpha_{3n+1} \alpha_{3n+3} = n(n+1)$
and the formula for $p_{n+1,n-1}$, we deduce that $\alpha_{3n+5} = x$.
And using this together with $\alpha_{3n+3} + \alpha_{3n+4} = 2n+2$
and the formula for $p_{n+1,n}$,
we deduce that $\alpha_{3n+4} \alpha_{3n+6} = (n+1)(n+2)$.
This completes the induction.

From here on the proof is identical to that of 
Proposition~\ref{prop.Sfrac.j=0.solutions}.

\medskip

{\bf Sufficiency.}
It is easy to verify that, for any $\kappa \in [0,1]$,
the parameters $\balpha$ defined by
\reff{def.alphas.j=2.alpha=1}/\reff{def.cn.j=0}
satisfy \reff{eq.proof.prop.Sfrac.j=2.solutions.alpha=1}
and thence \reff{eq.proof.prop.Sfrac.j=2.solutions.alpha=1.0}.
\qed

\section{When is $\bm{B_\xi^{-1} P B_\xi}$ $\bm{(r,1)}$-banded?}
   \label{app.banded}

We say that a matrix $A = (a_{ij})_{i,j \ge 0}$
is \textbfit{$\bm{(r,s)}$-banded}
if $a_{ij} = 0$ whenever $j < i-r$ or $j > i+s$.
Otherwise put, the nonzero elements $a_{ij}$ can occur
only on the diagonal, in the first $r$ bands below the diagonal,
and in the first $s$ bands above the diagonal.
We say that a matrix is \textbfit{lower-Hessenberg}
if it is $(\infty,1)$-banded.

Let $P = (p_{ij})_{i,j \ge 0}$ be a lower-Hessenberg matrix
with entries in a commutative ring $R$.
We are interested in the matrix $B_\xi^{-1} P B_\xi$,
where $B_\xi$ is the $\xi$-binomial matrix \reff{def.Bx}
and $\xi$ is an indeterminate,
understood as a matrix with entries in the polynomial ring $R[\xi]$.
Obviously $B_\xi^{-1} P B_\xi$ is lower-Hessenberg;
we want to know under what conditions it is $(r,1)$-banded
for some integer $r \ge 0$.
An obvious necessary condition is that $P$ must itself be $(r,1)$-banded,
since $P$ can be recovered from $B_\xi^{-1} P B_\xi$ by specializing
$\xi$ to zero.  The full necessary and sufficient condition
can be written in terms of
the superdiagonal sequence $(p_{n,n+1})_{n \ge 0}$
and the $m$th subdiagonal sequences $(p_{n,n-m})_{n \ge m}$ for $0 \le m \le r$,
and goes as follows:

\begin{proposition}
   \label{prop.banded}
Fix an integer $r \ge 0$,
and let $P = (p_{ij})_{i,j \ge 0}$ be an $(r,1)$-banded matrix
with entries in a commutative ring $R$ containing the rationals.
Then the following are equivalent:
\begin{itemize}
   \item[(a)]  $B_\xi^{-1} P B_\xi$ is $(r,1)$-banded.
        [That is, the $m$th subdiagonal of $B_\xi^{-1} P B_\xi$ vanishes
         for all $m \ge r+1$.]
   \item[(b)]  The $(r+1)$st subdiagonal of $B_\xi^{-1} P B_\xi$ vanishes.
   \item[(c)]  The superdiagonal satisfies
\be
   p_{n,n+1}  \;=\; f_{-1}(n)
 \label{eq.prop.app.B.a}
\ee
where $f_{-1}(\,\cdot\,)$ is a polynomial of degree at most $r$
with coefficients in $R$;
and the $m$th subdiagonals for $0 \le m \le r$ satisfy
\be
   p_{n,n-m}  \;=\; n(n-1) \cdots (n-m+1) \: f_m(n)
   \;,
   \label{eq.prop.app.B.b}
\ee
where in each case $f_m(\,\cdot\,)$ is a polynomial of degree at most $r-m$
with coefficients in $R$.
\end{itemize}
\end{proposition}

\proof
(a)$\implies$(b) is trivial;
we will prove (b)$\implies$(c) and (c)$\implies$(a).

We begin by computing the matrix element $(B_\xi^{-1} P B_\xi)_{k+t,k}$
for $k,t \ge 0$:
\begin{subeqnarray}
   \!\!\!\!
   (B_\xi^{-1} P B_\xi)_{k+t,k}  \;=\;  (B_{-\xi} P B_\xi)_{k+t,k}
   & = &
   \sum\limits_{m=-1}^r \sum\limits_{j=k}^{k+t-m} \binom{k+t}{j+m} (-\xi)^{(k+t)-(j+m)} \:
        p_{j+m,j} \: \binom{j}{k} \xi^{j-k}
        \hspace*{-1cm} \nonumber \\ \slabel{eq.appBmatelem.a} \\[1mm]
   & = &
   \sum\limits_{m=-1}^r \xi^{t-m} \sum\limits_{j=k}^{k+t-m} (-1)^{(t-m)-(j-k)}
       \binom{k+t}{j+m} \binom{j}{k} \: p_{j+m,j}
   \;.
     \hspace*{-1cm} \nonumber \slabel{eq.appBmatelem.b} \\
        \label{eq.appBmatelem}
\end{subeqnarray}
(Here we let the binomial coefficient $\binom{N}{-1} = 0$ and $p_{-1,0} =0$.
This happens in \reff{eq.appBmatelem} when $j=k=0$ and $m=-1$.)

With $l = j-k$, we can rewrite equation~\reff{eq.appBmatelem} as
\begin{subeqnarray}
   [\xi^{t-m}] \, (B_\xi^{-1} P B_\xi)_{k+t,k}
   & = &
   \sum_{l=0}^{t-m} (-1)^{(t-m)-l} \binom{k+t}{k+l+m} \binom{k+l}{k} \:
        p_{k+l+m,k+l}
     \\[2mm]
   & = &
   {(k+t)! \over k! \, (t-m)!}
   \sum_{l=0}^{t-m} (-1)^{(t-m)-l} \binom{t-m}{l}
        {p_{k+l+m,k+l} \over (k+l+m)! / (k+l)!}
   \;,
      \nonumber \\
\end{subeqnarray}
valid when $k,t \ge 0$ and $-1 \le m \le t$
(remembering that $p_{-1,0} =0$).

So, for $m \ge -1$, let $g_m$ be the function defined on the domain $\N$
(with values in $R$)
by
\be
   g_m(k)  \;\eqdef\;  {p_{k+m,k} \over (k+m)! / k!}
   \label{eq.prop.appB.gs}
\ee
with the convention that $g_{-1}(0) = 0$.
Then we have
\begin{subeqnarray}
   [\xi^{t-m}] \, (B_\xi^{-1} P B_\xi)_{k+t,k}
   & = &
   {(k+t)! \over k! \, (t-m)!}
   \sum_{l=0}^{t-m} (-1)^{(t-m)-l} \binom{t-m}{l} \: g_m(k+l)
   \qquad
       \\[2mm]
   & = &
   {(k+t)! \over k! \, (t-m)!} \: (\Delta^{t-m} g_m)(k)
  \slabel{eq.forward}
\end{subeqnarray}
where $\Delta$ denotes the forward difference operator
$(\Delta f)(n) = f(n+1) - f(n)$;
this is valid when $k,t \ge 0$ and $-1 \le m \le t$.
It follows that
$[\xi^{t-m}] \, (B_\xi^{-1} P B_\xi)_{k+t,k}$ vanishes for all $k \ge 0$
if and only if $g_m$ is a polynomial of degree $< t-m$.

Applying this now under hypothesis~(b) for $t=r+1$ and $m \in [-1,r]$,
we conclude that $g_m$ is a polynomial of degree $\le r-m$.
When $m\ge 0$, we set $n = k+m$, which gives us \reff{eq.prop.app.B.b}.
When $m = -1$, equation~\reff{eq.prop.appB.gs} gives us
\be
p_{k-1,k} \;=\; {g_{-1}(k) \over k}
  \quad\hbox{for $k \ge 1$}
\ee
where $g_{-1}$ is a polynomial of degree $\le r+1$
that satisfies $g_{-1}(0) = 0$.
Setting $n=k-1$ gives us \reff{eq.prop.app.B.a}.
We have therefore proven that (b) implies (c).


On the other hand, if $g_m$ is a polynomial of degree $\le r-m$
for all $m \in [0,r]$,
and $p_{k-1,k} = f_{-1}(k)$ where $f_{-1}$ is a polynomial
of degree at most $r$,
it follows from \reff{eq.forward} that
$(B_\xi^{-1} P B_\xi)_{k+t,k} = 0$ for all $t \ge r+1$.
So (c) implies (a).
\qed

\section{Total positivity for a variant class of quadridiagonal matrices}
   \label{sec.prodmat.TP.multivariate.quadridiagonal.OLD}


In this appendix we prove total positivity
for a class of quadridiagonal matrices.
The result is very similar to Theorem~\ref{thm.prodmat.TP.bis.gen.new2},
but for a slightly different class of matrices;
and the proof is also very similar, though slightly more difficult.
We present this variant result because we think that
it may be useful in future work.
Indeed, in one current project of ours
\cite{Deb-Sokal_schett_FPSAC,Deb-Sokal_schett},
the total positivity of the production matrix can be proven using
a special case of Theorem~\ref{thm.quadmat2},
but not, as far as we can tell, using
Theorem~\ref{thm.prodmat.TP.bis.gen.new2}.

Consider the quadridiagonal lower-Hessenberg matrix $P$,
defined by
\begin{subeqnarray}
   P  & \eqdef &  L_1 L_2 U \, +\, L_1 D_1 \,+\, L_2 D_2
           \slabel{eq.thm.prodmat.TP.bis.gen.new3.a} \\[2mm] 
      & = &    L_1(L_2 U \, +\, D_1) \,+\, L_2 D_2
           \slabel{eq.thm.prodmat.TP.bis.gen.new3.b} \\[2mm] 
      & = &    L_2(L_1 U + D_2) \,+\, L_1 D_1 
           \slabel{eq.thm.prodmat.TP.bis.gen.new3.c}
   \label{eq.thm.prodmat.TP.bis.gen.new3}
\end{subeqnarray}
where
\begin{subeqnarray}
   L_1  & = &  \alpha I \,+\, xL  \\[1mm]
   L_2  & = &  \beta I  \,+\, yL
   \label{eq.thm.prodmat.TP.bis.gen.new3.L1L2}
\end{subeqnarray}
and
\begin{itemize}
	\item $L$ is the lower-bidiagonal matrix with 
		the sequence $a_0,a_1,\ldots$ on the diagonal, 
		the sequence $b_1, b_2, \ldots$ on the subdiagonal,
		and zeroes elsewhere;
	\item $U$ is the upper-bidiagonal matrix with 
		the sequence $c_1, c_2,\ldots$ on the superdiagonal,
		the sequence $d_0, d_1, \ldots$ on the diagonal,
		and zeroes elsewhere;
	\item $D_1$ is the diagonal matrix with entries $e_0, e_1,\ldots\,$;
	\item $D_2$ is the diagonal matrix with entries $f_0, f_1,\ldots\,$;
\end{itemize}
and $\alpha$, $\beta$, $x$, $y$,
$\bfa = (a_n)_{n\geq 0}$, $\bfb = (b_n)_{n\geq 1}$, 
$\bfc = (c_n)_{n\geq 1}$, $\bfd = (d_n)_{n\geq 0}$,
$\bfe = (e_n)_{n\geq 0}$, $\bfff = (f_n)_{n\ge 0}$
are all indeterminates.
Note that \reff{eq.thm.prodmat.TP.bis.gen.new3.b} =
\reff{eq.thm.prodmat.TP.bis.gen.new3.c}
because by construction $L_1$ and $L_2$ commute:
\be
   L_1  L_2  \;=\;  L_2 L_1
   \;.
\ee
This commutation is what makes the present situation
more restrictive than that of Theorem~\ref{thm.prodmat.TP.bis.gen.new2};
in compensation, the diagonal matrices $D_1$ and $D_2$
can both act on the same side (here the right).

The entries in the $k$th column of $P = (p_{n,k})_{n\geq 0}$
are given by
\begin{subeqnarray}
	p_{k-1,k}  & = &   (\alpha + x a_{k-1})(\beta + y a_{k-1}) c_k  \\[1mm]
	p_{k,k}    & = &   (\alpha + x a_k)(\beta + y a_k) d_k
	\,+\, (\alpha + x a_k )y b_k c_k 
	\,+\, x b_k (\beta+ y a_{k-1}) c_k \nonumber  \\
	& & \,+\, (\alpha + x a_k) e_k \,+\, (\beta+ y a_k) f_k
	\\[1mm]
	p_{k+1,k}  & = &   (\alpha+ x a_{k+1}) y b_{k+1} d_k 
	\,+\, x b_{k+1} (\beta + y a_k) d_k
	\,+\, x b_{k+1} y b_{k} c_{k} \nonumber \\
	& & 	\,+\, x b_{k+1} e_k \,+\, y b_{k+1} f_k
	\\[1mm]
	p_{k+2,k}  & = &   x y b_{k+2} b_{k+1} d_k   \\[1mm]
        p_{n,k}    & = &   0 \qquad\textrm{if $n<k-1$ or $n>k+2$}
\label{eq.thm.quadmat2}
\end{subeqnarray}
where by definition $b_0 = c_0 = 0$
and $a_n = b_n = c_n = d_n = e_n = f_n = 0$ whenever $n < 0$.
Our main result is:

\begin{theorem}[Total positivity of the generalized production matrix]
\label{thm.quadmat2}
The matrix $P$ defined by 
\reff{eq.thm.prodmat.TP.bis.gen.new3}/%
\reff{eq.thm.prodmat.TP.bis.gen.new3.L1L2}/%
\reff{eq.thm.quadmat2}
is totally positive, coefficientwise in the indeterminates
$\alpha,\beta, x,y, \bfa, \bfb, \bfc, \bfd, \bfe, \bfff$.
\end{theorem}

%

\medskip

The proof of Theorem~\ref{thm.quadmat2}
follows the same pattern as that of Theorem~\ref{thm.prodmat.TP.bis.gen.new2},
but works with columns rather than rows
(because the matrices $D_1$ and $D_2$ act on the right).
More precisely,
Theorem~\ref{thm.quadmat2} will be proven as follows:
Define the matrix $Q = (q_{n,k})_{n,k\geq 0} = (\bq_0,\bq_1,\ldots)$ by
\be
   Q  \;\eqdef\; P \big|_{\bfff = 0} \;=\; L_1(L_2 U + D_1)
   \;,
 \label{def.quadmat2.Q}
\ee
so that
\begin{subeqnarray}
        q_{k-1,k}  & = &   (\alpha + x a_{k-1})(\beta + y a_{k-1}) c_k  \\[1mm]
        q_{k,k}    & = &   (\alpha + x a_k)(\beta + y a_k) d_k
        \,+\, (\alpha + x a_k )y b_k c_k 
	 \,+\, x b_k (\beta+ y a_{k-1}) c_k \nonumber \\
	& & \,+\, (\alpha + x a_k) e_k 
        \\[1mm]
        q_{k+1,k}  & = &   (\alpha+ x a_{k+1}) y b_{k+1} d_k
        + x b_{k+1} (\beta + y a_k) d_k
        + x b_{k+1} y b_{k} c_{k} \nonumber \\
	& &	+ x b_{k+1} e_k
        \\[1mm]
        q_{k+2,k}  & = &   x y b_{k+2} b_{k+1} d_k   \\[1mm]
        q_{n,k}    & = &   0 \qquad\textrm{if $n<k-1$ or $n>k+2$}
  \label{eq.quadmat2.Q.entries}
\end{subeqnarray}
Then
\be
P \;=\; Q \,+\, L_2 D_2.
\label{eq.Q.P.difference.quadmat2}
\ee
We will begin by proving (Lemma~\ref{lemma.TP.quadmat2.Q})
that $Q$ is coefficientwise totally positive;
this proof 
uses the factorization $Q = L_1(L_2 U +D_1)$
together with the tridiagonal comparison theorem.
We omit the proof as it is almost identical
to that of Lemma~\ref{lemma.TP.Q.new2}.
It follows that for every integer $m \ge 0$,
the matrix $(\bq_0,\ldots,\bq_m)$ is totally positive.

The rest of the proof shows how to restore the terms in $P$ involving $\bfff$.
In terms of the column vectors 
$(\bp_k)_{k\geq 0}, (\bq_k)_{k\geq 0}, (\bll_k)_{k\geq 0}$
associated to the matrices $P, Q, L_2$, 
equation~\reff{eq.Q.P.difference.quadmat2} 
can be rewritten as 
\be
\bp_k \;=\; \bq_k \,+\, f_k  \bll_k,
\ee
where
\be
\bll_k
\;=\;
\begin{bmatrix}
	\bzero_{k \times 1}   \\[1mm]
        \beta+y a_k   \\[1mm]
        y b_{k+1}  \\[1mm]
        \bzero_{\infty \times 1}
\end{bmatrix}
\quad\textrm{for $k \ge 0$}
   \;.
\ee
We will show (Lemma~\ref{lemma.TP.quadmat2.QPfinal})
that for every pair of integers $0 \le k \le m+1$,
the matrix $(\bq_0,\ldots,\bq_{k-1},\bp_k,\ldots,\bp_m)$ is totally positive;
and we will do this, for each fixed $m \ge 0$,
by induction on $k = m+1,m,m-1,\ldots,0$.
The base case $k=m+1$ of this induction is thus Lemma~\ref{lemma.TP.quadmat2.Q},
and the final case $k=0$ is Theorem~\ref{thm.quadmat2}.
Similar to the proof of Lemma~\ref{lemma.TP.QPfinal.new2},
the proof of Lemma~\ref{lemma.TP.quadmat2.QPfinal} 
will involve the following steps:

\bigskip

   {\bf Lemma~\ref{lemma.TP.quadmat2.QR}:}
       The matrix $(\bq_0,\ldots,\bq_{k-1},\bll_k)$ is totally positive.

\bigskip

   \hangindent=1.5cm
   {\bf Lemma~\ref{lemma.TP.quadmat2.QPR.if}:}
       If the matrix $(\bp_{k+1},\ldots,\bp_m)$ is totally positive,
       then so is \hfill\break
       $(\bq_0,\ldots,\bq_{k-1},\bll_k,\bp_{k+1},\ldots,\bp_m)$.

\bigskip

   \hangindent=1.5cm
   {\bf The induction step (Lemma~\ref{lemma.TP.quadmat2.QP.inductive}):}
       If the matrix $(\bq_0,\ldots,\bq_{k},\bp_{k+1},\ldots,\bp_m)$
       \hfill\break
       is totally positive,
       then so is $(\bq_0,\ldots,\bq_{k-1},\bp_{k},\ldots,\bp_m)$.

\bigskip

\noindent
Putting this all together will prove Lemma~\ref{lemma.TP.quadmat2.QPfinal}
and hence Theorem~\ref{thm.quadmat2}.
Here Lemmas~\ref{lemma.TP.quadmat2.QR} and \ref{lemma.TP.quadmat2.QPR.if}
are closely analogous to
Lemmas~\ref{lemma.TP.QR.new2} and \ref{lemma.TP.QPR.if.new2},
but their proofs are somewhat more difficult.

\bigskip

We now begin the proof of Theorem~\ref{thm.quadmat2}.

\begin{lemma}[Total positivity of $Q$]
   \label{lemma.TP.quadmat2.Q}
The matrix $Q$ defined by 
\reff{def.quadmat2.Q}/\reff{eq.quadmat2.Q.entries}
is totally positive, coefficientwise in the indeterminates
$\alpha,\beta, x,y, \bfa, \bfb, \bfc, \bfd, \bfe$.

In particular, for every integer $m \ge 0$,
the matrix $(\bq_0,\ldots,\bq_m)$ is coefficientwise totally positive.
\end{lemma}

Lemma~\ref{lemma.TP.quadmat2.Q}
follows from the factorization $Q = L_1(L_2 U + D_1)$
by the same argument as in Lemma~\ref{lemma.TP.Q.new2}.

\medskip

\begin{lemma}
   \label{lemma.TP.quadmat2.QR}
For each integer $k \ge 0$,
the matrix $(\bq_0,\ldots,\bq_{k-1},\bll_k)$ is totally positive.
\end{lemma}

\proof
We prove this by induction on $k$.
The base case $k=0$ is trivial.
So we need to show the inductive step:
if $(\bq_0,\ldots,\bq_{k-2},\bll_{k-1})$ is totally positive,
then so is $(\bq_0,\ldots,\bq_{k-1},\bll_k)$.

Define $\widetilde{\bq}_{k-1} \eqdef \bq_{k-1} |_{d_{k-1} = 0}$, so that
$\widetilde{\bq}_{k-1} = (\widetilde{q}_{n,k-1})_{n\geq 0}^\intercal$
where
\begin{subeqnarray}
	\widetilde{q}_{k-2,k-1}  & = &   (\alpha + x a_{k-2})(\beta + y a_{k-2}) c_{k-1}  \\[1mm]
	\widetilde{q}_{k-1,k-1}    & = &   
	 (\alpha + x a_{k-1} )y b_{k-1} c_{k-1}
	\,+\, x b_{k-1} (\beta+ y a_{k-2}) c_{k-1} \nonumber\\
	& & \,+\, (\alpha + x a_{k-1}) e_{k-1}
        \\[1mm]
	\widetilde{q}_{k,k-1}  & = &   
        x b_{k} y b_{k-1} c_{k-1}
	\,+\, x b_{k} e_{k-1}
        \\[1mm]
        \widetilde{q}_{n,k-1}    & = &   0 \qquad\textrm{if $n<k-2$ or $n>k$}\;.
  \label{eq.quadmat2.Qtilde.entries}
\end{subeqnarray}
Then
\be
\bq_{k-1} 
\; = \;
\widetilde{\bq}_{k-1} 
\,+\, d_{k-1} (\alpha+x a_{k-1}) \bll_{k-1}
\,+\, d_{k-1} x b_{k} \bll_k
\;.
\ee
We will successively handle the second and third terms
on the right-hand side of this formula.

By Lemma~\ref{lemma.TP.quadmat2.Q} the matrix
$(\bq_0,\ldots,\bq_{k-1})$ is totally positive;
and specializing this matrix to $d_{k-1} = 0$ yields
$M_1 \eqdef (\bq_0,\ldots,\bq_{k-2},\widetilde{\bq}_{k-1})$,
which is therefore also totally positive.
And by the induction hypothesis, the matrix
$M_2 \eqdef (\bq_0,\ldots,\bq_{k-2},\bll_{k-1})$ is totally positive.
Applying Lemma~\ref{lemma.TP.1} (or rather its transpose)
to the matrices $M_1$ and $M_2$,
we conclude that the matrix
$M_3 
\eqdef 
(\bq_0,\ldots,\bq_{k-2}, \widetilde{\bq}_{k-1} 
+ d_{k-1} (\alpha+x a_{k-1})   \bll_{k-1})$
is totally positive.
Note that $M_3$ has $k+1$ nonzero rows (i.e.\ $0 \le n \le k$);
all subsequent rows are zero.
On the other hand, the column vector $\bll_k$ begins with $k$ zeroes.
Therefore, Lemma~\ref{lemma.TP.3} implies that
\be
   M_4
   \;\eqdef\;
   (M_3 \,|\, \bll_k)
   \;=\;
   \big( \bq_0,\ldots,\bq_{k-2}, \widetilde{\bq}_{k-1} 
   + d_{k-1} (\alpha+x a_{k-1})   \bll_{k-1}  ,
   \bll_k \big)
\ee
is totally positive.
Right-multiplying $M_4$ by the lower-bidiagonal matrix that has
1 on the diagonal, 
$d_{k-1} x b_{k}$
in position $(k,k-1)$ and zeroes elsewhere
--- in other words, adding $d_{k-1} x b_{k}$ times the last column of $M_4$
to its next-to-last column ---
we obtain the matrix
\be
   (\bq_0,\ldots,\bq_{k-2}, 
   \widetilde{\bq}_{k-1} 
   + d_{k-1} (\alpha+x a_{k-1}) \bll_{k-1}
   + d_{k-1} x b_{k} \bll_k, \bll_k)
   \;=\;
   (\bq_0,\ldots,\bq_{k-2}, \bq_{k-1}, \bll_k)
\ee
and prove its total positivity, completing the inductive step.
\qed

\begin{lemma}
   \label{lemma.TP.quadmat2.QPR.if}
Fix integers $0 \le k \le m$.
If the matrix $(\bp_{k+1},\ldots,\bp_m)$ is totally positive,
then so is the matrix
$(\bq_0,\ldots,\bq_{k-1},\bll_k,\bp_{k+1},\ldots,\bp_m)$.
\end{lemma}

\proof
The case $k=m$ is Lemma~\ref{lemma.TP.quadmat2.QR};  so assume that $k < m$.
Let $\bt_{k+1}$ be obtained from $\bp_{k+1}$
by first specializing $b_{k+1}=0$,
and followed by the substitution $f_{k+1} \to f_{k+1} + x b_{k+1} c_{k+1}$:
that is,
\be
  \bt_{k+1}
   \;=\;
   \left. \Bigl(
   \left. \bp_{k+1} \right|_{b_{k+1}=0}
   \Bigr) \right|_{f_{k+1} \to f_{k+1} + x b_{k+1} c_{k+1}}
   \;.
\ee
The entries in the column vector $\bt_{k+1}$ are 
$\bt_{k+1} = (t_{n,k+1})_{n\geq 0}^\intercal$
where
\begin{subeqnarray}
	t_{k,k+1}  & = &   (\alpha + x a_{k})(\beta + y a_{k}) c_{k+1}  \\[1mm]
	t_{k+1,k+1}    & = &   (\alpha + x a_{k+1})(\beta + y a_{k+1}) d_{k+1}
	\,+\, x y a_{k+1} b_{k+1} c_{k+1}
	\,+\, \beta x b_{k+1} c_{k+1} \nonumber \\
	& & \,+\, (\alpha + x a_{k+1}) e_{k+1} \,+\, (\beta+ y a_{k+1}) f_{k+1}
        \\[1mm]
	t_{k+2,k+1}  & = &   (\alpha+ x a_{k+2}) y b_{k+2} d_{k+1}
	\,+\, x b_{k+2} (\beta + y a_{k+1}) d_{k+1}
        \,+\, x b_{k+2} y b_{k+1} c_{k+1} \nonumber\\
	& & \,+\, x b_{k+2} e_{k+1} \,+\, y b_{k+2} f_{k+1}
        \\[1mm]
	t_{k+3,k+1}  & = &   x y b_{k+3} b_{k+2} d_{k+1}   \\[1mm]
        t_{n,k+1}    & = &   0 \qquad\textrm{if $n<k$ or $n>k+3$}
\end{subeqnarray}
Note that these substitutions would not affect $\bp_\ell$ for $\ell > k+1$.
Next, let $\widetilde{\bp}_{k+1}$ be identical to $\bt_{k+1}$
except that the entry $t_{k,k+1}$ is now made equal to 0,
i.e. $\widetilde{\bp}_{k+1} = (\widetilde{p}_{n,k+1})_{n\geq 0}^\intercal$
where
\begin{subeqnarray}
	\widetilde{p}_{k+1,k+1}    & = &   (\alpha + x a_{k+1})(\beta + y a_{k+1}) d_{k+1}
        \,+\, x y a_{k+1} b_{k+1} c_{k+1}
        \,+\, \beta x b_{k+1} c_{k+1} \nonumber \\
	& & \,+\, (\alpha + x a_{k+1}) e_{k+1} + (\beta+ y a_{k+1}) f_{k+1} 
	\\[1mm]
	\widetilde{p}_{k+2,k+1}  & = &   (\alpha+ x a_{k+2}) y b_{k+2} d_{k+1}
        \,+\, x b_{k+2} (\beta + y a_{k+1}) d_{k+1}
        \,+\, x b_{k+2} y b_{k+1} c_{k+1} \nonumber\\
	& & \,+\, x b_{k+2} e_{k+1} \,+\, y b_{k+2} f_{k+1}
        \\[1mm]
	\widetilde{p}_{k+3,k+1}  & = &   x y b_{k+3} b_{k+2} d_{k+1}   \\[1mm]
	\widetilde{p}_{n,k+1}    & = &   0 \qquad\textrm{if $n<k+1$ or $n>k+3$}
\end{subeqnarray}
Notice that 
\be
\bp_{k+1} 
\;=\;
\widetilde{\bp}_{k+1}
\,+\,
(\alpha +  x a_{k}) c_{k+1} \bll_k
\;.
\label{eq.bp.bptilde.quadmat2}
\ee

By hypothesis $M_1 \eqdef (\bp_{k+1},\ldots,\bp_m)$ is totally positive;
this implies, by substitution, that
$M_2 \eqdef (\bt_{k+1},\bp_{k+2},\ldots,\bp_m)$ is totally positive;
and finally, this implies that
$M_3 \eqdef (\widetilde{\bp}_{k+1},\bp_{k+2},\ldots,\bp_m)$ is totally positive,
because the nonzero rows of $M_3$ form a submatrix of $M_2$.

Now observe that the matrix
$S \eqdef (\bq_0,\ldots,\bq_{k-1},\bll_k,
               \widetilde{\bp}_{k+1},\bp_{k+2},\ldots,\bp_m)$
consists of two blocks overlapping in a single row:
\be
   S
   \;=\;
   (\bq_0,\ldots,\bq_{k-1},\bll_k \,|\,
               \widetilde{\bp}_{k+1},\bp_{k+2},\ldots,\bp_m)
   \;=\;
   \left[
       \begin{array}{cc|c}
           *   &   \bzero_{k \times 1}  &  \bzero_{k \times (m-k)} \\[1mm]
           *   &   \beta+y a_k          &  \bzero_{1 \times (m-k)} \\[1mm]
           \hline
           *   &   y b_{k+1}            &  *                   \\[1mm]
           \hline
           \bzero_{\infty \times k}  &  \bzero_{\infty \times 1}  & *
       \end{array}
   \right]
\ee
where the asterisks stand for blocks of unspecified entries
(which may be zero or nonzero).
By Lemma~\ref{lemma.TP.quadmat2.QR}, the matrix $(\bq_0,\ldots,\bq_{k-1},\bll_k)$
is totally positive;
and we have just shown that the matrix
$(\widetilde{\bp}_{k+1},\bp_{k+2},\ldots,\bp_m)$ is totally positive.
So Lemma~\ref{lemma.TP.3} implies that the matrix $S$ is totally positive.

On the other hand, 
using
equation~\reff{eq.bp.bptilde.quadmat2}
we can conclude that
the matrix\\ 
$(\bq_0,\ldots,\bq_{k-1},\bll_k,\bp_{k+1},\ldots,\bp_m)$
can be obtained from $S$ by right-multiplying it
by the upper-bidiagonal matrix that has
1 on the diagonal, 
$(\alpha +  x a_{k}) c_{k+1}$ in position $(k,k+1)$ and zeroes elsewhere.
This proves that the matrix
$(\bq_0,\ldots,\bq_{k-1},\bll_k,\bp_{k+1},\ldots,\bp_m)$ is totally positive.
\qed

The next two lemmas are the following:

\begin{lemma}
   \label{lemma.TP.quadmat2.QP.inductive}
Fix integers $0 \le k \le m$.
If the matrix $(\bq_0,\ldots,\bq_{k},\bp_{k+1},\ldots,\bp_m)$
is totally positive,
then so is $(\bq_0,\ldots,\bq_{k-1},\bp_{k},\ldots,\bp_m)$.
\end{lemma}


\begin{lemma}
   \label{lemma.TP.quadmat2.QPfinal}
For every pair of integers $0 \le k \le m+1$,
the matrix $(\bq_0,\ldots,\bq_{k-1},\bp_k,\ldots,\bp_m)$ is totally positive.
\end{lemma}


We omit the proofs of Lemmas~\ref{lemma.TP.quadmat2.QP.inductive}
and \ref{lemma.TP.quadmat2.QPfinal}
as they are analogous to the proofs of
Lemmas~\ref{lemma.TP.QP.inductive.new2}
and~\ref{lemma.TP.QPfinal.new2}, respectively,
but using columns instead of rows.

This completes the proof of Theorem~\ref{thm.quadmat2}.

\bigskip

We conclude by posing the following open problem:

\begin{problem}
Find a combinatorial interpretation for the output matrix $A = \scro(P)$
generated by the production matrix
\reff{eq.thm.prodmat.TP.bis.gen.new3}/\reff{eq.thm.quadmat2},
or by interesting specializations thereof.
\end{problem}

\end{document}